\documentclass[a4paper]{article}
\usepackage[english]{babel}					
\usepackage[utf8]{inputenc}					
\usepackage{fullpage}						
\usepackage[T1]{fontenc}					   
\usepackage{amsmath}						   
\usepackage{amsthm}							
\usepackage{amsfonts}						
\usepackage{amssymb}							
\usepackage{bbm}								
\usepackage{stmaryrd}						
\usepackage{tikz}							
\usepackage{graphicx}						
\usepackage{array}							
\usepackage{multirow}						
\usepackage[colorlinks=true]{hyperref}
\usepackage{pdfsync}
\usepackage{dsfont}

\newtheorem{thm}{Theorem}[section]

\newtheorem{clm}[thm]{Claim}
\def\R{\mathbb{R}}

\newcommand{\footremember}[2]{%
    \footnote{#2}
    \newcounter{#1}
    \setcounter{#1}{\value{footnote}}%
}
\newcommand{\footrecall}[1]{%
    \footnotemark[\value{#1}]%
} 

\def\one#1{\mathds{1}_{#1}}
\def\BV{{\mathrm{BV}}}
\def\dist{\operatorname{dist}}
\def\ds{\displaystyle}
\def\H{{\mathcal{H}}}
\def\haus{\H^{d-1}}

\title{Multiphase mean curvature flows with high mobility contrasts: a phase-field approach, with applications to nanowires}
\author{Elie Bretin\footremember{insa}{Univ Lyon, INSA de Lyon, CNRS UMR 5208, Institut Camille Jordan, 20 avenue Albert Einstein, F-69621 Villeurbanne Cedex, France.} \and Alexandre Danescu\footremember{inl}{Univ Lyon, Institut de Nanotechnologie de Lyon, CNRS UMR 5270, Ecole Centrale de Lyon, 36 Av. Guy de Collongue, F-69134 Ecully, France.} \and Jos\'e Penuelas\footrecall{inl} \and  Simon Masnou\footremember{ul1}{Univ Lyon, Universit\'e Claude Bernard Lyon 1, CNRS UMR 5208, Institut Camille Jordan, 43 blvd. du 11 novembre 1918, F-69622 Villeurbanne cedex, France.}}

\date{September 15, 2017}

\providecommand{\keywords}[1]{\textbf{\textit{Keywords: }} #1}

\begin{document}
\maketitle

\begin{abstract}
The structure of many multiphase systems is governed by an energy that penalizes the area of interfaces between phases weighted by surface tension coefficients.
However, interface evolution laws depend also on interface  mobility  coefficients.
Having in mind some applications where highly contrasted or even degenerate mobilities are involved, for which classical phase field models are inapplicable, we propose a new effective phase field approach to approximate multiphase mean curvature flows with mobilities.
The key aspect of our model is to incorporate the mobilities not in the phase field energy (which is conventionally the case) but in the metric which determines the gradient flow. We show the consistency of such approach by a formal analysis of the sharp interface limit. We also propose an efficient numerical scheme which allows us to illustrate the advantages of the model on various examples, as the wetting of droplets on solid surfaces  or the simulation of nanowires growth generated by the so-called vapor-liquid-solid method.

\end{abstract}
\section{Introduction }
Many physical systems involve a collection of interfaces whose positions and shapes are constrained so as to minimize their total area, e.g., soap foams, immiscible fluids, polycrystalline materials, etc. A typical expression of the total area is
$$\sum_{i\not=j}\sigma_{ij}\mbox{Area}(\Gamma_{ij}),$$
where $\{\Gamma_{ij}\}_{i\not=j}$ denote the interfaces, and $\{\sigma_{ij}\}_{i,j}$ are the so-called surface tensions. Depending on the context, the area energy may be either isotropic (like for soap foams) or anisotropic (as for polycrystalline materials). Starting from a given interfacial configuration, the energy gradient flow toward a minimizingconfiguration follows two rules~\cite{Herring1999} :
\begin{itemize}
 \item Every interface evolves by mean curvature flow, i.e., at every point $x \in \Gamma_{ij}$  which is not a junction point between three of more interfaces, the normal velocity $V_{ij}$ of the interface $\Gamma_{ij}$ is proportional to its mean curvature, further denoted as $H_{ij}$:
\begin{equation*}
   \frac{1}{m_{ij}} V_{ij} = \sigma_{ij} H_{ij},
\end{equation*}
where $m_{ij}$ is called the {\em mobility} of interface $\Gamma_{ij}$ -- hereafter, without explicit mention, no summation over repeated subscripts $i,j$ (Einstein rule) is assumed. Such proportionality between the normal velocity and the mean curvature is a straightforward consequence of (a localized version of) the Clausius-Duhem inequality as derived in sharp interface theory, see~\cite{gurtin1999sharp,fried1996phase,jiang2012phase,cahn1960theory}.
\item The Herring's angle condition holds at every triple junction, {\em i.e.} if $x$ is a triple junction between phases $i$, $j$, and $k$ then
$$ \sigma_{ij} n_{ij} +  \sigma_{jk} n_{jk} +  \sigma_{ki} n_{ki} = 0,$$
where $n_{ij}$ denotes the unit normal at $x$ to $\Gamma_{ij}$, pointing from $\Omega_i$ to $\Omega_j$. Here, we denote by $\{\Omega_i\}_i$ the collection of relatively closed sets which form the phase partition, so that $\Gamma_{ij}=\Omega_i\cap\Omega_j$. The physical meaning of Herring's condition is the equilibrium of the triple point due to force balance. This condition holds if, and only if, the surface tensions satisfy the triangle inequality $\sigma_{ik}\leq \sigma_{ij} + \sigma_{jk}$, $\forall i,j,k$.
\end{itemize}

It is important to notice that surface tensions appear both in the geometric energy and in the gradient flow, whereas mobilities play a role only in the gradient flow. This observation is a key to the approach that we shall propose in this paper.

There is a vast literature on numerical methods for the approximation of mean curvature flows. Methods can be roughly classified into the following four categories (some of them are exhaustively reviewed and compared in~\cite{review_interface}):
\begin{itemize}
\item[(i)] parametric methods~\cite{Deckelnick2005,Barrett2008_2};
\item[(ii)] level set methods \cite{OsherSethian,Osherbook1,Osherbook2,Evans_spruck,chen_giga_goto};
\item[(iii)] convolution/thresholding type algorithms~\cite{BenceMerrimanOsher,Ishhi_pires_souganidis,Ruuth_efficient,Esedoglu_Otto};
\item[(iv)] phase field approaches \cite{Modica1977,Chen1992,Garcke199887,Garcke_amulti,Bretin_Masnou_multiphase}.
\end{itemize}

In this paper, we focus on phase field methods, which can be used in many various physical situations, and offer both nice theoretical properties and numerical efficiency. A phase field model able to approximate a multiphase mean curvature flow with mobilities, as above, has been proposed already in~\cite{Garcke_amulti}. However, as we shall see later, this model is not well suited from a numerical viewpoint for handling highly contrasted or degenerate mobilities. We will propose in this paper a new method which is more adapted, and both efficient and numerically accurate. Our motivation for this work comes from a real example where high contrast of mobilities can be observed: the growth of nanowires by the so-called {\it vapor-liquid-solid} method. We will discuss a model for simulating a simplified approximation of such growth.

\subsection{Nanowires grown by the VLS method}

Nanowires represent the classical prototype of one-dimensional nanostructures and offer unique mechanical~\cite{Wang2017}, electronic, and optical properties~\cite{Li200618,Calahorra2017, Royo2017}. Thanks to these properties, nanowires have been used as elementary building blocks for the realization of nanoscale devices~\cite{Huang2001,Krogstrup2013} and in fundamental condensed matter physics~\cite{Mourik2012}. Among the various fabrication techniques, the {\em vapor-liquid-solid} growth method (further denoted VLS) is the most widespread~\cite{WagnerEllis1964}. Roughly speaking, the VLS method uses liquid nano-droplets catalysts on crystalline surfaces (as Au droplets for Si nanowires growth~\cite{Vincent2012} or 3rd group elements droplets for III-V nanowires~\cite{Fontcuberta2008}) and controlled external flux of atoms in a vacuum chamber, so that, at a certain concentration inside the droplets, solid phase nucleates at the liquid/solid interface. Repeated nucleation generates nanostructures headed by nano-droplets. Even if the basic bricks of the growth process are well-understood (see~\cite{PhysRevLett.99.146101,PhysRevB.88.195304,Tersoff,Jacobsson}), a unified model able to simulate the entire growth process at realistic time and length scales is still missing.

The basic physical phenomena involved in the VLS growth are: (a) superficial diffusion of adatoms at the free-substrate surface (b) incorporation of adatoms either by direct flux, ({\em i.e.} across the droplet free-surface) or by surface diffusion ({\em i.e.} across the boundary of the interface between the liquid droplet and the solid substrate), (c) diffusion of adatoms in the liquid droplet and (d) solidification at the liquid-solid interface. Clearly, (a), (b), and (d) are surface phenomena.

In the specific context of crystalline phenomena, either binary (two phases) or multiphase (more than two phases) phase field models have been proposed and studied for instance in~\cite{cahn1960theory,bellettini_1996,karma1996phase,fried1996phase,gurtin1999sharp,jiang2012phase,bellettini2013,Garcke_amulti}. A phase field model was successfully used to describe solidification and growth of binary and/or ternary alloys (see recent results in~\cite{boettinger2002phase,poulsen2016early,hotzer2015large,korbuly2017orientation}) as it allows for morphology changes, overlap of diffusion fields and particle coalescence or splitting~\cite{chen2002phase}. To the best of our knowledge, the challenge of a unified model able to describe nanowires growth
by VLS has been only partially adressed in~\cite{WangPhD,schwalbach2012stability,wang2014three}. We will propose and study in this paper an isotropic quasi-static Allen-Cahn multiphase field model as a rough approximation of the real anisotropic case. This is a first important step because, in contrast with the binary case where the phase field theory is well-understood, the multiphase case remains an active research topic.

\begin{figure}[htbp]
\centering
\includegraphics[scale=0.3]{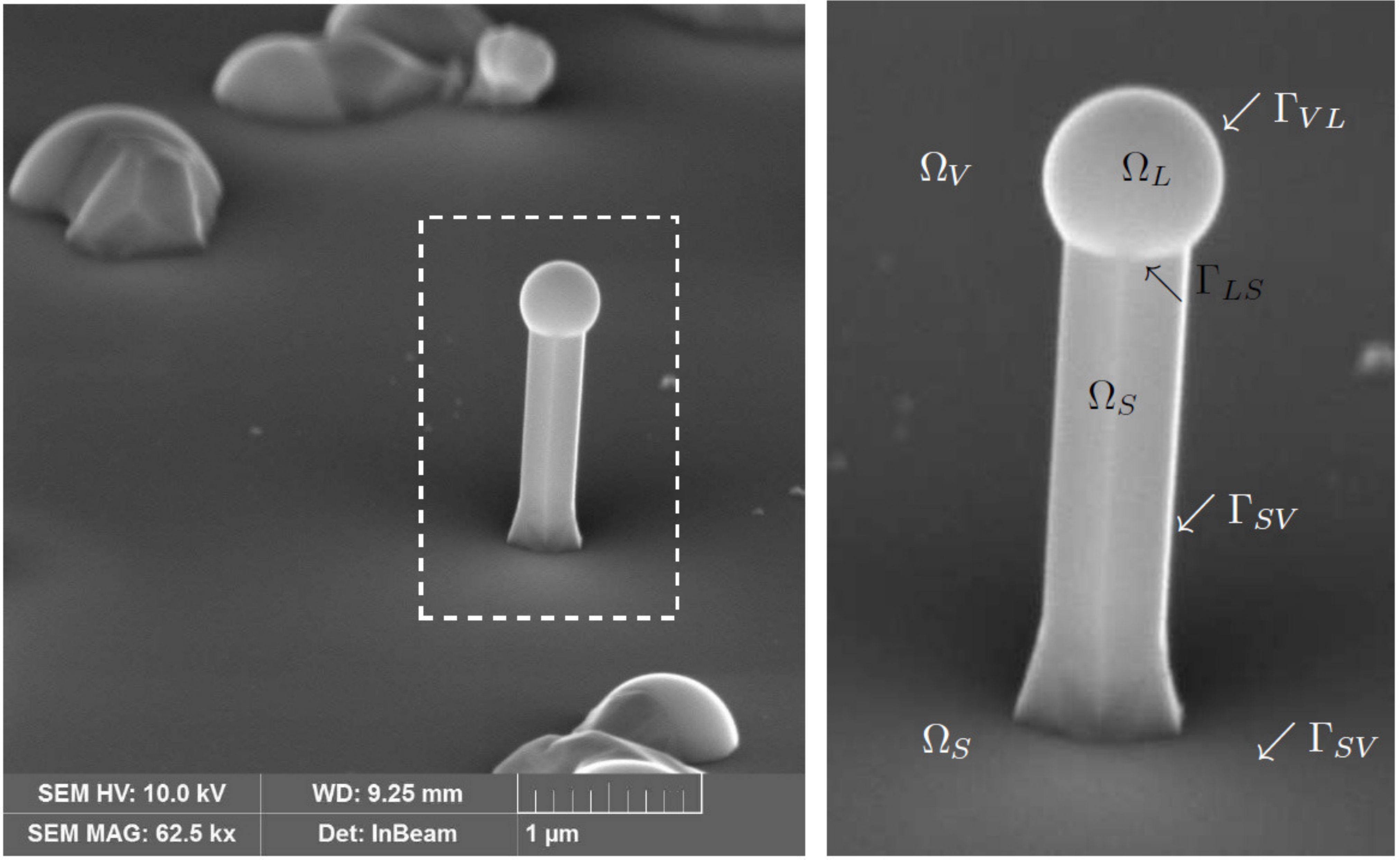}
\caption{Scanning electron microscopy image of an isolated nanowire on a cristalline substrate. The nanowire is made of GaAs which appears to be facetted, while the catalyst droplet is made of Ga and exhibits a spherical shape. The sample was grown on a Si(111) substrate by molecular beam epitaxy~\cite{Benali2017}. Left image includes an isolated nanowire (inset) and also several droplets (spherical caps) with crystalline material attached to it (left upper corner). This morphology is also a common situation in VLS growth. The right picture is a zoom of the inset in the left picture where we illustrate the domains $\Omega_S,$ $\Omega_L$ and $\Omega_V.$ In the actual (oblique) perspective the planar substrate, which is also a part of $\Omega_S$, covers the background of the image.}
\label{fig:fil_isole}
\end{figure}

We consider the VLS catalytic-growth process in the generic situation of steady growth conditions under constant {\em isotropic} external flux {\em of adatoms,} fixed substrate temperature, and for low density of nanowires. The shape of the nanowire and its near environnement can be modelled by using a partition of an open box $Q = \Omega_{V} \cup \Omega_{L} \cup \Omega_{S}\subset\R^d$ where $\Omega_{V}$, $\Omega_{L}$, and $\Omega_{S}$ are relatively closed subsets of $Q$ which represent the vapor, liquid and solid phases, and satisfy for all $i,j\in\{V,L,S\}$, $i\not=j$:
$$\Gamma_{ij}=\Omega_{i} \cap \Omega_{j}=\partial\Omega_{i}\cap\partial\Omega_{j} \cap Q,$$
with the convention $\Gamma_{ii}=\emptyset$. In Figure \ref{fig:fil_isole} we illustrate the sets $\Omega_{V},$ $\Omega_{L},$ $\Omega_{S}$ and the interfaces $\Gamma_{ij}$ by using an isolated nanowire obtained by the VLS catalytic-growth method. We notice the almost spherical shape of the catalytic droplet, which is a result of the minimization of the surface energy at the $LV$ interface, and we also notice the decrease of nanowire's radius in the first stages of the growth (see nanowire's foot).

Evolution of the nanowire shape is driven by the total interfacial energy
$$ J(\Omega_{V},\Omega_{L},\Omega_{S}) =   \sigma_{VL}\haus(\Gamma_{VL})  +  \sigma_{LS}\haus(\Gamma_{LS}) + \sigma_{SV}\haus(\Gamma_{SV}),$$
where $\haus$ is the $(d-1)$--dimensional Hausdorff measure~\cite{AFP}, and $\sigma_{LS}$, $\sigma_{VL}$ and $\sigma_{SV}$ represent
the surface tension coefficients between the liquid-solid, the vapor-liquid and
the solid-vapor phases, respectively.

As the $L^2$-gradient flow of the multiphase perimeter functional $J$ tends to minimize the surface energy, some additional constraints have to be imposed in order to include the essential qualitative features of the physical problem. Given an evolving partition $t\mapsto (\Omega_{L}(t),\Omega_{S}(t),\Omega_{V}(t))$ we shall assume that:

\begin{itemize}
\item The liquid phase volume is conserved, {\em i.e.}
$$\textrm{Vol}(\Omega_{L})'(t) =  0$$
which is indeed the case when the catalyser does not spread into the solid phase. This is the classical case of Au-catalysed growth of Si (or Ge) nanowires {\em although some studies have reported on the dissolution of the Au catalyst in the Si nanowires~\cite{Vincent2012}} but not that of III-V semiconductors self-catalysed nanowires\footnote{The case of self-catalysed nanowire growth will be discussed in a future work.}.
\item The velocity of the nanowire growth is proportional to the area of the solid-liquid interface, {\em i.e.}
$$\textrm{Vol}(\Omega_{S})'(t) =  c_{S} \haus(\Gamma_{SL}), \qquad \textrm{Vol}(\Omega_{V})'(t) =  -c_{S} \haus(\Gamma_{SL}).$$
Apparently restrictive, the first part of this assumption covers both the cases of nanowire size-dependent growth velocity valid for small nanowire radius and the case of size-dependent growth velocity valid for large nanowire radius (which is obviously covered by taking $c_S=0$). For theoretical models predicting these two regimes, the reader is referred to \cite{PhysRevB.88.195304}. The second part of the above assumptions is motivated by the physical realistic requirement that, since the liquid phase is conserved, the solid phase grows at the expense of the vapor phase.
\end{itemize}

By using classical Lagrange multipliers which account for both the partition and volume constraints, we obtain the following system
$$
\begin{cases}
\frac{1}{m_{LS}} v_{LS}(x,t) & =  \sigma_{LS} H_{LS}(x)  + \mu_{L}(t) +  \mu_{S}(t) + \lambda(x,t), \\
\frac{1}{m_{VL}} v_{VL}(x,t) & =  \sigma_{VL} H_{VL}(x)  + \mu_{L}(t) +  \mu_{V}(t) + \lambda(x,t) , \\
\frac{1}{m_{SV}} v_{SV}(x,t) & =  \sigma_{SV} H_{SV}(x)  + \lambda(x,t),
\end{cases}
$$
where $m_{\cdot\cdot}$ and $H_{\cdot\cdot}$ denote the mobilities and the mean curvature vectors at interfaces. As the solidification is located only at the solid-liquid interface $\Gamma_{LS}$, the  mobilities $m_{SV}$ and $m_{VL}$ should be chosen to be several orders of magnitude lower than $m_{LS}$. This physical requirement is particularly difficult to account for within classical phase field models.

\subsection{Classical phase field approximation (without mobilities)}

\subsubsection{Perimeter gradient flow and Allen-Cahn equation}
We start with the case of a single set $\Omega\subset\R^d$ which evolves following the {\em normal velocity law}
$$V_n = m H.$$
By a simple time rescaling, we can assume without loss of generality that $m=1.$ As above, $H$ denotes the mean curvature of $\partial\Omega$ and the evolution equation coincides with the $L^2$-gradient flow of the perimeter functional
$$P(\Omega) = \haus(\partial \Omega).$$
We recall that the basic principle of the phase field method is to replace the discontinuous characteristic function $\one{\Omega}$ by a smooth approximation $u$, and the singular perimeter energy by the smooth Van der Waals-Cahn-Hilliard functional \cite{CAHN:1958}
$$ P_{\varepsilon}(u) = \int_Q \left( \varepsilon \frac{|\nabla u|^2}{2} + \frac{1}{\varepsilon} W(u) \right)  dx,$$
where $\varepsilon$ is a small parameter and $W$ is a suitable double-well potential, e.g. $W(s) = \frac{1}{2} s^2 (1-s)^2.$
In this framework Modica and Mortola \cite{Modica1977} proved that $P_\varepsilon$  tends to $c_W \overline{P}$ in the sense of $\Gamma$-convergence for the $L^1$ topology, where $c_W = \int_0^1 \sqrt{2 W(s)} ds$ and
$$ \overline{P}(u) = \begin{cases}
                      |D u|(Q) \quad \text{ if }    u \in\BV(Q,\{0,1 \}), \\
                      + \infty \quad \text{otherwise}.
                     \end{cases}
$$
Here, $\BV$ denotes the space of functions with bounded variation in $Q\subset\R^d$, see~\cite{AFP}, and $\BV(Q,\{0,1 \})$ is the set of $\BV$ functions which take values in $\{0,1\}$. For $u\in\BV(Q)$, $|Du|(Q)$ denotes the total variation of $u$ in $Q$, defined as
$$ |D u |(Q) = \int_{Q} |D u| = \sup \left\{ \int_\Omega u \operatorname{div}g dx,\;  g \in C^1_0(Q,\R^d),\;\quad |g| \leq 1 \right\}. $$
In particular, when $\Omega$ is a set with finite perimeter $P(\Omega)$ in $Q$~\cite{AFP,Mag12a}, the characteristic function  $\one{\Omega}$ of $\Omega$ satisfies  $\overline{P}(1_{\Omega}) = P(\Omega)$. Modica and Mortola showed that $\one{\Omega}$ can be approximated by the smooth functions $u^{\varepsilon}(x) = q\left( \frac{\operatorname{dist}(x,\Omega)}{\varepsilon} \right)$ which satisfy $ \lim_{\varepsilon \to 0} P_{\varepsilon}(u^{\varepsilon}) = c_W P(\Omega)$ (as usual in the theory of $\Gamma$-convergence, $\epsilon\to 0$ must be intended in the sequential sense, i.e. it is related to a sequence $(\epsilon_n)_n$ such that $\lim\limits_{n\to\infty}\epsilon_n=0$). In the definition of $u^\varepsilon$,  $\dist(x,\Omega)$ denotes the signed distance function to $\Omega$ and $q$ is the so-called {\it optimal profile} associated with the potential $W$, defined by
$$q =  \mathop{\rm argmin} \limits_{p } \left\{ \int_{\R} \sqrt{ {W}(p(s))} |p'(s)| ds,\;  p(-\infty) = 1,\, p(0) = 1/2,\, p(+\infty) = 0 \right\},$$
where $p$ ranges over all Lipschitz continuous functions $p : \R \to \R$. A simple derivation of the Euler equation associated with this minimization problem shows that
\begin{equation} \label{def_q}
 q'(s) = -\sqrt{2 {W}(q(s))} \quad\text{and}\quad q''(s) = W'(q(s)), \quad  \text{ for all } s \in \R,
\end{equation}
which implies that $q(s) = (1- \tanh(s))/2$ in the case of the standard double well potential ${W}(s) = \frac1 2 s^2 (1-s)^2 $ considered above.\\

The classical Allen-Cahn equation \cite{Allen1979} is obtained as the $L^2$-gradient flow of the Van der Waals--Cahn--Hilliard energy $P_\varepsilon$, and up to time-rescaling reads as
$$   u_t =  \Delta u - \frac{1}{\varepsilon^2}W'(u).$$
For this equation, existence and uniqueness of a solution are well-known, as well as a comparison principle, see for example~\cite[Chap 14, 15]{Ambrosio2000}. A smooth set $\Omega$ evolving by mean curvature flow can be approximated by
$$ \Omega^\varepsilon(t) =  \left\{x \in \R^d,\; u^\varepsilon(x,t) \geq \frac{1}{2} \right\},$$
where $u^\varepsilon$ solves the Allen-Cahn equation with initial condition
$$u^{\varepsilon}(x,0) = q \left( \frac{\dist(x,\Omega(0))}{\varepsilon}\right).$$
A formal asymptotic expansion of $u^{\varepsilon}$ near the interfaces~\cite{BellettiniBook2013} shows that, at least formally, $u^{\varepsilon}$ is quadratically close to the optimal profile, i.e.
$$ u^{\varepsilon}(x,t) = q \left( \frac{\dist(x,\Omega^{\varepsilon}(t))}{\varepsilon}\right) + O(\varepsilon^2),$$
with associated normal velocity $V^{\varepsilon}$  satisfying
$$ V^{\varepsilon} = H + O(\varepsilon^2).$$
More rigorously, the convergence of $\partial \Omega_{\varepsilon}(t)$ to $\partial \Omega(t)$ has been proved for smooth motions by \cite{Chen1992, Mottoni1995, Bellettini1995} with a quasi-optimal convergence order $O(\varepsilon^2 | \log \varepsilon |^2)$.

\subsubsection{Multiphase field approximation in the additive case}
Using the inclusion-exclusion principle for a three-phase material, the nanowire surface energy can be expressed as
\begin{eqnarray*}
  J(\Omega_V,\Omega_L,\Omega_S) &=& \sigma_{ LS} \haus(\Gamma_{LS}) +
                                   \sigma_{ SV} \haus(\Gamma_{SV})  +
                                    \sigma_{ VL} \haus(\Gamma_{VL}) , \\
                                &=& \sigma_{V}\haus(\partial \Omega_V)  +
                                    \sigma_{L}\haus(\partial \Omega_L)  +
                                    \sigma_{S}\haus(\partial \Omega_S) ,
\end{eqnarray*}
where
$\sigma_{S} = \frac{1}{2} (\sigma_{ SV} + \sigma_{ LS} -  \sigma_{ VL}),$
$\sigma_{L} = \frac{1}{2} (\sigma_{ LS} + \sigma_{ VL} -  \sigma_{ SV})$ and
$\sigma_{V} = \frac{1}{2} (\sigma_{ VL} + \sigma_{ SV} -  \sigma_{ LS})$
are non-negative numbers due to the triangle inequality. This reformulation is important because it replaces interfacial areas by the area of {\it volume boundaries}, which opens the way to a phase-field approximation (recall that a phase field approximates characteristic functions of volume sets). The same principle applies for the general $N$-phase case when the surface tensions are {\it additive}, i.e., there exist nonnegative numbers $\{\sigma_i\}_i$ such that $\sigma_{ij} = \sigma_i + \sigma_j$ for any $i\not=j$ (when $N\leq 3$, {\it any} collection of surface tensions satisfying the triangle inequality is additive). In this very case, the $N$-phase perimeter functional can be written as
\begin{eqnarray*}
   P(\Omega_1, \Omega_2, \cdots,\Omega_N) &=& \frac{1}{2} \sum_{1\leq i < j\leq N}  \sigma_{ij}  \haus(\Gamma_{ij}) = \sum_{i}^{N}  \sigma_{i}  \haus(\partial \Omega_i).
\end{eqnarray*}

It appears immediately that $P$ can be approximated by the multiphase Cahn-Hilliard energy defined for all ${\boldsymbol u} = (u_1, u_2, \dots, u_N)$ by
$$P_{\varepsilon}({\boldsymbol u}) =
\begin{cases}
\frac{1}{2}\ds \sum_{i=1}^{N} \int_{Q} \sigma_i \left( \varepsilon \frac{|\nabla u_i|^2}{2} + \frac{1}{\varepsilon} W(u_i) \right) dx, &\text{if~} \sum_{i=1}^{N} u_i = 1 \text{,} \\
+ \infty & \text{otherwise.}
\end{cases}
$$
If ${\boldsymbol u}^{\varepsilon}$ is such that $\liminf P_{\varepsilon}({\boldsymbol u}^{\varepsilon} )< \infty$ then, by Modica-Mortola's Theorem \cite{Modica1977}, up to a subsequence, ${\boldsymbol u^{\varepsilon}} \to {\boldsymbol u} = (\one{\Omega_1},\one{\Omega_2}, \dots, \one{\Omega_N})$ and the constraint $\sum_{i=1}^{N} u_i = 1$ ensures that ${\boldsymbol \Omega} = \{\Omega_1,\Omega_2, \dots, \Omega_N \}$ is a partition of $Q$.\\

The $\Gamma$-convergence of $P_{\varepsilon}$  to $c_{W} P$ is was established in~\cite{oudet2011} for the particular case where $\sigma_i = 1$, $\forall i$. More general $\Gamma$-convergence results were obtained in \cite{baldo1990,Bretin_Masnou_multiphase} for inhomogeneous surface tensions $\sigma_{ij}$, while multiphase field models in the context of anisotropic surface tensions were introduced and analyzed in \cite{Garcke_amulti,Garcke199887}.

The $L^2$-gradient flow of $P_{\varepsilon}$ leads to the following Allen-Cahn system
\begin{equation} \label{eqn:allencahn_simple}
 \partial_t u^{\epsilon}_k = \sigma_k \left[ \Delta u^{\epsilon}_k - \frac{1}{\varepsilon^2} W'(u^{\epsilon}_k) \right] + \lambda^{\epsilon}, \quad   \forall k=1, \dots, N \text{,}
\end{equation}
where the Lagrange multiplier field $\lambda^{\epsilon}$ encodes the partition constraint $ \sum_{k=1}^{N} u^{\epsilon}_k = 1$. \\

By using the method of matched asymptotic expansions developed in \cite{fife,pego,belpaoqo,Loreti_march,Garcke199887}, we will prove in this paper the folllowing result:
\begin{clm} \label{claim:allencahn_simple}
Denoting $\Omega^{\epsilon}_i=\left\{ x \in Q ; u_i(x,t) \geq \frac{1}{2} \right\}$,  the solution ${\boldsymbol u}^{\varepsilon}$  of \eqref{eqn:allencahn_simple} expands formally near the interface $\Gamma_{ij}$ as
 $$
\begin{cases}
 u_i^{\varepsilon} &= q \left( \frac{ dist(x,\Omega^{\varepsilon}_i)}{\varepsilon} \right) + O(\varepsilon),\\
 u_j^{\varepsilon} &=  1 - q \left( \frac{ dist(x,\Omega^{\varepsilon}_i)}{\varepsilon} \right) + O(\varepsilon),\\
 u_k^{\varepsilon} &=  O(\varepsilon), \text{ for k } \in \{1,2, \dots, N\}\setminus\{i,j\}.\\
\end{cases}
$$
Moreover, the associated normal velocity satisfies $V_{ij}^{\varepsilon} = \frac{1}{2}\sigma_{ij} H_{ij} + O(\varepsilon).$
\end{clm}

These formal results indicate that the phase-field model~\eqref{eqn:allencahn_simple} is consistent but its solutions converge not faster than a linear order in $\varepsilon$.\\

In order to improve the convergence order we consider the slightly modified phase field system
\begin{equation} \label{eqn:allencahn_modif}
  \partial_t u^{\epsilon}_k = \sigma_k \left[ \Delta u^{\epsilon}_k - \frac{1}{\varepsilon^2} W'(u^{\epsilon}_k) \right] + \lambda^{\epsilon} \sqrt{2W(u_k)} \quad   \forall k=1, \dots, N,
\end{equation}
where, again, the Lagrange multiplier field $\lambda^{\epsilon}$ encodes the partition constraint $\sum_{k=1}^{N} u^{\epsilon}_k = 1.$
The idea to localize the Lagrange multiplier $\lambda$ near the diffuse interface has been recently proposed in \cite{Bretin_Denis_oudet_Jacquo} in order to improve the accuracy of the two-phases model. In our case, we will show that, at least formally,

\begin{clm} \label{claim:allencahn_modif} Near $\Gamma_{ij},$ the solution ${\boldsymbol u}^{\varepsilon}$  of \eqref{eqn:allencahn_modif} expands formally as
$$
\begin{cases}
u_i^{\varepsilon} &= q \left( \frac{ dist(x,\Omega^{\varepsilon}_i)}{\varepsilon} \right) + O(\varepsilon^2),\\
u_j^{\varepsilon} &=  1 - q \left( \frac{ dist(x,\Omega^{\varepsilon}_i)}{\varepsilon} \right) + O(\varepsilon^2),\\
u_k^{\varepsilon} &=  O(\varepsilon^2), \text{ for k } \in \{1,2, \dots, N\}\setminus\{i,j\},\\
\end{cases}
$$
with $V_{ij}^{\varepsilon} = \frac{1}{2}\sigma_{ij} H_{ij} + O(\varepsilon).$
\end{clm}
Therefore, solutions to the new model~\eqref{eqn:allencahn_modif}, converges quadratically (at least formally) to the optimal profiles.

\subsection{Incorporation of mobilities}
\subsubsection{The energetic viewpoint}
It was proposed in \cite{Garcke_amulti} to incorporate the mobilities $m_{ij}$ directly in the Cahn-Hilliard energy by considering a model of the form
$$P_{\epsilon}({\boldsymbol u}) =  \int_{Q} \epsilon f({\boldsymbol u},\nabla {\boldsymbol u}) + \frac{1}{\epsilon} {W}({\boldsymbol u}) dx,$$
where $f({\boldsymbol u},\nabla {\boldsymbol u}) = \sum_{i<j} m_{ij} \sigma_{ij} \left| u_i \nabla u_j - u_j \nabla u_i \right|^2$ and the multi-well potential ${W}$ is defined by
$${W}({\boldsymbol u}) = 9 \sum_{i,j=1,  i<j}^{N} \frac{\sigma_{ij}}{m_{ij}} u_i^2 u_j^2  + \sum_{i<j<k} \sigma_{ijk} u_i^2 u_j^2 u_k^2.$$
The term  $\sum_{i<j<k} \sigma_{ijk} u_i^2 u_j^2 u_k^2$ can be regarded as a penalization, with the coefficient $\sigma_{ijk}$ being chosen sufficiently large to ensure
the convergence of $P_{\epsilon}$ to the multiphase field perimeter $P$. As the mobility appears in the energy (and not only in the flow of $P_{\epsilon}$), it can be
expected that the size of the diffuse interface $\Gamma_{ij}$ depends on the mobility $m_{ij}$. This can be easily seen already for the two-phases case. Take indeed the modified Cahn-Hilliard energy
$$ P_{\varepsilon}^m(u) = P_{m\epsilon}(u)= \int_Q \left(m\varepsilon \frac{|\nabla u|^2}{2} + \frac{1}{m\varepsilon} W(u) \right)  dx,$$
for which the classical $L^2$-gradient flow reads as (up to time rescaling):
$$   u_t =  m\Delta u - \frac{1}{m\varepsilon^2}W'(u).$$
The matched asymptotic expansion method~\cite{fife,pego,belpaoqo} applied to this equation gives a solution of the form:
$$u(x,t)=q \left(\frac {\operatorname{dist}(x,\Omega^\varepsilon)} {\varepsilon m} \right).$$
with the associated velocity law:
$$V^\varepsilon=mH+O(\epsilon^2).$$
Obviously, the mobility plays a role in the size of the interface. Generalizing to the multiphase case, each mobility $m_{ij}$ will impact explicitly the size of the diffuse interface associated with $\Gamma_{ij}$. From the numerical point of view, this approach raises a significant limitation for high contrast mobilities and, in particular, it cannot be used for degenerate, i.e., vanishing mobilities.

Another issue with this model is related to the difference of nature between surface tensions and mobilities: whereas surface tensions are geometric parameters which appear in the sharp energy, mobilities are typically evolution parameters which play a role out of equilibrium. This is why we propose to incorporate mobilities not in the {\it geometric energy} but rather in the {\it metric} used for defining the gradient flow.

\subsubsection{A new approach: the metric viewpoint}
We propose a novel approach to incorporate the mobility so as to handle the special case of degenerate mobilities.
The idea is to introduce the gradient flow of $P_{\epsilon}$ with respect to a weighted scalar product
$$\langle  {\boldsymbol u},{\boldsymbol v} \rangle_{L_A^2(Q,\R^N)} =  \int_Q  ({\boldsymbol A}{\boldsymbol u}) \cdot {\boldsymbol v} dx$$
where the matrix ${\boldsymbol A}$ depends on the mobilities $m_{ij}$. The advantage of such approach can be easily seen in the binary case: defining the new scalar product
$$\langle  u,v\rangle_{L_m^2(Q)} =  \int_Q  \frac 1 m u\,v dx,$$
the $L^2_m$-gradient flow of the classical Cahn-Hilliard energy is
$$u_t =  m(\Delta u - \frac{1}{\varepsilon^2}W'(u)).$$
whose associated optimal profile is
$$u(x,t)=q \left(\frac {\operatorname{dist}(x,\Omega^\varepsilon)} {\varepsilon } \right).$$
with the evolution law
$$V^\varepsilon=mH+O(\epsilon^2).$$
Of course, the new metric cannot be explicitly defined for vanishing mobilities. However, because the size of the interface does not depend on m, arbitrary small values of m can be used without any loss in numerical efficiency. We shall now explain how, in the multiphase case, ${\boldsymbol A}$ can be first easily defined for a specific class of mobilities (which will be called {\it harmonically additive}), and then for more general mobilities.

\paragraph{Harmonically additive mobilities}~\\
Let us first assume an additivity property for the mobility coefficients, i.e, there exist some non negative coefficients $m_{i}$ such that
$$ \frac{1}{m_{ij}} = \frac{1}{m_i} + \frac{1}{m_j}.$$
This assumption has no clear physical justification, yet, it makes sense in a few situations. For instance, in the case of three phases, if the mobility coefficients satisfy the harmonic triangle inequality, i.e. $\frac 1{m_{ij}}\leq \frac 1{m_{ik}}+\frac 1{m_{kj}}$, then the addivity property is satisfied. From the modeling viewpoint, such assumption has a clear consequence: it yields a second order approximation for a suitable choice of the metric, namely ${\boldsymbol A} = {\boldsymbol M}^{-1}$ with ${\boldsymbol M}$ given by
$$ M_{ij} =  \begin{cases}
               m_i & \text{ if }  i=j \\
               0& \text{otherwise}
              \end{cases}
$$
The associated $L^2_{\boldsymbol A}$-gradient flow gives the following Allen-Cahn system:
\begin{equation} \label{eqn:allencahn_mob_additive}
 \partial_t u^{\epsilon}_k =   m_k \left[ \sigma_k \left( \Delta u^{\epsilon}_k - \frac{1}{\epsilon^2} W'(u^{\epsilon}_k) \right) +  \lambda^{\epsilon} \sqrt{2 W(u_k)} \right], \quad \forall k \in \{1,2,\cdots,N\}
\end{equation}
where the Lagrange multiplier field $\lambda^{\epsilon}$ is again associated to the constraint  $\sum u_i^{\epsilon} = 1$.

For this particular category of mobilities, we will derive formally the following result:
\begin{clm} \label{claim:allencahn_mob_additive} Around the interface $\Gamma_{ij}$,  the solution ${\boldsymbol u}^{\varepsilon}$  to \eqref{eqn:allencahn_mob_additive} has the following form (at least formally):
$$
\begin{cases}
 u_i^{\varepsilon} &= q \left( \frac{ dist(x,\Omega^{\varepsilon}_i)}{\varepsilon} \right) + O(\varepsilon^2),\\
 u_j^{\varepsilon} &=  1 - q \left( \frac{ dist(x,\Omega^{\varepsilon}_i)}{\varepsilon} \right) + O(\varepsilon^2),\\
 u_k^{\varepsilon} &=  O(\varepsilon^2), \text{ for k } \in \{1,2, \dots, N\}\setminus\{i,j\},\\
\end{cases}
$$
with $\frac{1}{m_{ij}} V_{ij}^{\varepsilon} = \sigma_{ij} H_{ij} + O(\varepsilon)$.
\end{clm}
Therefore, this multiphase field model with harmonically additive mobilities has quadratic convergence in~$\epsilon$.

\paragraph{General mobilities}~\\
The above additivity assumption is clearly not always satisfied, for instance $m_{12}=1$, $m_{13}=1$, and $m_{23}<0.5$ is a natural choice of mobility coefficients for which the harmonic inequality fails. However, even for general mobilities, one can still prove a convergent property, yet less optimal, i.e. of order one rather than two. Define indeed
$$ A_{ij} =
\begin{cases}
-\frac{1}{m_{ij}} & \text{if}~  i\neq j \\
0                 & \text{if} ~ i=j.
\end{cases}
$$
whose associated $L^2$-gradient flow is
\begin{equation} \label{eqn:allencahn_mob_gen}
  {\boldsymbol A} \partial_t {\boldsymbol u }^{\epsilon}=   \sigma \Delta {\boldsymbol u}^{\epsilon} - \frac{1}{\epsilon^2} W'({\boldsymbol u}^{\epsilon})  + \lambda^{\epsilon} \sqrt{2 W({\boldsymbol u}^{\epsilon})},
\end{equation}
where, for all $k \in \{1,2,\cdots,N \},$
$$(\sigma \Delta {\boldsymbol u})_k = \sigma_k \Delta u_k, \quad W'({\boldsymbol u})_k = W'(u_k),\quad \text{ and } \quad (\sqrt{2 W({\boldsymbol u})})_k  = \sqrt{2W(u_k)}.$$
The Allen-Cahn system~\eqref{eqn:allencahn_mob_gen} is well-posed as soon as ${\boldsymbol A}$ is semi-definite positive on $(1,1,\cdots,1)^{\perp},$ which in turn imposes some
restriction on the choice of the mobility $m_{ij}$ (see~\cite{Bretin_Masnou_multiphase} for a similar discussion about surface tensions).\\

We will show, at least formally, that the solution  ${\boldsymbol u}^{\varepsilon}$ to~\eqref{eqn:allencahn_mob_gen} has the following form.
\begin{clm} \label{claim:allencahn_mob_gen} Around the interface $\Gamma_{ij}$,  the solution ${\boldsymbol u}^{\varepsilon}$  to \eqref{eqn:allencahn_mob_gen} is, at least formally, of the form
$$
\begin{cases}
 u_i^{\varepsilon} &= q \left( \frac{ dist(x,\Omega^{\varepsilon}_i)}{\varepsilon} \right) + O(\varepsilon),\\
 u_j^{\varepsilon} &=  1 - q \left( \frac{ dist(x,\Omega^{\varepsilon}_i)}{\varepsilon} \right) + O(\varepsilon),\\
 u_k^{\varepsilon} &=  O(\varepsilon^2), \text{ for k } \in \{1,2, \dots, N\}\setminus\{i,j\},\\
\end{cases}
$$
with
$$  \frac{1}{m_{ij}} V_{ij}^{\varepsilon} = \sigma_{ij} H_{ij} + O(\varepsilon).$$
\end{clm}
We conclude that this multiphase field model is of order one only. In practice, of course, as soon as the harmonic additivity is satisfied, we shall opt for the model described in the previous section.

\subsubsection{Application to the modeling of nanowire growth}\label{sub:appmodel}

We recall that the nanowire growth can be modeled as the evolution of a partition ${\boldsymbol{\Omega}} = (\Omega_L,\Omega_S,\Omega_V)$ with the folllowing velocities at boundaries:
$$
\begin{cases}
\frac{1}{m_{ LS}} v_{ LS} &=  \sigma_{ LS} H  + \mu_L(t) + \mu_S(t)  +  \lambda, \\
\frac{1}{m_{ VL}} v_{ VL} &=  \sigma_{ VL} H  + \mu_L(t) + \mu_V(t) +  \lambda , \\
\frac{1}{m_{ SV}} v_{ SV} &=  \sigma_{ SV} H  + \lambda.
\end{cases}
$$
where $\lambda(x,t),$  $\mu_L(t)$, $\mu_V(t)$  and  $\mu_S(t)$ correspond to the Lagrange multipliers associated to the partition constraint and to the volume constraint, respectively.
We further assume that the mobility coefficients satisfy
$$m_{ LS}= \frac{\delta}{1 + \delta}, \quad  m_{ VL}=\frac{1}{2},  \text{ and } \quad  m_{ SV} = \frac{\delta}{1 + \delta}$$
with $\delta\ll 1$. This choice is motivated by interfaces $LS$ and $SV$ having much smaller mobilities than the $VL$ interface. Moreover, it corresponds to the harmonically additive
case with $m_S = \delta$ and  $m_L = m_V = 1.$ \\

The phase field approximation of this model is given by the  $M^{-1}$ $L^2$-gradient flow of the multiphase Cahn Hilliard energy
$${P_{\epsilon}}({\boldsymbol u}) =  \sigma_L P_{\epsilon}(u_L) +  \sigma_V P_{\epsilon}(u_V) +  \sigma_S P_{\epsilon}(u_S), \text{ with }  {\boldsymbol u} = (u_L,u_S,u_V),$$
i.e.
$$
\begin{cases}
 \frac{1}{m_L} \partial_t u_L(x,t) &= \sigma_L \left[ \Delta u_L - \frac{1}{\epsilon^2} W'(u_L)\right] + \lambda\sqrt{2 W(u_L)} + \mu_L \sqrt{2 W(u_L)} \\
 \frac{1}{m_S} \partial_t u_S(x,t) &= \sigma_S \left[ \Delta u_S - \frac{1}{\epsilon^2} W'(u_S)\right] + \lambda\sqrt{2 W(u_S)} + \mu_S u_L u_S  \\
 \frac{1}{m_V} \partial_t u_V(x,t) &= \sigma_S \left[ \Delta u_V - \frac{1}{\epsilon^2} W'(u_V)\right] + \lambda\sqrt{2 W(u_V)} + \mu_V u_L u_V \\
\end{cases}
$$
where $\lambda(x,t),$  $\mu_L(t)$, $\mu_V(t)$, and  $\mu_S(t)$ encode both the partition constraint
$$  u_L + u_S + u_V = 1, $$
and the volume constraints :
$$  \int_{Q} \partial_t  u_L(x,t) dx = 0, \text{ and }    \int_{Q}  \partial_t u_S(x,t) dx =  -  \int_{Q}  \partial_t u_V(x,t) dx = \frac{c_S}{\epsilon}\int_{Q} u_L(x,t)u_S(x,t) dx.$$

\subsection{Outline of the paper}

The paper is organized as follows: Claims 1-4 are proven in Section 2 using the method of matched asymptotic expansions.
We describe in Section 3 a numerical method to approximate the solutions of a multiphase mean curvature flows with mobilities and possible additional volume constraints. We provide various examples of such multiphase flows in order to illustrate the influence of mobility. We show in particular, with several examples of a droplet wetting on a solid surface, that our method is well-suited numerically for handling degenerate or highly contrasted mobilities. Finally, we apply our method to simulate an isotropic approximation of a nanowire grown by the VLS method. Our numerical results related to nanowires are confirmed by a theoretically derived optimal profile, whose derivation is new to the best of our knowledge.

\section{Asymptotic expansion of solutions to the Allen-Cahn systems}
This section is devoted to the (formal) identification of sharp interface limits of solutions ${\boldsymbol u}^\varepsilon=(u_1^\varepsilon,\dots,u^{\varepsilon}_{N})$ to the Allen-Cahn systems \eqref{eqn:allencahn_simple}, \eqref{eqn:allencahn_modif}, \eqref{eqn:allencahn_mob_additive}, and  \eqref{eqn:allencahn_mob_gen}. To this aim, we use the formal method of matched asymptotic expansions
proposed in  \cite{fife,pego,belpaoqo,Loreti_march}, which we apply around each interface $\Gamma_{ij}$. ~\\

\subsection{Preliminaries}
\paragraph{Outer expansion far from $\Gamma_{ij}$:} We assume that the {\em outer expansion} of $u^{\varepsilon}_{k}$, {\em i.e.,} the expansion far from the front $\Gamma_{ij}$ has the form:
$$u^{\varepsilon}_k(x,t) = u^0_k(x,t) + \varepsilon u^1_k(x,t) +  O(\varepsilon^2), \text{ for all } k \in \{1,2, \cdots,N\}.$$
In particular and analogously to \cite{Loreti_march}, it is not difficult to see that if $E_i(t) = \{ x \in \Omega, u^{\epsilon}_i \geq \frac{1}{2}\},$ then
$$u^0_i(x,t) =
\begin{cases}
 1 & \text{ if } x \in E_i(t) \\
 0 & \text{otherwise}
\end{cases}, \quad  u^0_j(x,t) =
\begin{cases}
 0 & \text{ if } x \in E_i(t) \\
 1 & \text{otherwise}
\end{cases}
$$
and $ u^1_i = u^1_j = 0,$ $u^0_k = u^1_k = 0$ for all $k \in \{1,2, \cdots N\}\setminus \{i,j\}.$

\paragraph{Inner expansions around $\Gamma_{ij}$:} In a small neighborhood of $\Gamma_{ij}$, we define the stretched normal distance to the front as
$ z = \frac{1}{\varepsilon}{d_i(x,t)},$ where $d_i(x,t)$ denotes the signed distance to $E_i(t)$ such that $d_i(x,t)<0$ in $E_i(t)$.
The {\em inner expansions} of $u^{\varepsilon}_k(x,t)$ and $\lambda^{\epsilon}(x,t)$, {\em i.e.} expansions close to the front, are assumed of the form
$$ u^{\varepsilon}_k(x,t) = U^{\varepsilon}_{k}(z,x,t) = U^{0}_{k}(z,x,t) + \varepsilon U^{1}_k(z,x,t) + O(\varepsilon^2),  \text{ for all } k \in \{1,2, \cdots,N\},$$
and
$$  \lambda^{\epsilon}(x,t) =  \Lambda^{\varepsilon}(z,x,t)  =  \varepsilon^{-2} \Lambda^{-2}(z,x,t) +  \varepsilon^{-1} \Lambda^{-1}(z,x,t) + O(1).$$
Moreover, if $n$ denotes the unit normal to $\Gamma_{ij}$ and $V^{\epsilon}_{ij}$ the normal velocity to the front, for $x\in\Gamma_{ij}$
$$ V^{\varepsilon}_{ij} = \partial_t d_i(x,t) = V^{0}_{ij} + \varepsilon V^1_{ij} + O(\varepsilon^2) , \quad n = \nabla d_i(x,t).$$
where $\nabla$ refers to the spatial derivative only.

Following~\cite{pego,Loreti_march} we assume that
$U^{\epsilon}_k(z,x,t)$ does not change when $x$ varies normal to $\Gamma_{ij}$ with $z$ held fixed, or equivalently
$( \nabla U^{\epsilon}_k )_{z={\rm const.}}\cdot n = 0$. This amounts to requiring that the blow-up with respect to the parameter $\varepsilon$ is coherent with the flow.

Following~\cite{pego,Loreti_march}, it is easily seen that
$$
\begin{cases}
 \nabla u_k = \nabla_x U_k + \varepsilon^{-1} n \partial_z U_k, \\
 \Delta u_k =\Delta_x U_k + \varepsilon^{-1} \Delta d_i ~\partial_z U_k + \varepsilon^{-2} \partial^2_{zz} U_k, \\
 \partial_t u_k = \partial_t U_k - \varepsilon^{-1} V_{ij}^{\epsilon} \partial_z U_{k}.
\end{cases}
$$
Recall also that in a sufficiently small neighborhood of $\Gamma_{ij}$, according to Lemma 14.17 in~\cite{GT}, we have
$$ \Delta d_i (x,t) = \sum_{k=1}^{d-1} \frac{\kappa_k(\pi(x))}{1 + \kappa_k(\pi(x)) d_i(x,t)} = \sum_{k=1}^{d-1} \frac{\kappa_k(\pi(x))}{1 + \kappa_k(\pi(x)) \varepsilon z },$$
where $\pi(x)$ is the projection of $x$ on $\Gamma_{ij}$ and $\kappa_k$ are the principal curvatures on $\Gamma_{ij}$.
In particular this implies that
$$ \Delta d_i (x,t) = H_{ij} - \varepsilon z \|A_{ij}\|^2 + O(\varepsilon^2), $$
where $H_{ij}$ and $\|A_{ij}\|^2$ denote, respectively, the mean curvature and the squared $2$-norm of the second fundamental form on $\Gamma_{ij}$ at $\pi(x)$.

\paragraph{Matching conditions between outer and inner expansions:} The matching conditions (see \cite{Loreti_march} for more details) can be written as:
$$   \lim_{z \to + \infty} U^0_i(z,x,t) = 0, \lim_{z \to - \infty} U^0_i(z,x,t) = 1, \quad \lim_{z \to \pm \infty} U^1_i(z,x,t) = 0,$$
$$   \lim_{z \to + \infty} U^0_j(z,x,t) = 1, \lim_{z \to - \infty} U^0_j(z,x,t) = 0, \quad \lim_{z \to \pm \infty} U^1_j(z,x,t) = 0,$$
and
$$    \lim_{z \to \pm \infty} U^0_k(z,x,t) =  \lim_{z \to \pm \infty} U^1_k(z,x,t) = 0, \text{for all } k \in \{1,2, \cdots N\}\setminus \{i,j\}.$$

\subsection{Analysis of the classical additive Allen-Cahn system (no mobilities)}

We first consider the classical additive Allen-Cahn system~\eqref{eqn:allencahn_simple}, i.e. the set of equations
$$ \partial_t u^{\epsilon}_k =  \sigma_k \left( \Delta u^{\epsilon}_k - \frac{1}{\epsilon^2} W'(u^{\epsilon}_k) \right) +  \lambda^{\epsilon}$$
where the Lagrange multiplier field $\lambda^{\epsilon}$ encoding the pointwise constraint $\sum_{j=1}^N u_j = 1$ can be explicitly computed as
\begin{equation} \label{def_lambda_allencahn_simple}
\lambda^{\epsilon} = - \frac{1}{N} \sum_{j=1}^N   \sigma_j \left( \Delta u^{\epsilon}_j - \frac{1}{\epsilon^2} W'(u^{\epsilon}_{j}) \right).
\end{equation}
We focus first on inner expansions of $u^{\varepsilon}_k(x,t)$ and $\lambda^{\epsilon}(x,t)$, i.e. expansions close to the front $\Gamma_{ij}$. Injecting into  \eqref{eqn:allencahn_simple} and \eqref{def_lambda_allencahn_simple} the following expressions:
$$ u^{\varepsilon}_k(x,t) = U^{\varepsilon}_{k}(z,x,t) = U^{0}_{k}(z,x,t) + \varepsilon U^{1}_k(z,x,t) + O(\varepsilon^2),  \text{ for all } k \in \{1,2, \cdots,N\},$$
and
$$  \lambda^{\epsilon}(x,t) =  \Lambda^{\varepsilon}(z,x,t)  =  \varepsilon^{-2} \Lambda^{-2}(z,x,t) +  \varepsilon^{-1} \Lambda^{-1}(z,x,t) + O(1).$$
leads to the following terms at various orders.

\paragraph{Order $\varepsilon^{-2}$ :} Identifying the terms of order $\varepsilon^{-2}$ in \eqref{eqn:allencahn_simple} and \eqref{def_lambda_allencahn_simple} gives
$$
     \sigma_k \left( \partial^2_{zz} U^0_k - W'(U^0_k) \right)  + \Lambda^{-2} = 0, \text{ for all } k \in \{1,2, \cdots N\},
$$
and
$$  \Lambda^{-2}  = - \frac{1}{N} \sum_{k=1}^{N} \left(  \partial^2_{zz} U^0_k - W'(U^0_k)  \right).$$
Moreover, the boundary conditions obtained from the matching conditions and the equality $U^0_i(0,x,t) =\frac 1 2$ give (recall that $q$ is the optimal profile defined as the solution to~\eqref{def_q}):
$$
\begin{cases}
 U^0_i(z,x,t)  &= q(z), \\
  U^{0}_j(z,x,t) &= q(-z) = 1-q(z), \\
   U^{0}_k(z,x,t) &= 0, \text{ for all } k \in \{1,2,\cdots,N\} \setminus\{i,j\} \\
    \Lambda^{-2} = 0.
\end{cases}
$$

\paragraph{Order $\varepsilon^{-1}$:} Matching the terms of order $\varepsilon^{-1}$ in \eqref{eqn:allencahn_simple} and \eqref{def_lambda_allencahn_simple} gives
$$ V_{ij} \partial_z U^0_{k} =  \sigma_k \left[ \partial^2_{zz} U^{1}_k - W''(U^0_k)U^1_k  +   H_{ij}  \partial_z U^0_{k} \right]  +  \Lambda^{-1},  $$
and
$$  \Lambda^{-1}  = - \frac{1}{N} \sum_{k=1}^{N}  \sigma_k \left[  \partial^2_{zz} U^1_k - W''(U^0_k) U^1_k + H_{ij} \partial_z U^0_{k} \right].$$
In particular, for $k=i$ and $k=j$ we obtain respectively
 \begin{align}
       V^0_{ij} q'(z) &=  \sigma_i \left( \partial^2_{zz} U^{1}_i - W''(q(z)) U^1_i \right) +  \sigma_i q'(z) H_{ij}  +  \Lambda^{-1}, \nonumber \\
     - V^0_{ij} q'(z) &=  \sigma_j \left( \partial^2_{zz} U^{1}_j - W''(q(z)) U^1_j \right) -  \sigma_j q'(z) H_{ij}  +  \Lambda^{-1}, \nonumber
   \end{align}
so that
$$ 2 V^0_{ij} q'(z)  =   \sigma_i \left( \partial^2_{zz} U^{1}_i - W''(q(z)) U^1_i \right) - \sigma_j \left( \partial^2_{zz} U^{1}_j - W''(q(z)) U^1_j \right) +  (\sigma_i + \sigma_j) q'(z) H_{ij}.$$
Multiplying this equation by $q',$ integrating over $\R$ and taking into account the matching conditions we obtain
$$V^0_{ij} =  \frac{1}{2} (\sigma_i + \sigma_j) H_{ij} = \frac{1}{2} \sigma_{ij} H_{ij},$$
which shows that the first order term of the interface velocity matches with the expected velocity law. Moreover, the Lagrange multiplier $ \Lambda^{-1}$ satisfies
$$ \int_{\R} \Lambda^{-1}(z,x,t) q'(z) dz = \frac{c_W}{2} (\sigma_j - \sigma_i) H_{ij}(x,t),$$
where  $c_W = \int_{\R} (q'(s))^2 ds = \int_0^{1} \sqrt{2 W(s)} ds$. This shows that $ \Lambda^{-1}$ is expected of the form
$$ \Lambda^{-1}(z,x,t) = \frac{1}{2} [\sigma_j - \sigma_i] H_{ij}(x,t)~\eta(z).$$

For each $k\in\{1,2,\ldots, N\}\setminus\{i,j\},$ the functions $U^1_k$ are defined as the solutions to
$$
\begin{cases}
  \sigma_i \left( \partial^2_{zz} U^{1}_i - W''(q(z)) U^1_i \right) = \frac{1}{2} [\sigma_j - \sigma_i] H_{ij}(x,t) \left( q'(z) - \eta(z) \right), \\
  \sigma_j \left( \partial^2_{zz} U^{1}_j - W''(q(z)) U^1_j \right) = \frac{1}{2} [\sigma_j - \sigma_i] H_{ij}(x,t) \left( q'(z) - \eta(z) \right), \\
  \sigma_k \left( \partial^2_{zz} U^{1}_k - W''(0)    U^1_k \right) = \frac{1}{2} [\sigma_j - \sigma_i] H_{ij}(x,t)  \eta(z).
\end{cases}
$$
with additional Dirichlet limit conditions at $z=\pm\infty$.  The profile $\eta$ can be also obtained by imposing the additional constraint  $\sum_{k=1}^{N}U^{1}_k = 0$.
Finally, the asymptotic expansion shows that the second term $U^{1}_k$ does not vanish in general as soon as $\sigma_i \neq \sigma_j$ which proves Claim~\ref{claim:allencahn_simple}.

\subsection{Analysis of the modified additive Allen-Cahn system (no mobilities)}
We now analyze the modified additive Allen-Cahn system~\eqref{eqn:allencahn_modif}, i.e. the system
$$ \partial_t u^{\epsilon}_k =  \sigma_k \left( \Delta u^{\epsilon}_k - \frac{1}{\epsilon^2} W'(u^{\epsilon}_k) \right) +  \lambda^{\epsilon} \sqrt{2 W(u^{\epsilon}_k)}$$
with $\lambda^{\epsilon}$ defined as
\begin{equation}\label{def_lambda_allencahn_modif}
\lambda^{\epsilon} = \frac{\sum_{j=1}^N \left(   \sigma_j \left( \Delta u^{\epsilon}_j - \frac{1}{\epsilon^2} W'(u^{\epsilon}_{j}) \right)  \right)}{ \sum_{j=1}^N \sqrt{2 W(u^{\epsilon}_j)}}.
\end{equation}

\paragraph{Order $\varepsilon^{-2}: $} The first order equations \eqref{eqn:allencahn_modif} and \eqref{def_lambda_allencahn_modif}  read now respectively as
$$
     \sigma_k \left( \partial^2_{zz} U^0_k - W'(U^0_k) \right)  + \Lambda^{-2} \sqrt{2 W(U^0_k)} = 0, \text{ for all} \quad k \in \{1,2, \cdots N\},
$$
so that
$$\Lambda^{-2}  \sum_{k=1}^{N}  \sqrt{2 W(U^0_k)}  = - \sum_{k=1}^{N} \sigma_k \left(  \partial^2_{zz} U^0_k - W'(U^0_k)  \right).$$
As previously, from the boundary conditions we obtain
$U^0_i(z,x,t)  = q(z),$ $U^{0}_j(z,x,t) = 1 - q(z),$ $U^{0}_k(z,x,t) = 0$ and  $\Lambda^{-2} = 0.$

\paragraph{Order $\varepsilon^{-1}$:} Matching terms of order $\varepsilon^{-1}$ gives
$$ V_{ij} \partial_z U^0_{k} =  \sigma_k \left[ \partial^2_{zz} U^{1}_k - W''(U^0_k)U^1_k  +   H_{ij}  \partial_z U^0_{k} \right]  +  \Lambda^{-1} \sqrt{2 W(U^0_k)} , \text{ for all} \quad k \in \{1,2, \cdots N\} $$
and
$$   \sum_{k=1}^{N}  \sqrt{2 W(U^0_k)} \Lambda^{-1}  = - \sum_{k=1}^{N}  \sigma_k \left[  \partial^2_{zz} U^1_k - W''(U^0_k) U^1_k + H_{ij} \partial_z U^0_{k} \right].$$
In particular,  for all $k \in \{1,2, \cdots N \}\setminus\{i,j\}$, we have
$\sigma_k \left( \partial^2_{zz} U^{1}_k - W''(0) U^1_k \right)=0$ which proves, by using the additional boundary conditions, that $U^{1}_k = 0.$
Moreover,   as $\sqrt{2 W(q)} = -q'$,  the equations for $k=i$ and $k=j$ read
$$
\begin{cases}
    V_{ij} q'(z) =  \sigma_i \left( \partial^2_{zz} U^{1}_i - W''(q(z)) U^1_i \right) +  \sigma_i q'(z) H_{ij}  -  \Lambda^{-1}(z,x,t) q'(z), \\
   -V_{ij} q'(z) =  \sigma_j \left( \partial^2_{zz} U^{1}_j - W''(q(z)) U^1_j \right) -  \sigma_j q'(z) H_{ij}  -  \Lambda^{-1}(z,x,t) q'(z),
\end{cases}
$$
while \eqref{def_lambda_allencahn_modif} gives
$$
 2 q'(z) \Lambda^{-1}(z,x,t) =  \sigma_i \left( \partial^2_{zz} U^{1}_i - W''(q(z)) U^1_i \right) +  \sigma_i q'(z)H_{ij}+
                             +  \sigma_j \left( \partial^2_{zz} U^{1}_j - W''(q(z)) U^1_j \right) -  \sigma_j q'(z)H_{ij}.
$$
As a consequence we obtain
$$
   V_{ij} q'(z) =  \frac{1}{2} (\sigma_i + \sigma_j) q'(z) H_{ij}  +  \frac{1}{2}\sigma_i \left( \partial^2_{zz} U^{1}_i - W''(q(z)) U^1_i \right) -
                                                                       \frac{1}{2}\sigma_j \left( \partial^2_{zz} U^{1}_j - W''(q(z)) U^1_j \right).
$$
Finally, multiplying this last equation by $q'$ and integrating over $\R$ leads to the interface law
$ V_{ij} =  \frac{1}{2} (\sigma_i + \sigma_j) H_{ij}.$
From the Fredholm alternative we deduce that
$ U^{1}_i = U^1_j = 0$ and $\Lambda^{-1}(z,x,t) = \frac{1}{2} \left( \sigma_i - \sigma_j\right)$
so that Claim~\ref{claim:allencahn_modif} is proved.

\subsection{Analysis of the Allen-Cahn system with harmonically additive mobilities}
We assume in this section that mobility coefficients are harmonically additive, {\em i.e.}, there exist coefficients
$m_i>0$ such that
$$ \frac{1}{m_{ij}} = \frac{1}{m_i} + \frac{1}{m_j},$$
and we consider the Allen-Cahn system
$$ \frac{1}{m_k} \partial_t u^{\epsilon}_k =  \sigma_k \left( \Delta u^{\epsilon}_k - \frac{1}{\epsilon^2} W'(u^{\epsilon}_k) \right) +  \lambda^{\epsilon}
\sqrt{2 W(u^{\epsilon}_k)},  \text{ for all } k \in \{1,2,\cdots,N\},$$
where
$$  \lambda^{\epsilon}  = \frac{\sum_{k} \left(   m_k \sigma_k  \left( \Delta u^{\epsilon}_k(x,t) - \frac{1}{\epsilon^2} W'(u^{\epsilon}_{k}(x,t)) \right)  \right)}{ \sum_{k} m_k \sqrt{2 W(u_k)}}$$
is the Lagrange multiplier field associated to the pointwise constraint $\sum_{k=1}^N{\boldsymbol u}_k = 1.$ The analysis below follows closely the matching conditions already used in both previous subsections:

\paragraph{Order $\varepsilon^{-2}$:} Identifying the terms of order $\varepsilon^{-2}$ gives for all  $k \in \{1,\cdots,N\}$:
$$\sigma_k \left( \partial^2_{zz} U^0_k - W'(U^0_k) \right)  + \Lambda^{-2} \sqrt{2 W(U^0_k)} = 0,$$
and
$$\left[ m_k \sum_{k=1}^{N}  \sqrt{2 W(U^0_k)} \right] \Lambda^{-2}  = - \sum_{k=1}^{N} m_k \sigma_k \left(  \partial^2_{zz} U^0_k - W'(U^0_k)  \right).$$
and leads to $U^0_i(z,x,t)  = q(z),$ $U^{0}_j(z,x,t) = q(-z),$ $U^{0}_k(z,x,t) = 0$ and $\Lambda^{-2} = 0.$

\paragraph{Order $\varepsilon^{-1}$:} Matching the next order terms shows that
$$ \frac{1}{m_k} V_{ij} \partial_z U^0_{k} =  \sigma_k \left[ \partial^2_{zz} U^{1}_k - W''(U^0_k)U^1_k  +   H_{ij}  \partial_z U^0_{k} \right]  +  \Lambda^{-1} \sqrt{2 W(U^0_k)}  $$
and
$$ \left[ m_k \sum_{k=1}^{N}  \sqrt{2 W(U^0_k)} \right] \Lambda^{-1}  = - \sum_{k=1}^{N}  m_k \sigma_k \left[  \partial^2_{zz} U^1_k - W''(U^0_k) U^1_k + H_{ij} \partial_z U^0_{k} \right].$$
Notice that for all $k \in \{1,2, \cdots N \}\setminus\{i,j\}$, we have $m_k \sigma_k\left( \partial^2_{zz} U^{1}_k - W''(0) U^1_k \right)=0$ from which, by using matching boundary conditions, we deduce that
$U^{1}_k = 0.$ Moreover, as $\sqrt{2 W(q)} = -q'$,  equations for $U_i^1$ and $U_j^1$ become
\begin{equation}
\begin{cases}
 \frac{1}{m_i} V_{ij} q'(z) &=  \sigma_i \left( \partial^2_{zz} U^{1}_i - W''(q(z)) U^1_i \right) +  \sigma_i q'(z) H_{ij}  -  \Lambda^{-1}(z,x,t) q'(z), \\
 \frac{1}{m_j} V_{ij} q'(z) &=  \sigma_j \left( \partial^2_{zz} U^{1}_j - W''(q(z)) U^1_j \right) -  \sigma_j q'(z) H_{ij}  -  \Lambda^{-1}(z,x,t) q'(z)
\end{cases}
\label{Ui_and_Uj_in_Allen-Cahn_add_mobility}
\end{equation}
and
\begin{equation}
\begin{split}
 \left( m_i +  m_j \right) q'(z) \Lambda^{-1}(z,x,t) &=  m_i \sigma_i  \left( \partial^2_{zz} U^{1}_i - W''(q(z)) U^1_i \right) +   m_i \sigma_i q'(z) H_{ij}+  \\
                                                     &+  m_j \sigma_j  \left( \partial^2_{zz} U^{1}_j - W''(q(z)) U^1_j \right) -   m_j \sigma_j q'(z) H_{ij}.
\end{split}
                                                     \label{lambda_AC_additive_mobility}
\end{equation}
Combining (\ref{Ui_and_Uj_in_Allen-Cahn_add_mobility}) and (\ref{lambda_AC_additive_mobility}) we obtain
\begin{eqnarray*}
\left(\frac{1}{m_i} + \frac{1}{m_j} \right)  V_{ij} q'(z) =  (\sigma_i + \sigma_j) q'(z) H_{ij}  +   \sigma_i \left( \partial^2_{zz} U^{1}_i - W''(q(z)) U^1_i \right) - \sigma_j \left( \partial^2_{zz} U^{1}_j - W''(q(z)) U^1_j \right).
\end{eqnarray*}
Multiplying this equation by $q'$ and integrating over $\R$ leads to the interface evolution
$$  \left(\frac{1}{m_i} + \frac{1}{m_j} \right)   V_{ij} =   (\sigma_i + \sigma_j) H_{ij}$$
or, equivalently $\frac{1}{m_{ij}}  V_{ij} =   \sigma_{ij} H_{ij}.$
Additionally, it appears that $ U^{1}_i = U^1_j = 0,$ and $\Lambda^{-1} = \frac{1}{2} \left( \sigma_i - \sigma_j\right),$ which finally proves Claim~\ref{claim:allencahn_mob_additive}.

\subsection{Analysis of the Allen-Cahn system  with general mobilities}
We now consider the more general situation described by the Allen-Cahn system~\eqref{eqn:allencahn_mob_gen}:
$$ - M \partial_t {\boldsymbol u}^{\epsilon}=   \sigma \Delta {\boldsymbol u}^{\epsilon} - \frac{1}{\epsilon^2} W'({\boldsymbol u}^{\epsilon})  + \lambda \sqrt{2 W({\boldsymbol u}^{\epsilon})},$$
where $M_{ij} = \frac{1}{m_{ij}}$ for $1\leq i\neq j\leq N,$ $M_{ii} = 0$ for $i\in\{1 ,2,\cdots N\},$ with Lagrange multiplier
$$ \lambda^{\epsilon}(x,t) = \frac{ \sum_{k=1}^{N} \left( M^{-1}\left(  \sigma \Delta {\boldsymbol u}^{\epsilon} - \frac{1}{\epsilon^2} W'({\boldsymbol u}^{\epsilon})\right) \right)_k}{ \sum_{k=1}^{N} \left(M^{-1} \sqrt{2 W({\boldsymbol u}^{\epsilon})} \right)_k}.$$

\paragraph{Order $\varepsilon^{-2}$: } For all $k \in \{1,2, \cdots,N\}$ we obtain $\sigma_k \left( \partial^2_{zz} U^0_k - W'(U^0_k) \right)  + \Lambda^{-2} \sqrt{2 W(U^0_k)} = 0,$
and
$$  \left[ \sum_{l=1}^{N}  \sum_{k=1}^{N} M^{-1}_{lk} \sqrt{2 W(U^0_k)} \right] \Lambda^{-2}  = - \sum_{l=1}^{N} \sum_{k=1}^{N}  M^{-1}_{lk}  \sigma_k \left(  \partial^2_{zz} U^0_k - W'(U^0_k)  \right).$$
Again, we can deduce that
$U^0_i(z,x,t)  = q(z)$ $U^{0}_j(z,x,t) = 1-q(z),$ $U^{0}_k(z,x,t) = 0$ and  $\Lambda^{-2} = 0.$

\paragraph{Order $\varepsilon^{-1}$:} Matching the next order terms shows that, for all $k \in \{1,\cdots,N\}$,
 $$ - V_{ij} \sum_{l=1}^{N} M_{k,l}  \partial_z U^0_{l} =  \sigma_k \left[ \partial^2_{zz} U^{1}_k - W''(U^0_k)U^1_k  +   H_{ij}  \partial_z U^0_{k} \right]  +  \Lambda^{-1} \sqrt{2 W(U^0_k)}.$$
In particular, if $k=i$ and $k=j$ we obtain
\begin{equation}
\begin{cases}
          V_{ij}  \frac{1}{m_{ij}} q'(z) = \sigma_i \left[ \partial^2_{zz} U^{1}_i - W''(q)U^1_i  +   H_{ij}  q'(z) \right]  -  \Lambda^{-1}(z,x,t) q'(z), \\
         -V_{ij}  \frac{1}{m_{ij}} q'(z) = \sigma_j \left[ \partial^2_{zz} U^{1}_j - W''(q)U^1_j  -   H_{ij}  q'(z) \right]  -  \Lambda^{-1}(z,x,t) q'(z),
    \end{cases}
    \label{V_ij_non_additive}
\end{equation}
while for all $k \in \{1,2 \cdots N\}\setminus\{i,j\}$
\begin{equation}
V_{ij} \left(  \frac{1}{m_{kj}} - \frac{1}{m_{ki}} \right) q'(z) =   \sigma_k \left[ \partial^2_{zz} U^{1}_k - W''(0)U^1_k \right].
\label{V_ij_non_additive_k_component}
\end{equation}
From (\ref{V_ij_non_additive}) and (\ref{V_ij_non_additive_k_component}) we deduce that
\begin{eqnarray*}
2 \frac{1}{m_{ij}} V_{ij}  q'(z) &=& \sigma_i \left[ \partial^2_{zz} U^{1}_i - W''(q)U^1_i  +   H_{ij}  q'(z) \right] \\
                                 &&- \sigma_j \left[ \partial^2_{zz} U^{1}_j - W''(q)U^1_j  -   H_{ij}  q'(z) \right] +  H_{ij} (\sigma_i + \sigma_j)  q'(z).
\end{eqnarray*}
Again, multiplying the expression above by $q',$ integrating over $\R$ and taking into account the conditions at $z\rightarrow\pm\infty$ leads to the interface evolution
$$\frac{1}{m_{ij}} V_{ij} =  \frac{1}{2} H_{ij} (\sigma_i + \sigma_j) = \frac{1}{2}  \sigma_{ij} H_{ij}.$$
Finally, we notice that as soon as $m_{ki} \neq m_{kj},$ $U^1_k$ can not vanish, since
$$V_{ij} \left(  \frac{1}{m_{kj}} - \frac{1}{m_{ki}} \right) q'(z) =   \sigma_k \left[ \partial^2_{zz} U^{1}_k - W''(0)U^1_k \right].$$
Thus, Claim~\ref{claim:allencahn_mob_gen} is proved.

\section{Numerical scheme and simulations}
We introduce in this section a numerical scheme to approximate solutions to the following systems:
\begin{itemize}
 \item the Allen-Cahn system with mobilities but without volume constraints :
 $$ \partial_t  u^{\epsilon}_k =   m_k \left[ \sigma_k \left( \Delta u^{\epsilon}_k - \frac{1}{\epsilon^2} W'(u^{\epsilon}_k) \right) +  \lambda^{\epsilon} \sqrt{2 W(u_k)} \right],$$
with $\sum_{k=1}^{N} u_k = 1$.
 \item the Allen-Cahn system with mobilities and volume contraints  :
 $$ \partial_t  u^{\epsilon}_k =   m_k \left[ \sigma_k \left( \Delta u^{\epsilon}_k - \frac{1}{\epsilon^2} W'(u^{\epsilon}_k) \right) +  \lambda^{\epsilon} \sqrt{2 W(u_k)}  + \mu_k(t) G_k({\boldsymbol u }) \right],$$
with $\sum_{k=1}^{N} u_k = 1$ and $ \int_{Q} u_k dx =  {\textrm Vol}_k(t).$
\end{itemize}

The multiphase field model for VLS growth can be regarded as an Allen-Cahn system with additional volume contraints for the three phases $ {\boldsymbol u }= (u_S,u_L,u_V).$ In this special case the potentials $G_k$ are defined by
$$ G_L({\boldsymbol u}) = \sqrt{2 W(u_L)},  \quad G_S({\boldsymbol u}) =  u_S u_L,\quad \text{ and }    \quad G_V({\boldsymbol u}) =  u_V u_L.$$
We consider the solution to the Allen-Cahn system for times $t \in [0,T]$, in a computation box $Q$ with periodic boundary conditions and with the initial
condition ${\boldsymbol u}(x,0) = {\boldsymbol u}^0$. We also assume that ${\boldsymbol u}^0$ satisfies
the partition constraint $\sum_{k=1}^{N} u^0_k = 1$.

We propose a standard Fourier spectral splitting scheme \cite{Chen_fourier} to compute numerically the solution
to the above Allen-Cahn systems. We recall that the Fourier $K$-approximation of a function $u$ defined in a box $Q = [0,L_1]\times \cdots \times [0,L_d]$ is given by
$$u^{K}(x) = \sum_{{\boldsymbol k}\in{\llbracket -\frac K 2,\frac K 2-1\rrbracket}^d} c_{\boldsymbol k} e^{2i\pi{\boldsymbol\xi}_k\cdot x},$$
where ${\boldsymbol k} = (k_1,\dots,k_d)$ and ${\boldsymbol \xi_k} = (k_1/L_1,\dots,k_d/L_d)$. In this formula, $c_{\boldsymbol k}$ runs over the $K^d$ first discrete Fourier coefficients of $u$. The inverse discrete Fourier transform of $c_{\boldsymbol k}$ leads to $u^{K}_{\boldsymbol k} =   \textrm{IFFT}[c_{\boldsymbol k}]$ where $u^{K}_{\boldsymbol k}$
denotes the value of $u$ at the points $x_{\boldsymbol k} = (k_1 h_1, \cdots, k_d h_d)$  where $h_{\alpha} = L_{\alpha}/K$ for $\alpha\in\{1,\cdots,d\}$. Conversely,
$c_{\boldsymbol k}$ can be computed as the discrete Fourier transform of $u^K_{\boldsymbol k},$ {\em i.e.}, $c_{\boldsymbol k} = \textrm{FFT}[u^K_{\boldsymbol k}].$
\subsection{Scheme overview}
We now introduce a time discrete sequence ${\boldsymbol u}^n$ for the approximation of ${\boldsymbol u}$ at times  $n \delta_t$. This sequence is defined as follows for each problem under study:
\paragraph{Allen-Cahn system without volume constraints:}
\begin{itemize}
\item[\textbf{Step $1$:}] $L^2$-gradient flow of the Cahn-Hilliard energy without constraints: let ${\boldsymbol u}^{n+1/2}$ be an approximation of ${\boldsymbol v}(\delta_t)$ where ${\boldsymbol v} = (v_1,\dots, v_N)$ is the solution to:
$$
\begin{cases}
   \partial_t v_k(x,t) &= m_k \sigma_k \left[\Delta v_k (x,t) - \frac{1}{\varepsilon^2} W'(v_k(x,t))\right] \quad \forall (x,t) \in Q \times [0,\delta_t] \text{,}\\
  {\boldsymbol v}(x,0) &= {\boldsymbol u}^{n}(x), \quad \forall x \in Q \text{ with periodic boundary conditions.}
\end{cases}
$$
\item[\textbf{Step $2$:}] Projection onto the partition constraint: for all $k\in \{ 1,2, \dots, N \}$ define $u^{n+1}_k$ by
$$ u_k^{n+1} =  u^{n+1/2}_k + m_k  \lambda^{n+1} \sqrt{2 W(u^{n+1/2}_k)}$$
where
$$\lambda^{n+1} = \frac{1 - \sum_{k}  u^{n+1/2}_k}{\sum_k  m_k \sqrt{2 W(u^{n+1/2}_k)} } $$
\end{itemize}
\paragraph{Allen-Cahn system with volume constraints:}
\begin{itemize}
\item[\textbf{Step $1$:}] $L^2$-gradient flow of the Cahn-Hilliard energy without constraint: let ${\boldsymbol u}^{n+1/2}$ be an approximation of ${\boldsymbol v}(\delta_t)$ where ${\boldsymbol v} = (v_1,v_2, \dots, v_N)$ is the solution of
$$
\begin{cases}
   \partial_t v_k(x,t) &= m_k \sigma_k \left[\Delta v_k (x,t) - \frac{1}{\varepsilon^2} W'(v_k(x,t))\right] \quad \forall (x,t) \in Q \times [0,\delta_t] \text{,}\\
  {\boldsymbol v}(x,0) &= {\boldsymbol u}^{n}(x), \quad \forall x \in Q \text{ with periodic boundary conditions.}
\end{cases}
$$
\item[\textbf{Step $2$:}] Projection onto both the partition and volume constraints: for all $k\in \{ 1,2, \dots, N \}$ define $u^{n+1}_k$ as
$$ u_k^{n+1} =  u^{n+1/2}_k + m_k \lambda^{n+1} \sqrt{2 W(u^{n+1/2}_k)} + m_k \mu_k^{n+1}  G_k({\boldsymbol u}^{n+1/2} ),$$
where $\lambda^{n+1}$ and $\mu_i^{n+1}$ encode the discrete constraints
$\sum_{k=1}^{N} u_k^{n+1} = 1$ and $\int_Q  u_k^{n+1}  =  V_k^{n+1} = {\textrm Vol}_k((n+1)\delta_t).$
\end{itemize}

\subsection{Solving step one with a semi-implicit Fourier spectral scheme}
To compute ${\boldsymbol u}^{n+1/2}$ we use a semi-implicit numerical method. To this end, we consider the equation
$$ \left( I_d - \delta_t m_k \sigma_k \left( \Delta  -  \alpha/\varepsilon^2 I_d \right)   \right) u_k^{n+1/2} =     u_k^{n}  - \frac{\delta_t m_k \sigma_k}{\varepsilon^2} \left( W'(u_k^n) - \alpha u_k^{n}\right),$$
where $\alpha$ is a positive stabilization parameter. It is known that the Cahn-Hilliard energy decreases unconditionally \cite{EDTK_Eyre,ShenWWW12} as soon as
the explicit part, {\em i.e.}, $s \to W'(s) - \alpha s$ is the derivative of a concave function. For $W(s) = \frac{1}{2} s^2 (1-s)^2$ this is true for $\alpha>2.$
We notice also that even without the stabilization parameter ({\em i.e.}, for $\alpha = 0$) the semi-implicit
scheme is still stable under the classical condition $\delta_t \leq \frac{C}{\varepsilon^2}$ where $C = \sum_{s \in [0,1]} |W''(s)|.$
Finally, recall that the equation is computed in $Q$ with periodic boundary conditions, and, consequently, the inverse of the operator
$\left( I_d - m_k \delta_k\delta_t \left( \Delta  -  \alpha/\varepsilon^2 I_d \right) \right) $ can be easily computed in Fourier space \cite{Chen_fourier}  using Fast Fourier Transform.

\subsection{Solving step two in the case of volume constraints}
The case without volume constraints can be of course easily deduced from what follows. We look for solutions $\lambda^{n+1} $ and $\mu_k^{n+1}$ such that
$ \sum_{i=1}^{N} u_k^{n+1} = 1$ and $\int_Q  u_k^{n+1}  = V_k^{n+1},$
and where, for all $i\in\{1,2,\ldots,N\},$
$$u_i^{n+1} =  u^{n+1/2}_i + \lambda^{n+1} m_i \sqrt{2 W(u^{n+1/2}_i)} +  \mu_i^{n+1} m_i G_i({\boldsymbol u^{n+1/2}}).$$

Let us introduce $\overline{\lambda}^{n+1}_i = \int_Q  m_i \sqrt{2 W(u^{n+1/2}_i)} \lambda^{n+1} dx$, and notice that integration over $Q$ of the above equation leads to
$$\mu_i^{n+1} =  \frac{ \left[ V_i^{n+1} - \int_Q u^{n+1/2}_i dx \right] -  \overline{\lambda}^{n+1}_i }{ \int_Q m_i G_i({\boldsymbol u^{n+1/2}}) dx},$$
and
$$ \lambda^{n+1} =   \frac{ \left[1 - \sum_{k} u^{n+1/2}_k \right]  -  \sum_k \mu_{k}^{n+1} m_k G_k({\boldsymbol u^{n+1/2}})  }{\sum_{k} m_k \sqrt{2 W(u^{n+1/2}_k) }}.$$
 Now, the coefficients  $\overline{\lambda}^{n+1}_i$ satisfy
\begin{eqnarray*}
\overline{\lambda}^{n+1}_i &=& \int_Q \lambda^{n+1} m_i \sqrt{2 W(u_i^{n+1/2})} dx =\\
                           &=&  \int_Q  \frac{   m_i \sqrt{2 W(u_i^{n+1/2})} \left[1 - \sum_{k} u^{n+1/2}_k \right]}{\sum_{k} m_k \sqrt{2 W(u^{n+1/2}_k) }} dx      - \int_Q  \frac{\sum_k \mu_{k}^{n+1}  m_k G_k({\boldsymbol u^{n+1/2}})  m_i \sqrt{2 W(u_i^{n+1/2})}  }{\sum_{k} m_k\sqrt{2 W(u^{n+1/2}_k) }}dx =\\
                           &=& \sum_k \left[ \int_Q  \frac{ m_i m_k G_k({\boldsymbol u^{n+1/2}})    \sqrt{2 W(u_i^{n+1/2})}/\int_Q m_k G_k({\boldsymbol u^{n+1/2}})   dx }{\sum_{j} m_j \sqrt{2 W(u^{n+1/2}_j) }}dx \right]\overline{\lambda}^{n+1}_k =\\
                           &+& \int_Q \left( \left[  \left[1 - \sum_{k} u^{n+1/2}_k \right] - \sum_{k} \frac{  \left( V_k^{n+1} - \int_Q u_k^{n+1/2} dx\right)  m_k G_k({\boldsymbol u^{n+1/2}})     }{ \int_Q m_k G_k({\boldsymbol u^{n+1/2}})   dx  }   \right] \frac{ m_i \sqrt{2 W(u^{n+1/2}_i)}}{ \sum_{k=1}^{N} m_k \sqrt{2 W(u^{n+1/2}_k)}} \right) dx.
\end{eqnarray*}
In particular, this implies that $\overline{{\boldsymbol \lambda}} = (\overline{\lambda}_1,\overline{\lambda}_2, \dots, \overline{\lambda}_N)$ solves the linear system
$$ (I_d - A) \overline{\lambda} = b,$$
where
$$A_{ik} =  \left[\int_Q \frac{m_i \sqrt{2 W(u^{n+1/2}_i)} m_k G_k({\boldsymbol u^{n+1/2}})  /\int_Q m_k G_k({\boldsymbol u^{n+1/2}})  dx     }{ \sum_{j=1}^{N}  m_j \sqrt{2 W(u^{n+1/2}_j)}} dx \right]$$
and
$$b_i =   \int_Q \left( \left[  \left[1 - \sum_{k} u^{n+1/2}_k \right] - \sum_{k} \frac{  \left( V_k^{n+1} - \int_Q  u_k^{n+1/2} dx\right) m_k G_k({\boldsymbol u^{n+1/2}})     }{ \int_Q m_k G_k({\boldsymbol u^{n+1/2}})  dx  }   \right] \frac{m_i \sqrt{2 W(u^{n+1/2}_i)}}{ \sum_{k=1}^{N}  m_k \sqrt{2 W(u^{n+1/2}_k)}} \right) dx.$$
Notice that since $\sum_{i} A_{ki} = 1$ and the matrix $(Id - A)$ is not invertible. Otherwise, it is not difficult to see that  $\sum b_i = 0 $ as  $\sum m_i^{n+1} = |Q|$. This means
that the linear system $(I_d - A) \overline{\lambda} = b$ admits at least one solution and, for stability reasons, we assume in addition that $\sum_i \overline{\lambda}_{i} = 0$.

\subsection{Validation of the approach for highly contrasted mobilities}
\begin{figure}[htbp]
\centering
	\includegraphics[width=3.7cm]{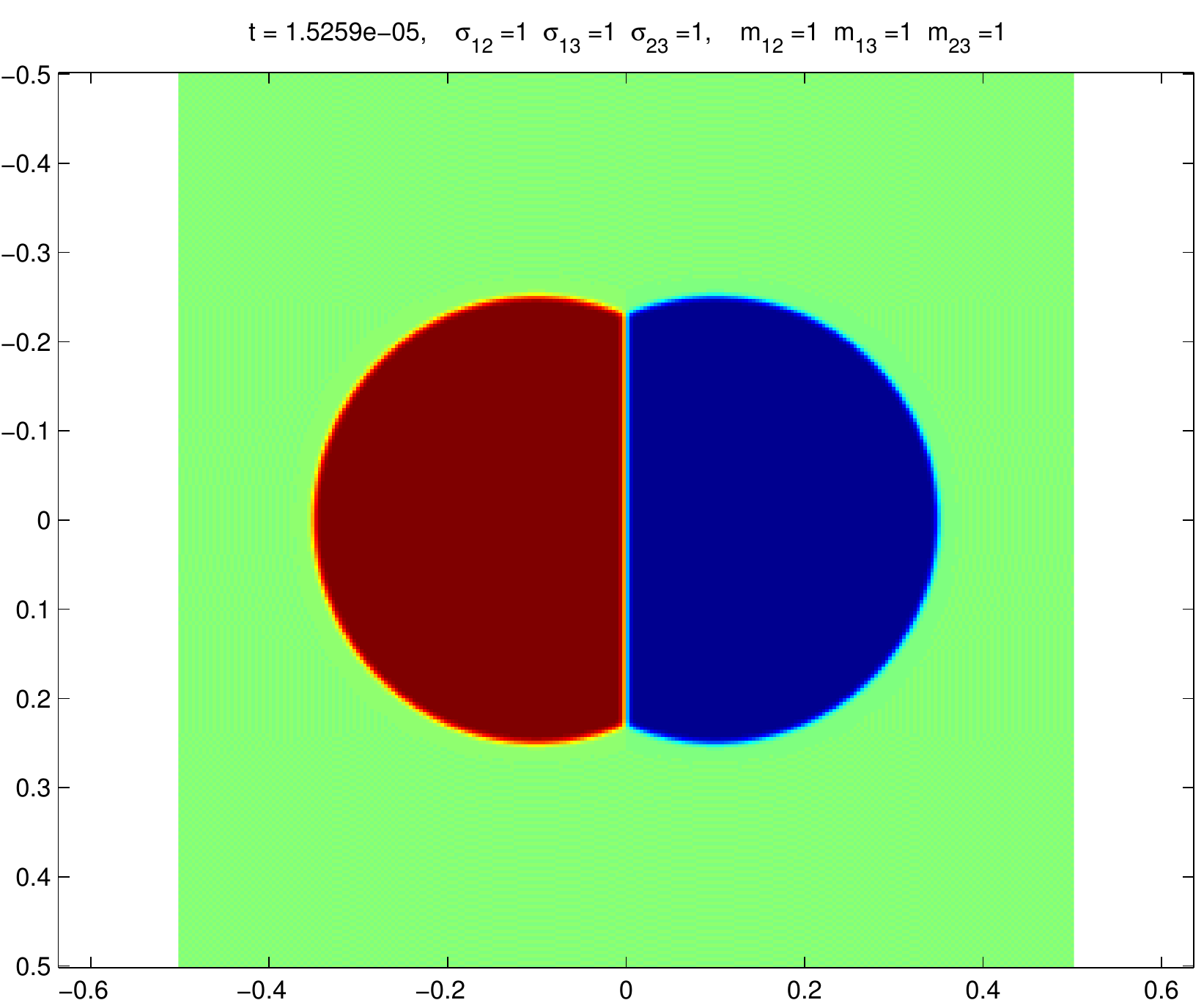}
        \includegraphics[width=3.7cm]{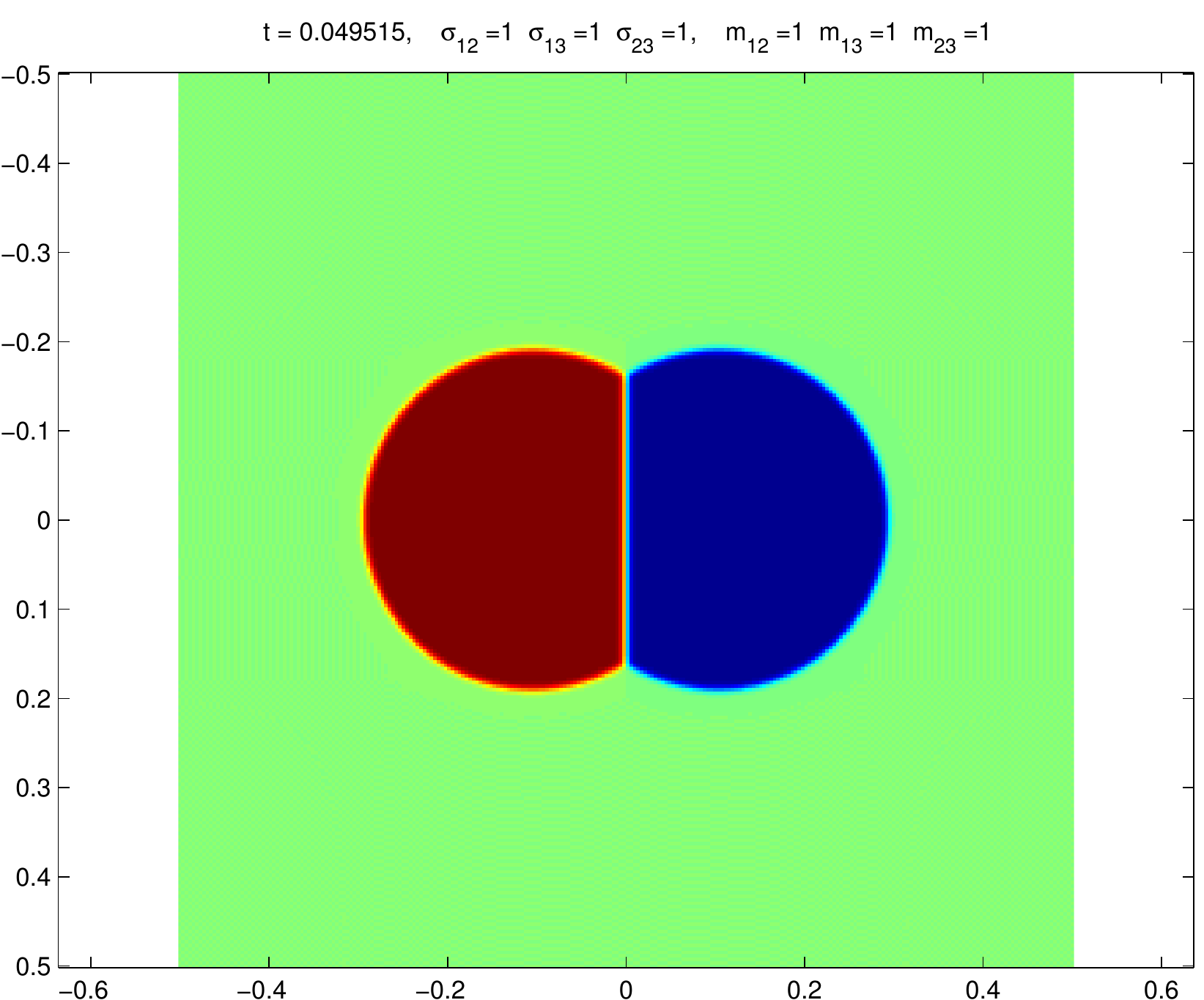}
        \includegraphics[width=3.7cm]{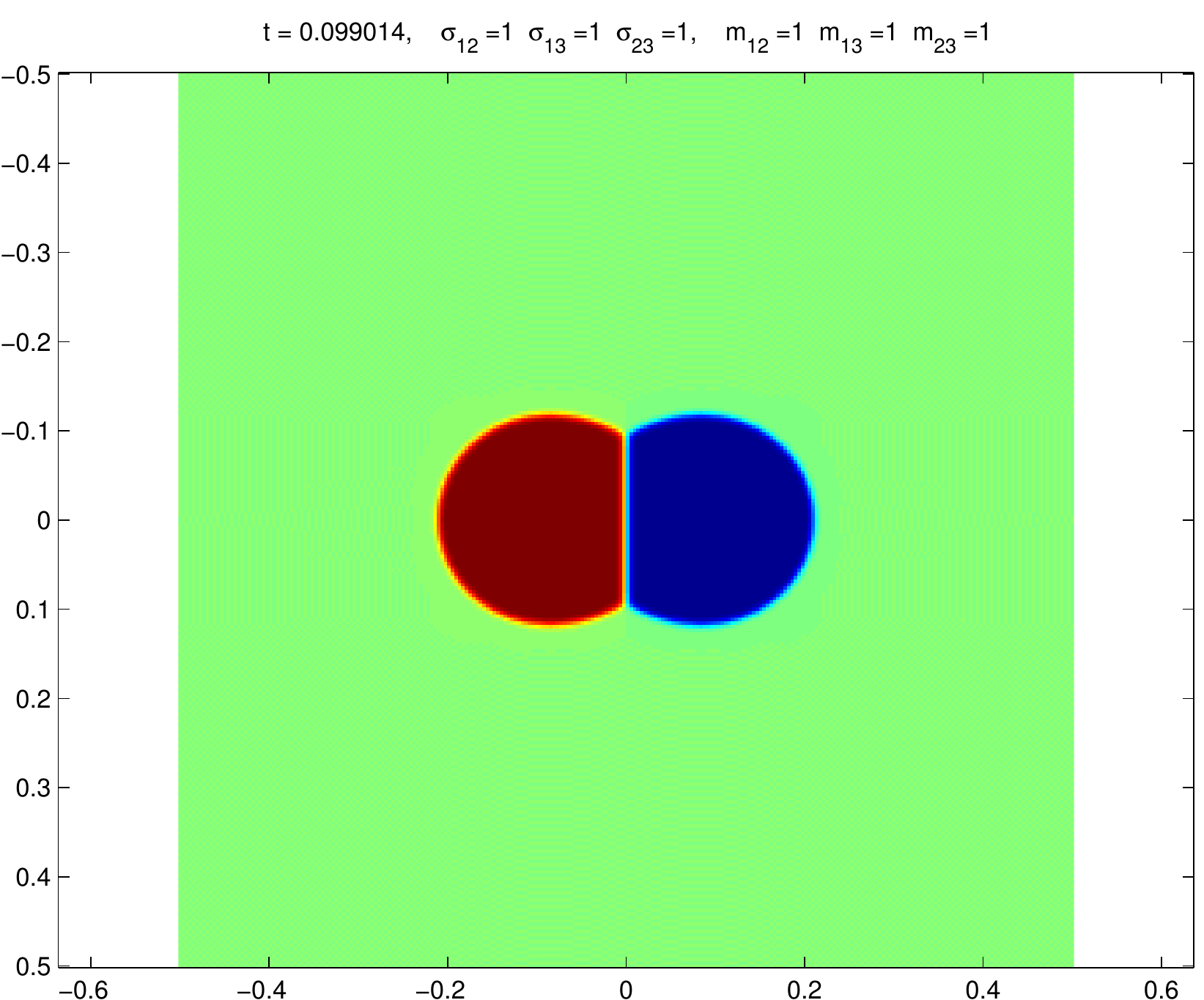}
        \includegraphics[width=3.7cm]{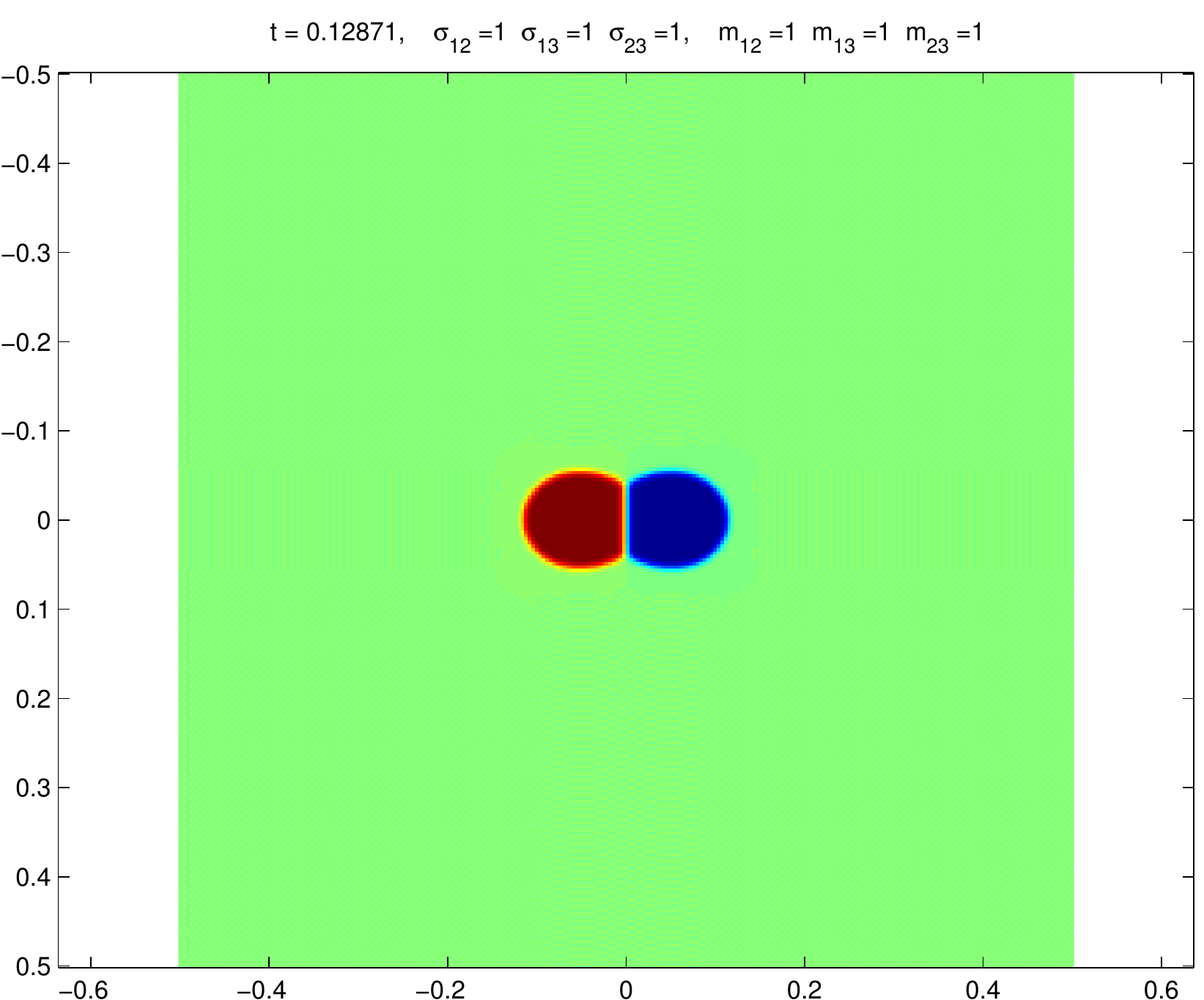} \\
        \includegraphics[width=3.7cm]{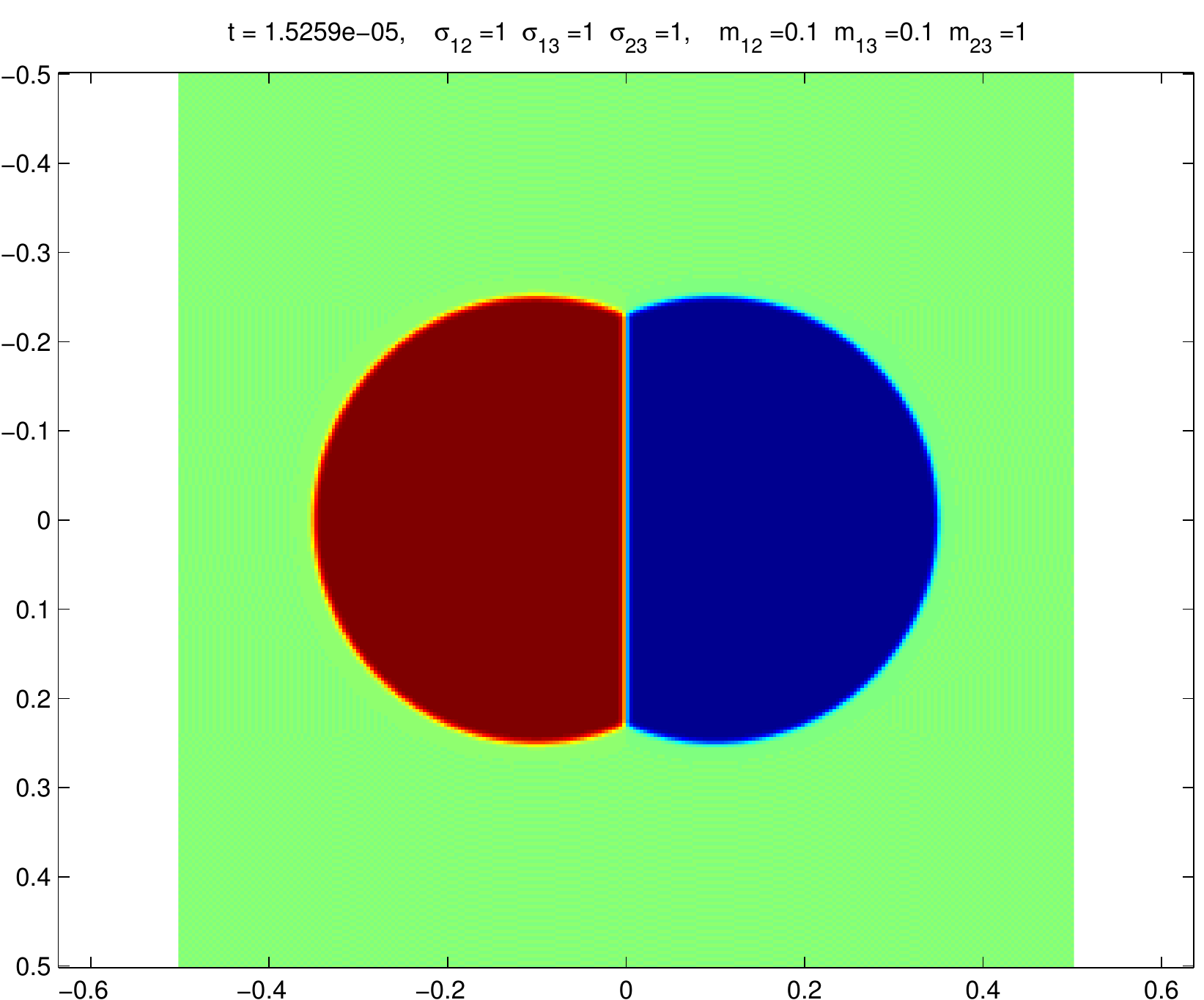}
        \includegraphics[width=3.7cm]{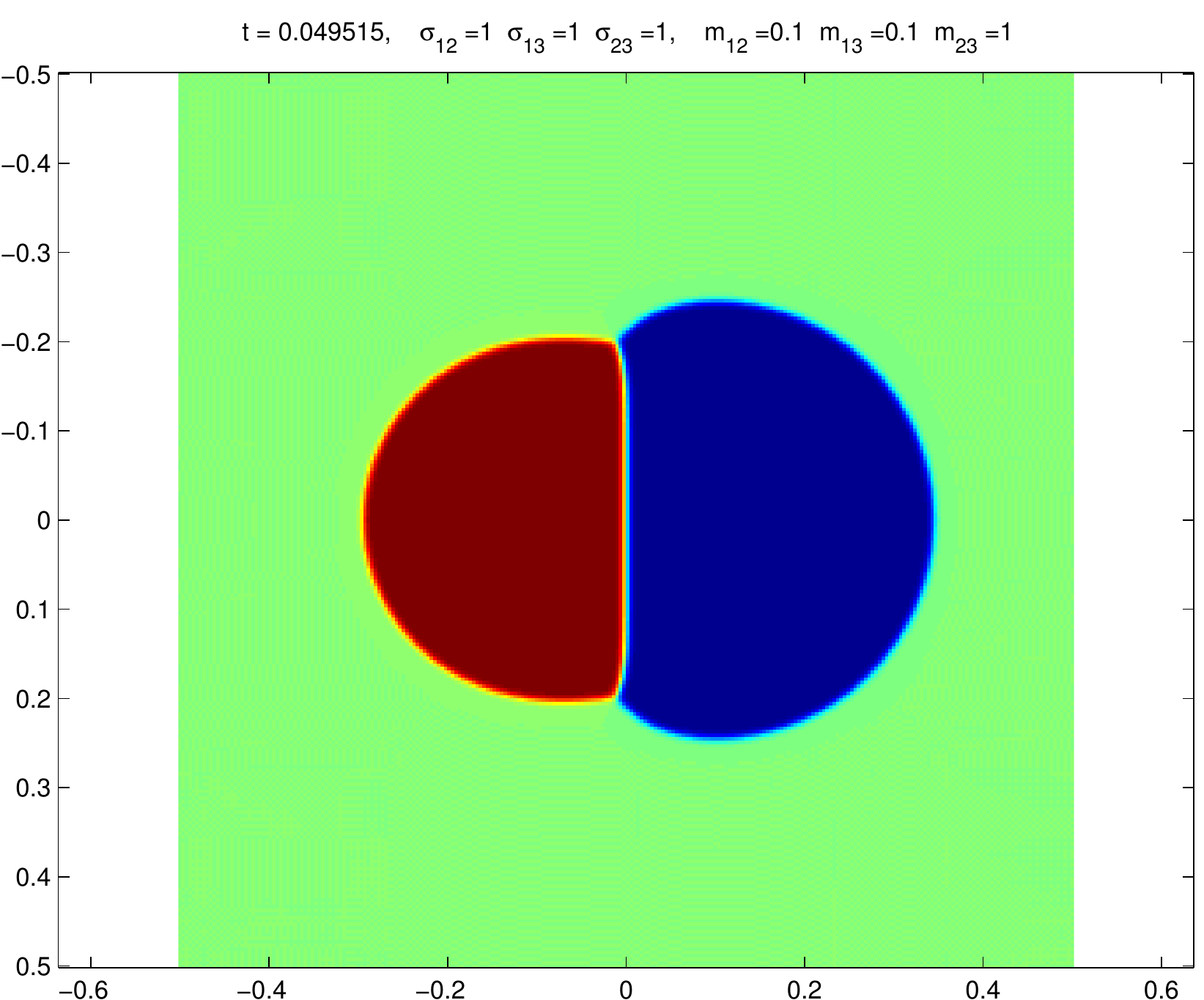}
        \includegraphics[width=3.7cm]{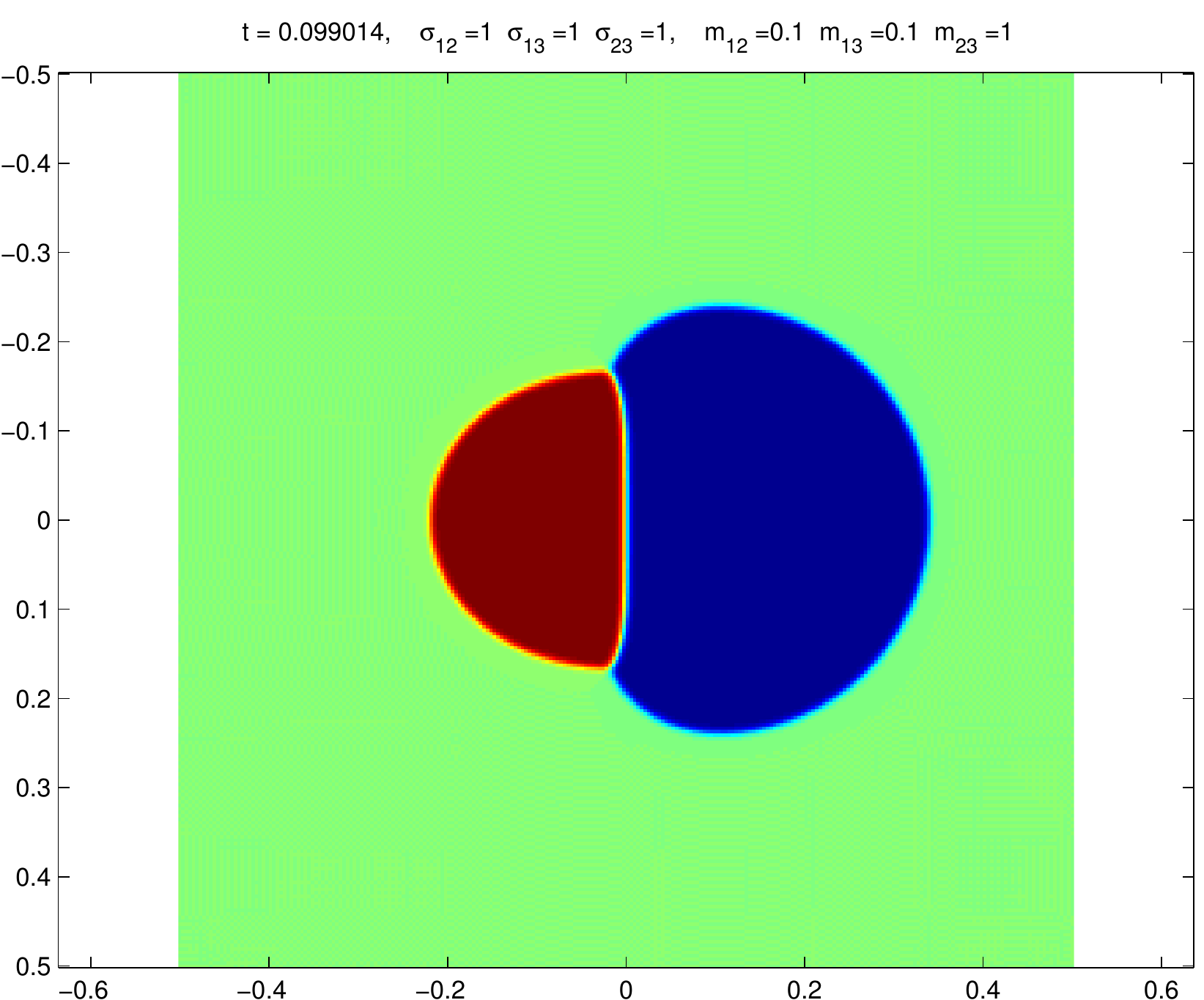}
        \includegraphics[width=3.7cm]{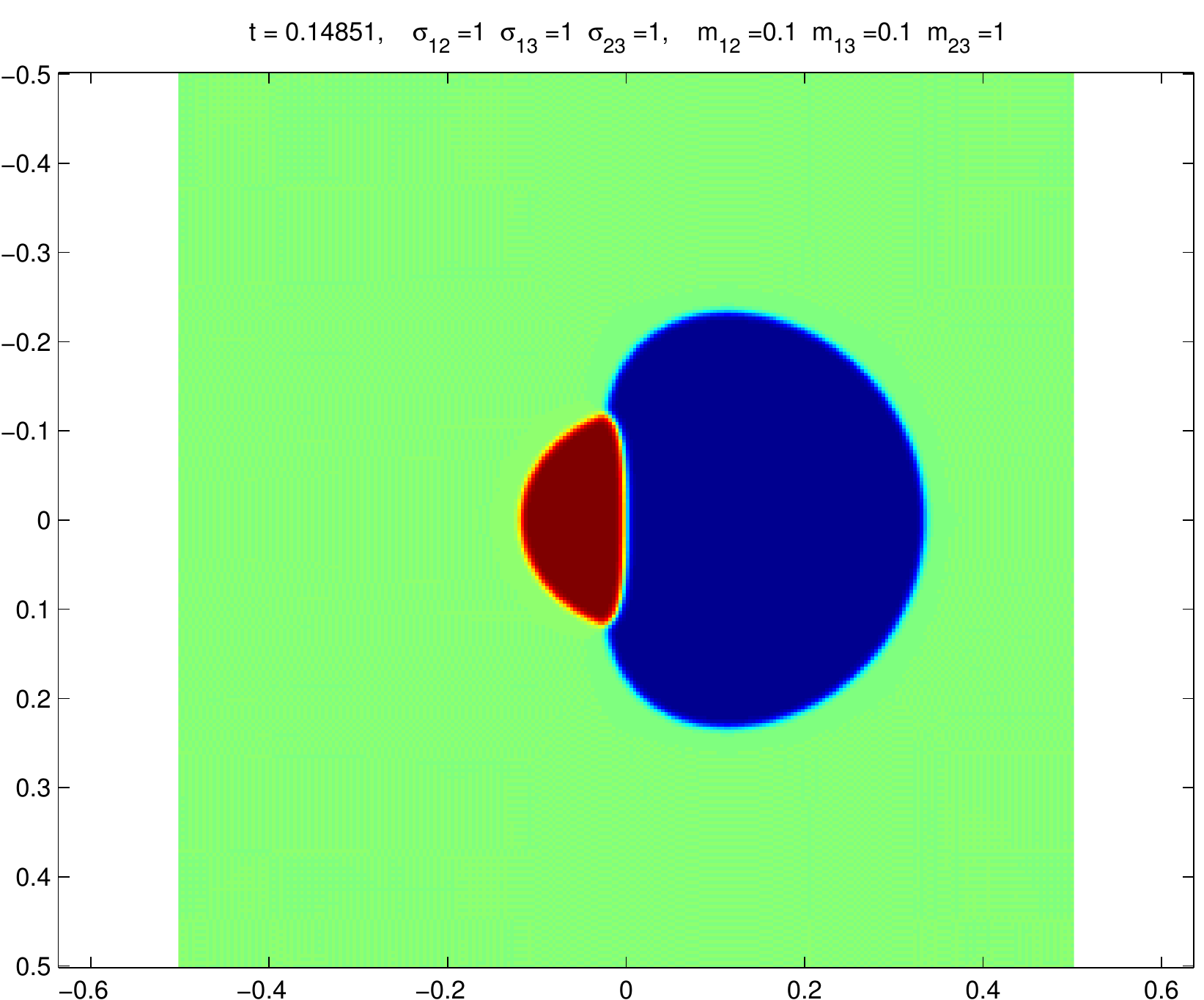} \\
        \includegraphics[width=3.7cm]{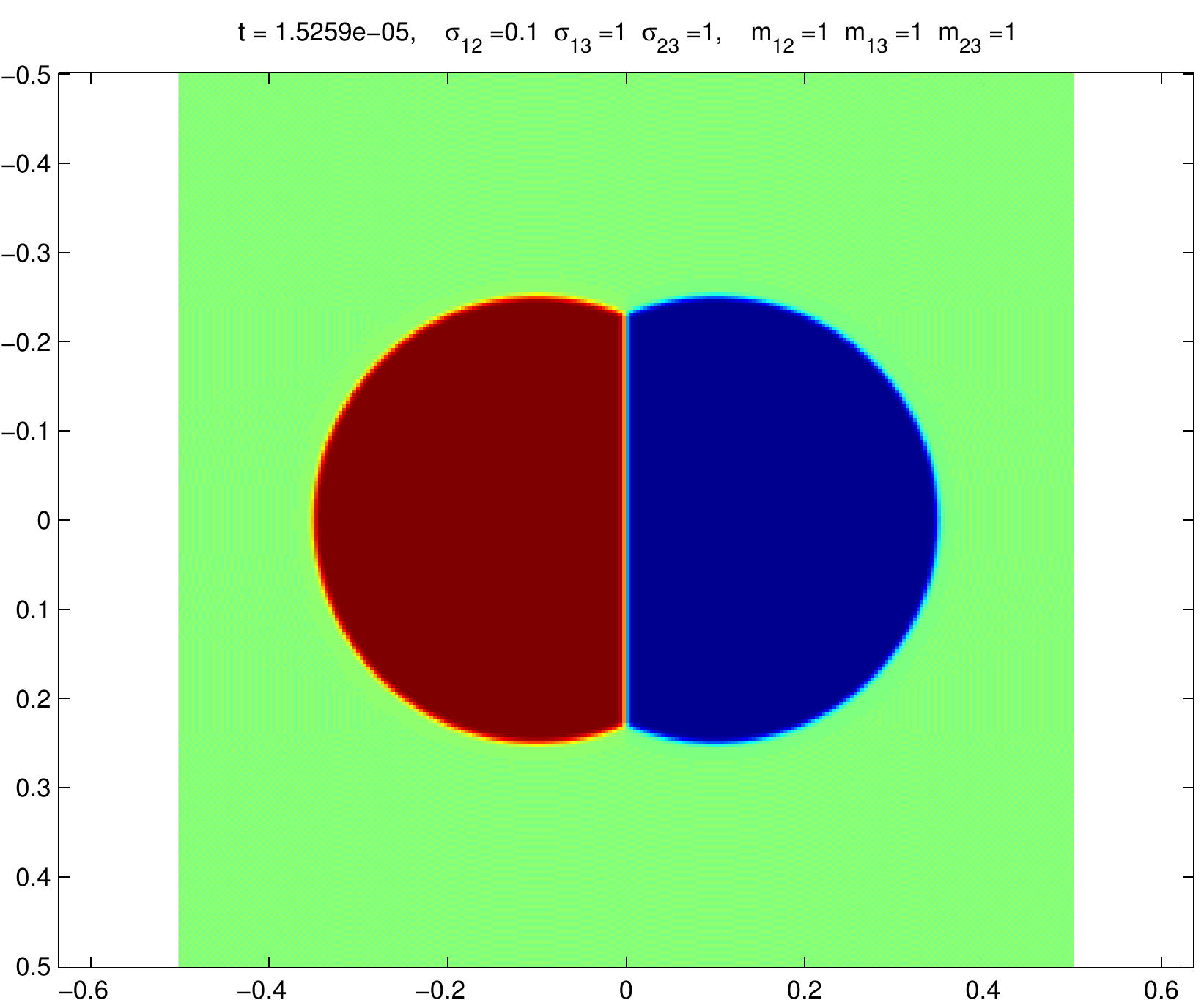}
        \includegraphics[width=3.7cm]{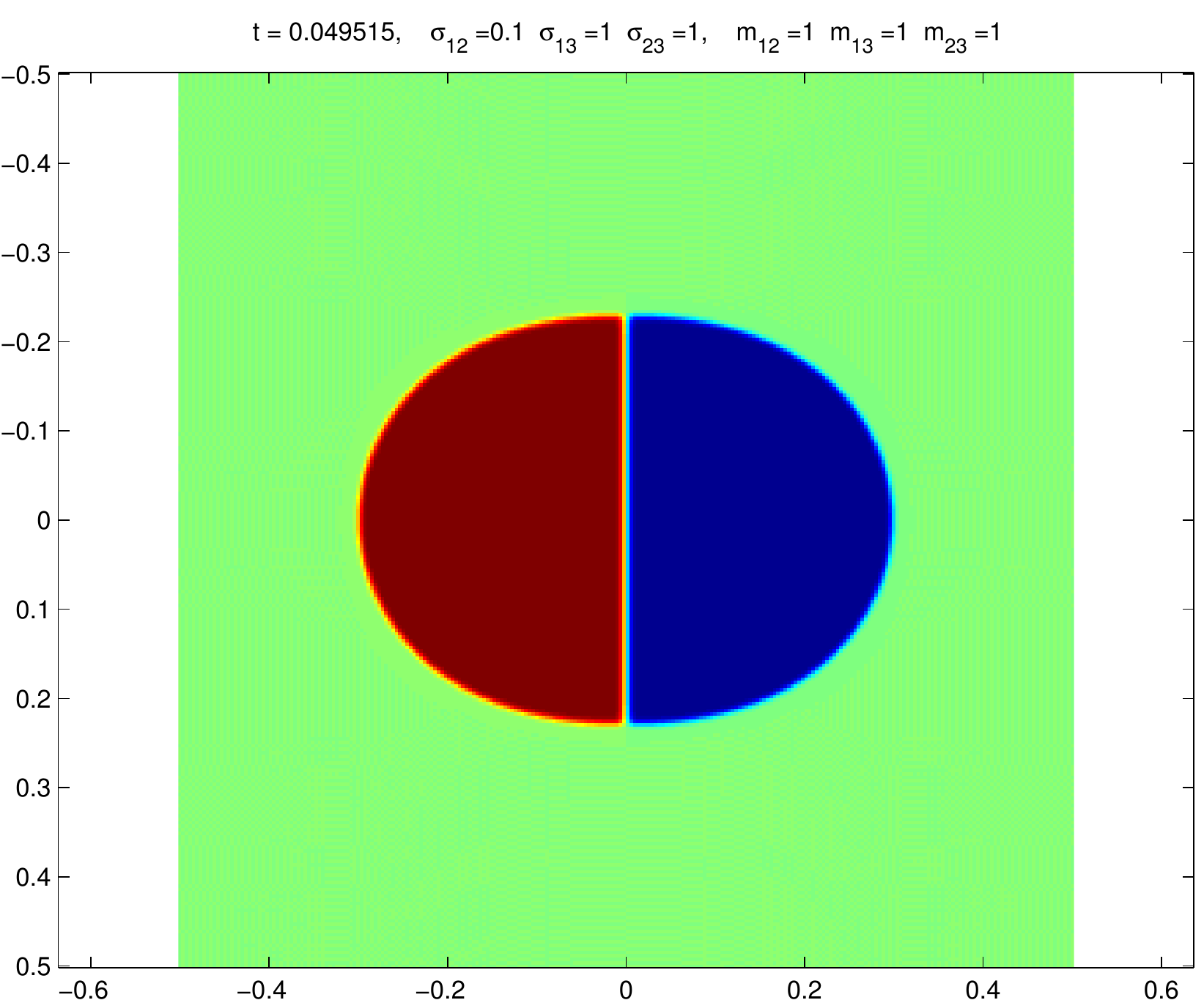}
        \includegraphics[width=3.7cm]{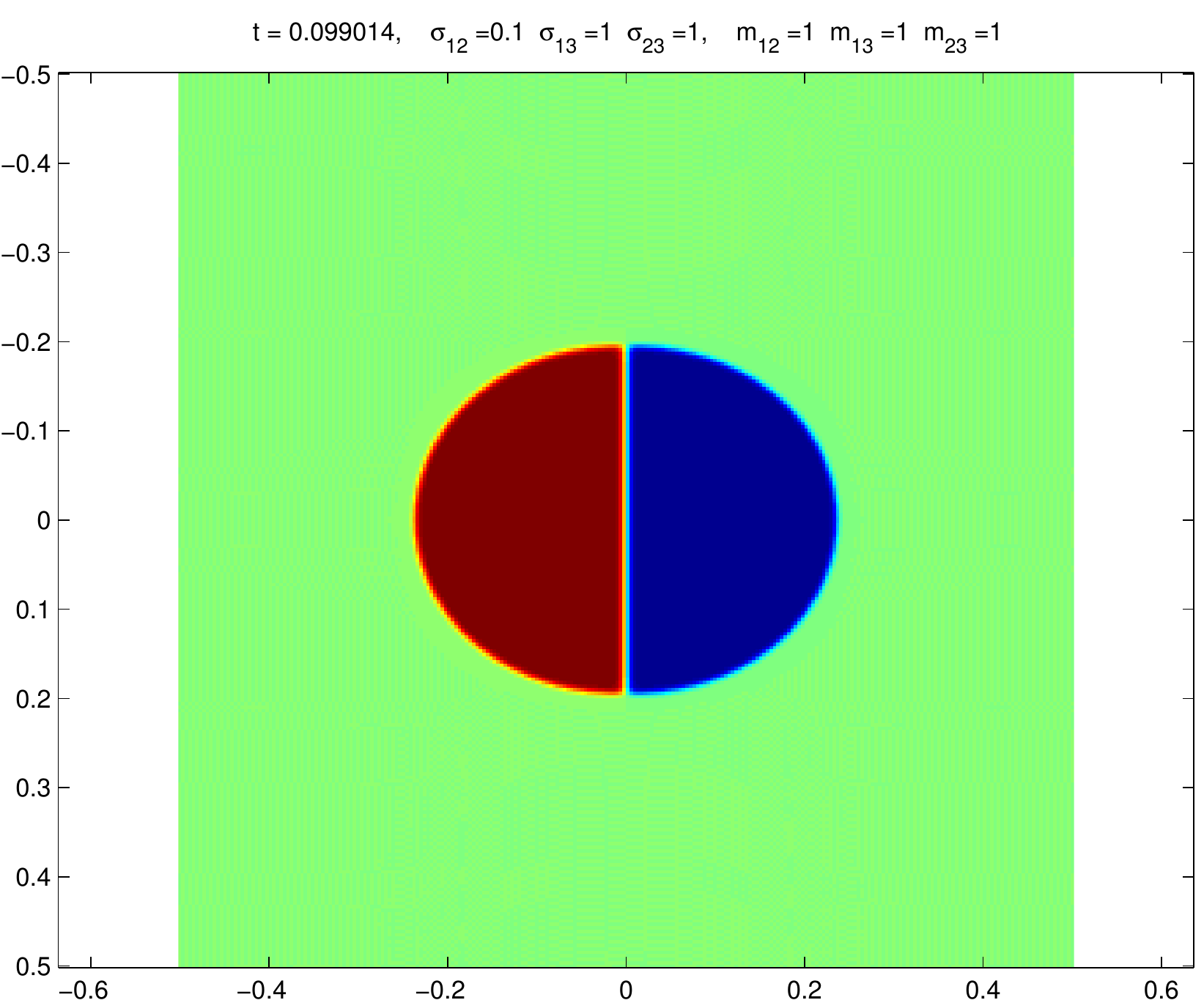}
        \includegraphics[width=3.7cm]{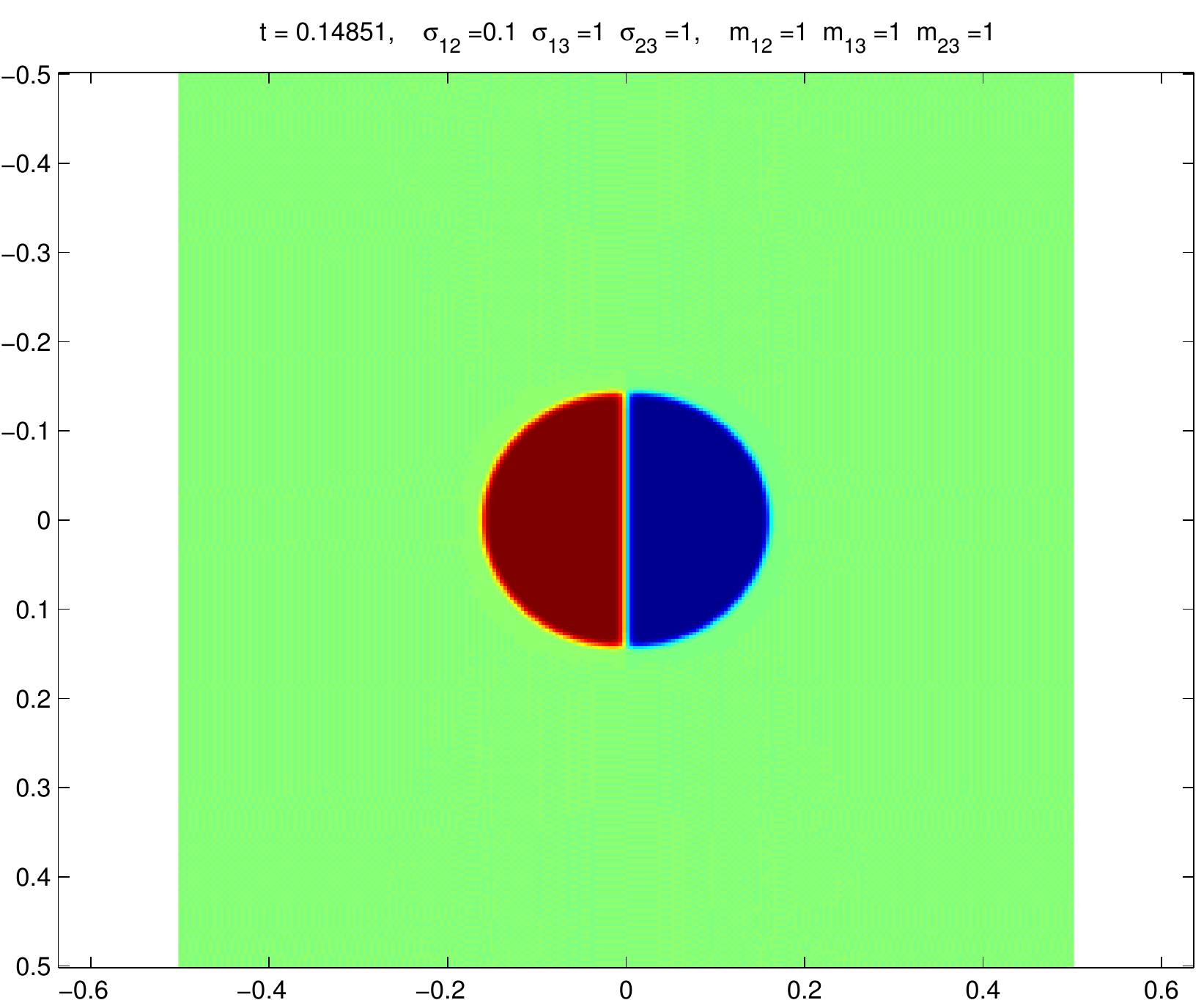} \\
        \includegraphics[width=3.7cm]{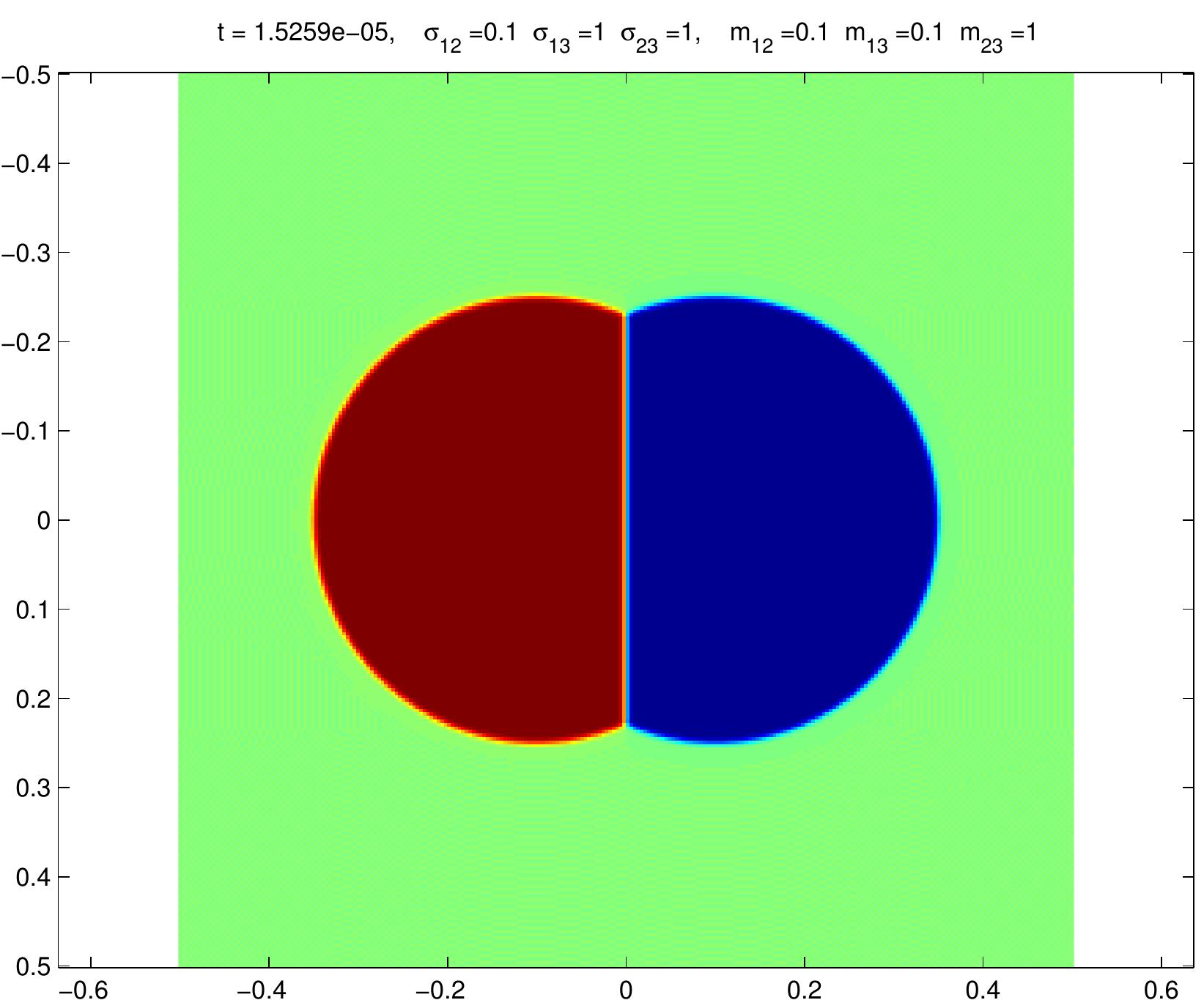}
        \includegraphics[width=3.7cm]{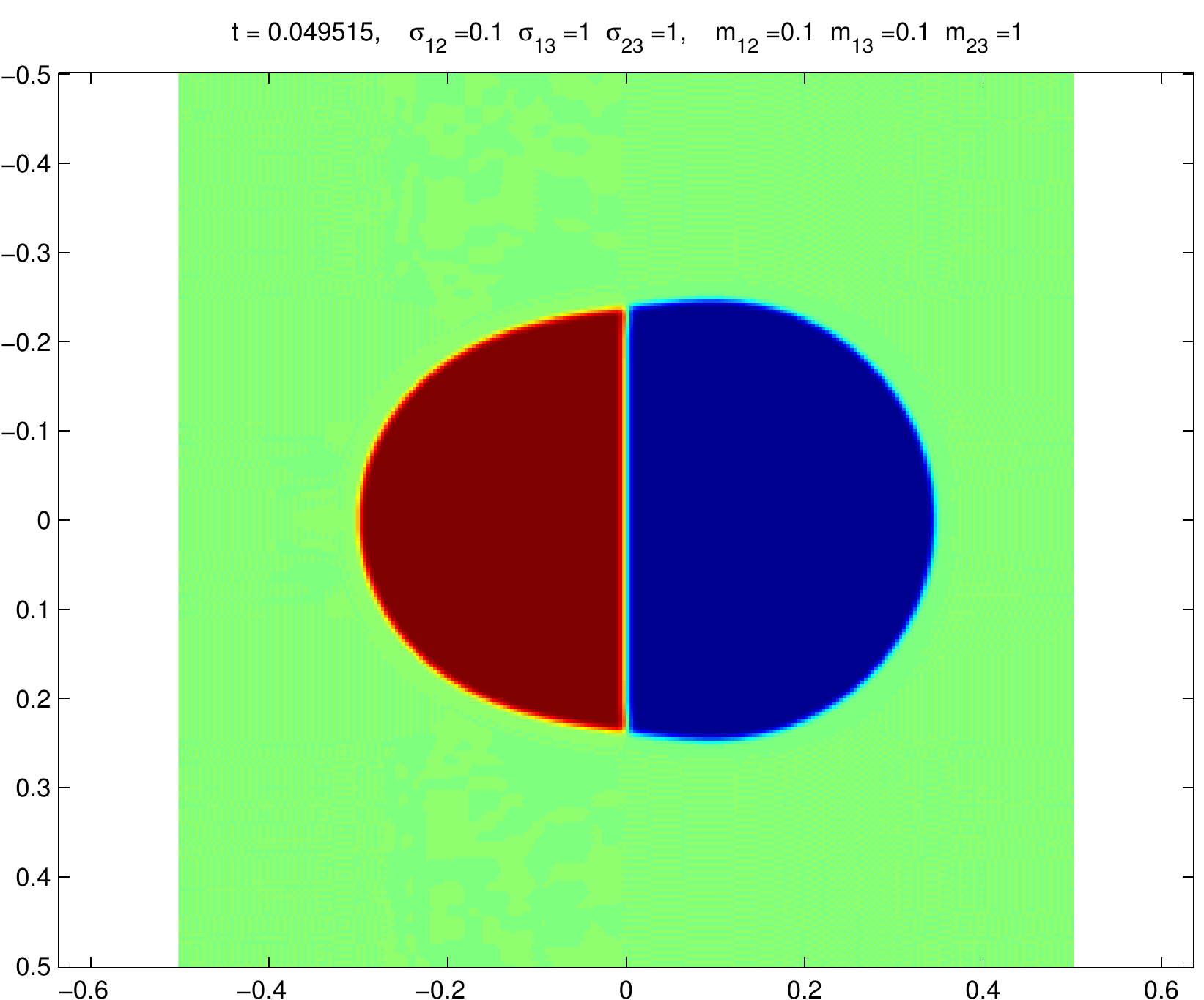}
        \includegraphics[width=3.7cm]{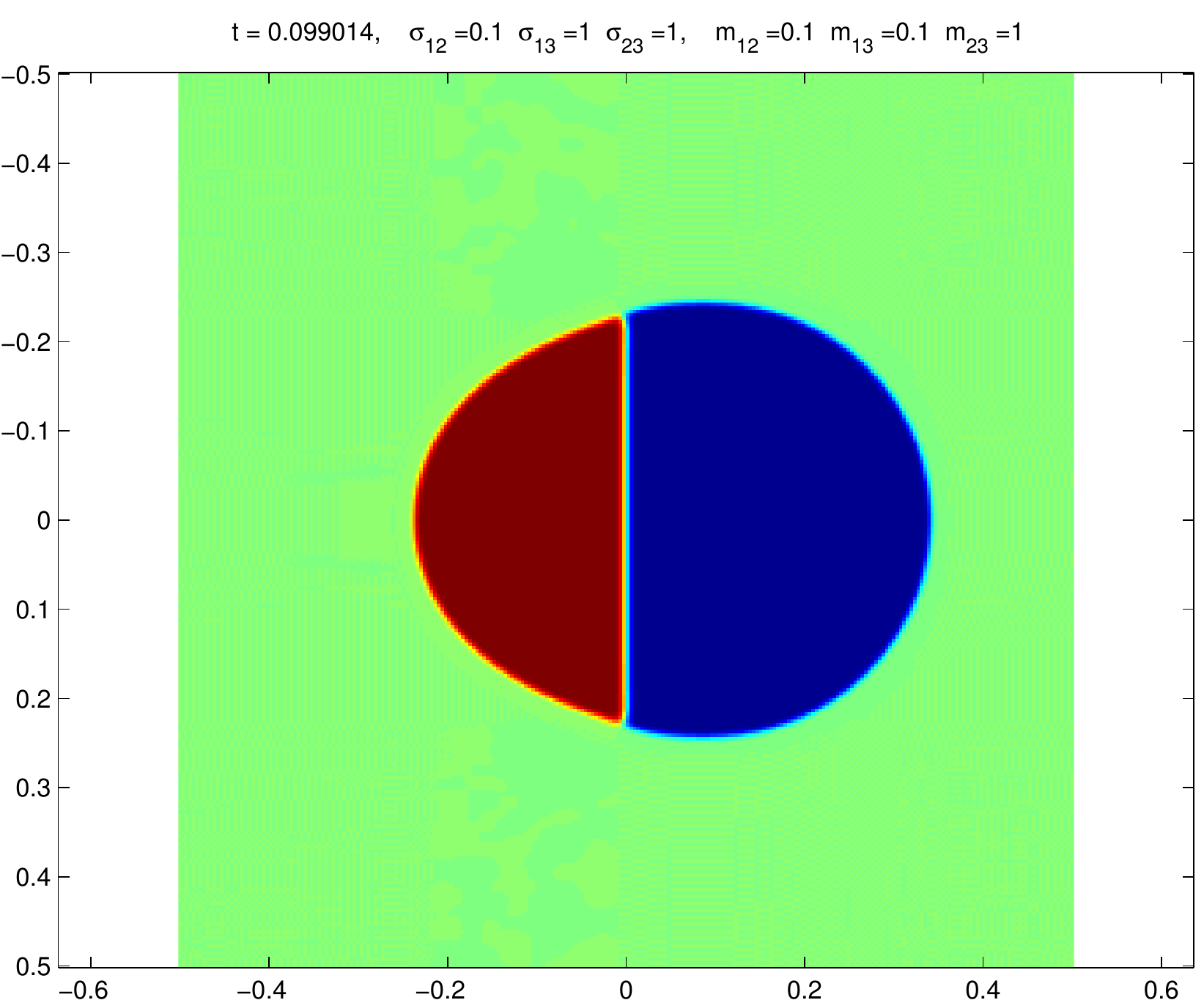}
        \includegraphics[width=3.7cm]{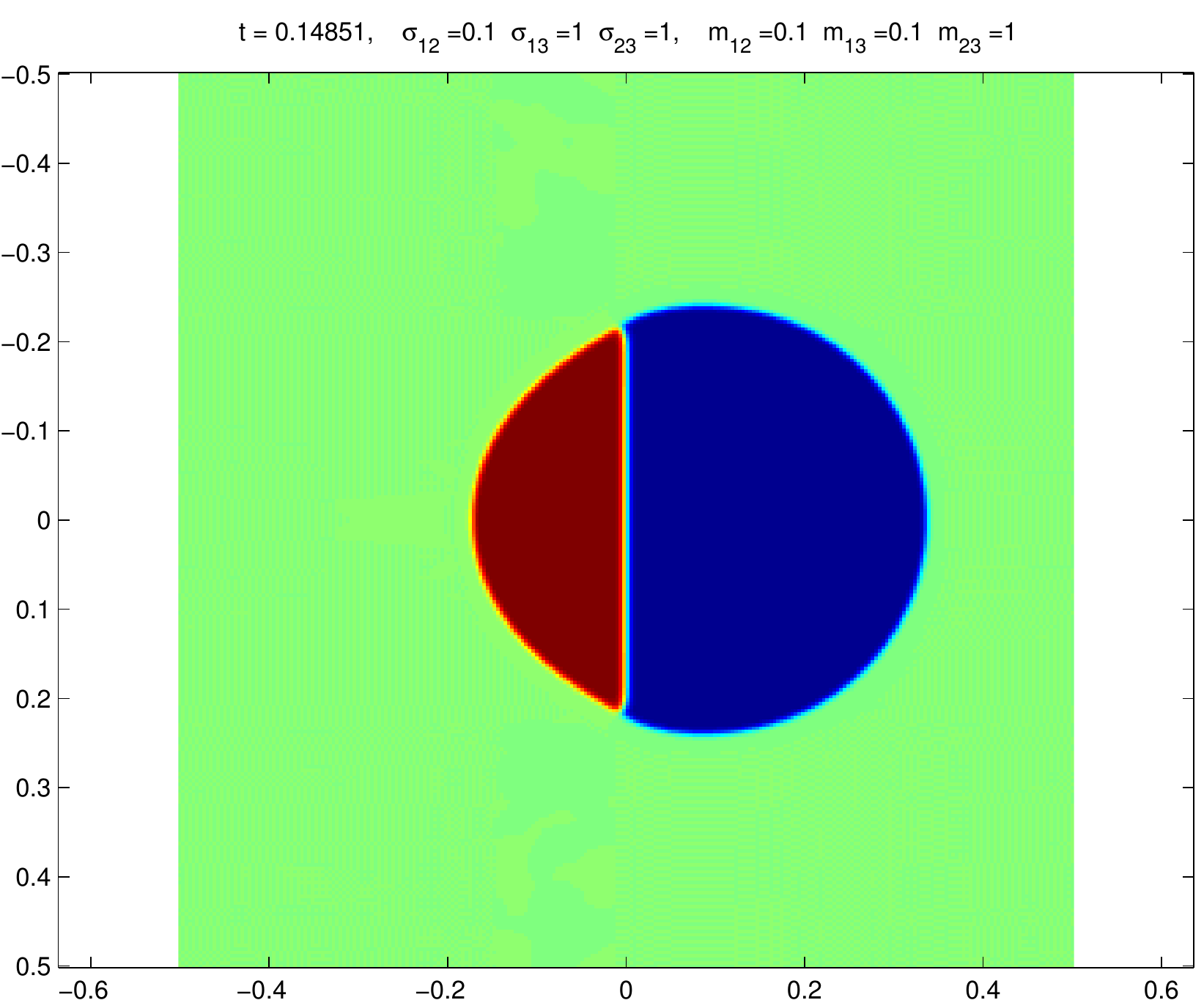} \\
\caption{Examples of multiphase mean curvature flows for different choices of surface tensions and mobility coefficients. Columns from left to right are screenshots of solutions at times $t=0$, $t=0.05$, $t=0.1$ and $t = 0.15$, respectively. From top to bottom, we use respectively $({\boldsymbol \sigma}_1,{\boldsymbol m}_1)$, $({\boldsymbol \sigma}_1,{\boldsymbol m}_2)$,  $({\boldsymbol \sigma}_2,{\boldsymbol m}_1)$, and $({\boldsymbol \sigma}_2,{\boldsymbol m}_2)$, with ${\boldsymbol \sigma}_1 = (1,1,1), {\boldsymbol \sigma}_2 = (0.1,1,1),  {\boldsymbol m}_1 = (1,1,1) \text{ and }  {\boldsymbol m}_2 = (0.1,0.1,1)$.}
\label{Mobility_mean_curvature_flow_2D}
\end{figure}

\paragraph{Experimental consistency}~\\
Figure \eqref{Mobility_mean_curvature_flow_2D} illustrates numerical results obtained with different sets of surface tension coefficients
${\boldsymbol \sigma} = (\sigma_{12},\sigma_{13},\sigma_{23})$ and mobility ${\boldsymbol m} = (m_{12},m_{13},m_{23})$ :
$${\boldsymbol \sigma}_1 = (1,1,1), {\boldsymbol \sigma}_2 = (0.1,1,1),  {\boldsymbol m}_1 = (1,1,1) \text{ and }  {\boldsymbol m}_2 = (0.1,0.1,1).$$
The phases $\Omega_1$, $\Omega_2$ and $\Omega_3$ are represented in blue, red and green colors, respectively. All numerical experiments have been performed with the following numerical parameters: $N = 2^8$, $\epsilon = 1/N$, $\delta_t = 1/N^2$ and $L_1 = L_2 = 1$. In particular, we can notice on this simple example the influence
of  the surface tensions ${\boldsymbol \sigma}$ on the evolution of the triple points (so as to satisfy Herring's condition) and the influence of the mobilities only on the velocity of each interface.
Another important remark is that the width of the diffuse interface depends only on $\epsilon$; it does not depend neither on the surface tensions ${\boldsymbol \sigma}$ nor on mobilities ${\boldsymbol m}$. We believe that this is a major advantage of our approach.

\paragraph{Simulation of wetting phenomena}~\\
Our phase field model can also handle the case of the evolution of a liquid phase on a fixed solid surface by simply imposing a null mobility of the SV and SL interfaces:
$$ m_{LS} = m_{SV} = 0.$$
Two centuries ago, Young \cite{Young_angle_contact} established the optimal shape of a drop in equilibrium on a solid surface.
In particular, Young's law prescribes the contact angle $\theta$ of the liquid on the solid, i.e.
$$ \cos(\theta) = \frac{\sigma_{SV} - \sigma_{LS}}{\sigma_{VL}}$$
which represents the horizontal component of the force balance at the triple point.
The wetting phenomenon was modeled by Cahn~\cite{Cahn_Allen_Cahn_angle} in a phase-field setting.
Cahn proposed to extend the Cahn-Hilliard energy by adding a surface energy term which describes the liquid-solid interaction.
This approach has been used in \cite{Turco2009} for numerical simulations of one droplet, but it cannot be used for angles $\theta \geq \frac{\pi}{2}$.  A different approach~\cite{penalization1} using the smoothed boundary method proposes to compute the Allen-Cahn equation
using generalized Neumann boundary conditions in order to force the correct contact angle condition. Note that an extension of these approaches to many droplets can be found in \cite{multi_droplet}.
More recently, two of the authors of the current paper have proposed  a multiphase field model \cite{Bretin_Masnou_multiphase} which allows freezing the solid phase to approximate droplets' wetting. This approach is equivalent to using null mobilities, {\em i.e.}, $m_{SV} = m_{LS}=0$ and its main advantages are simplicity and accuracy. In particular, it does not impose in any way the contact angle, which is instead implicitly prescribed just by energy minimization.~\\

\begin{figure}[htbp]
\centering
	\includegraphics[width=3.7cm]{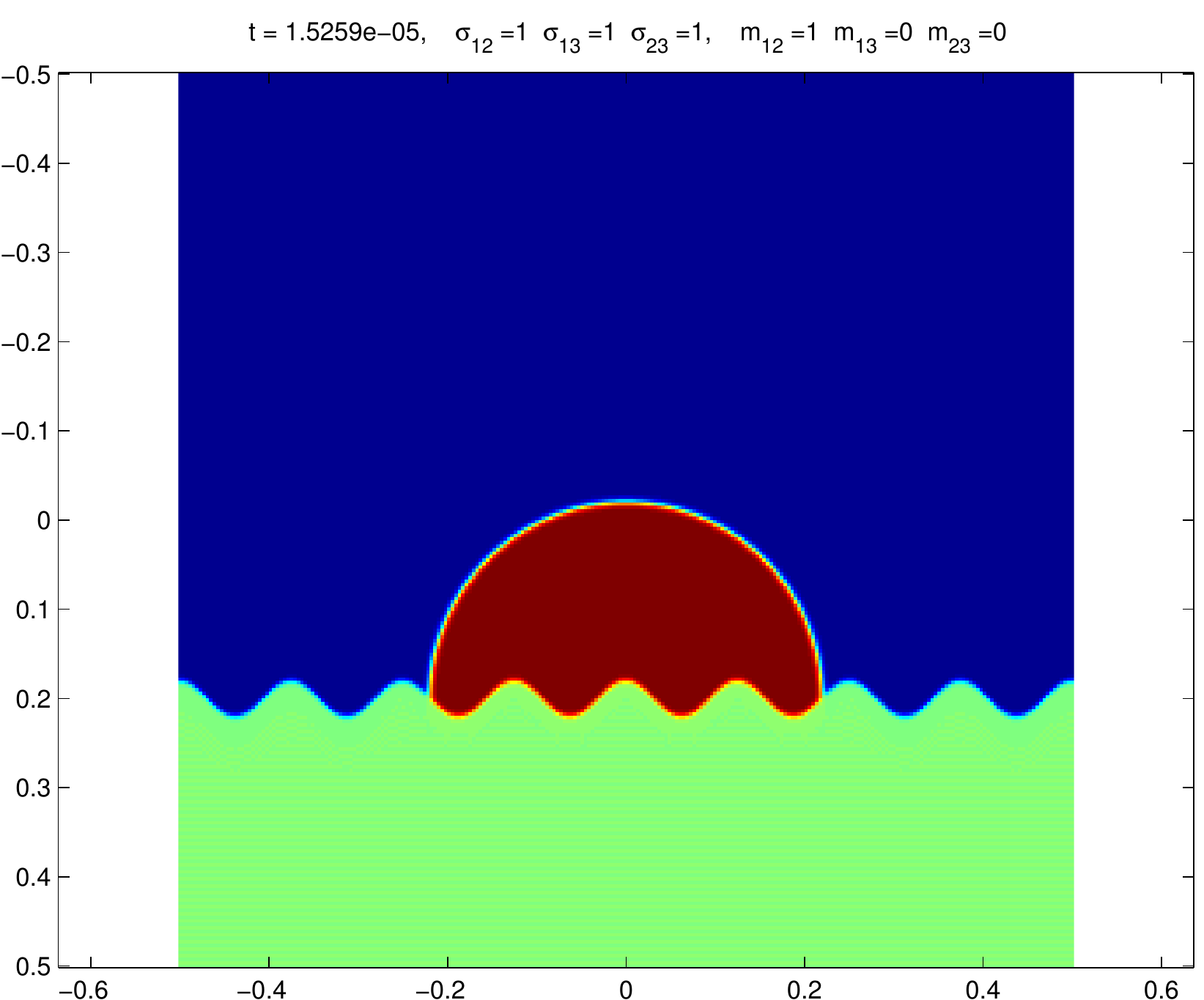}
        \includegraphics[width=3.7cm]{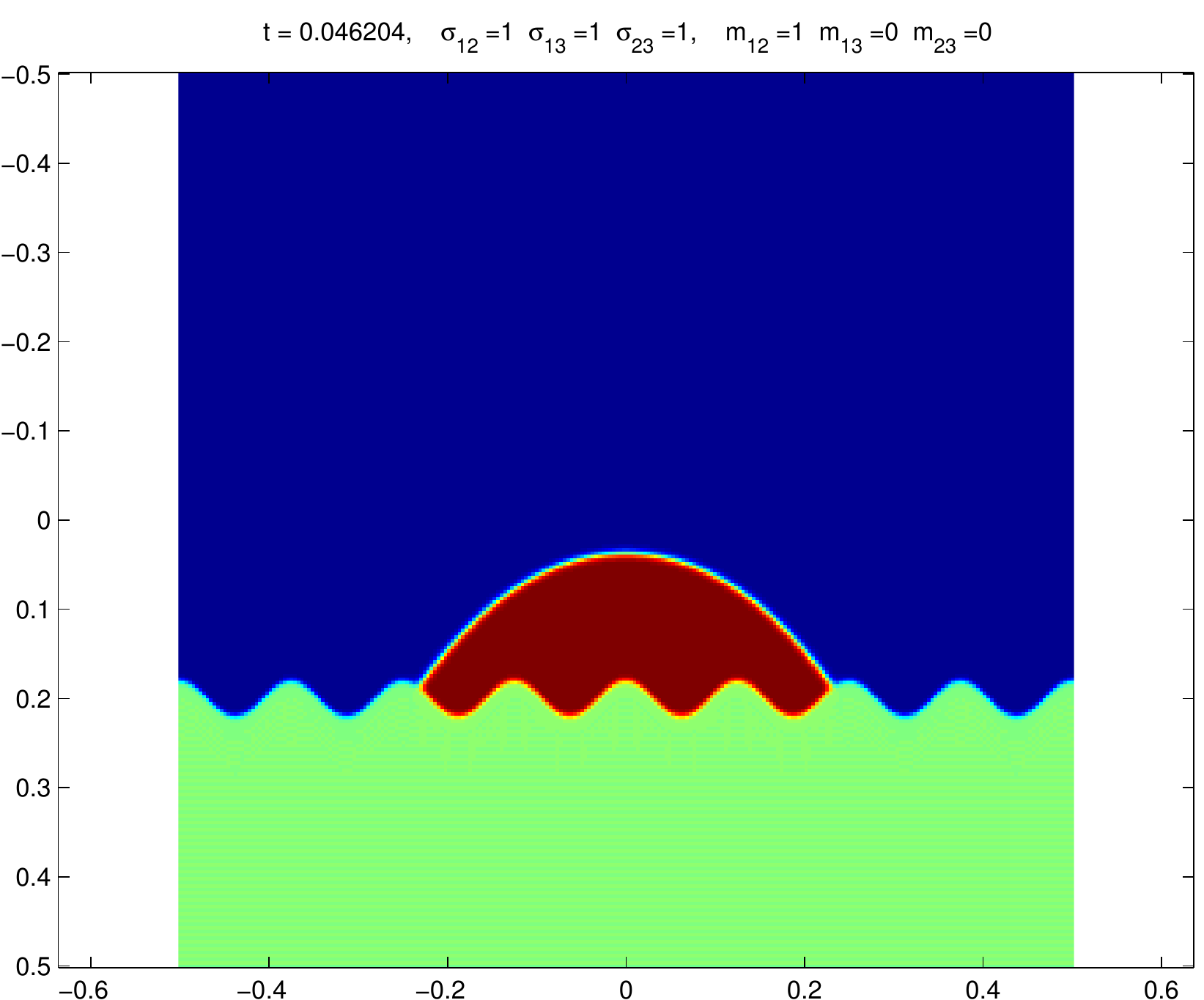}
        \includegraphics[width=3.7cm]{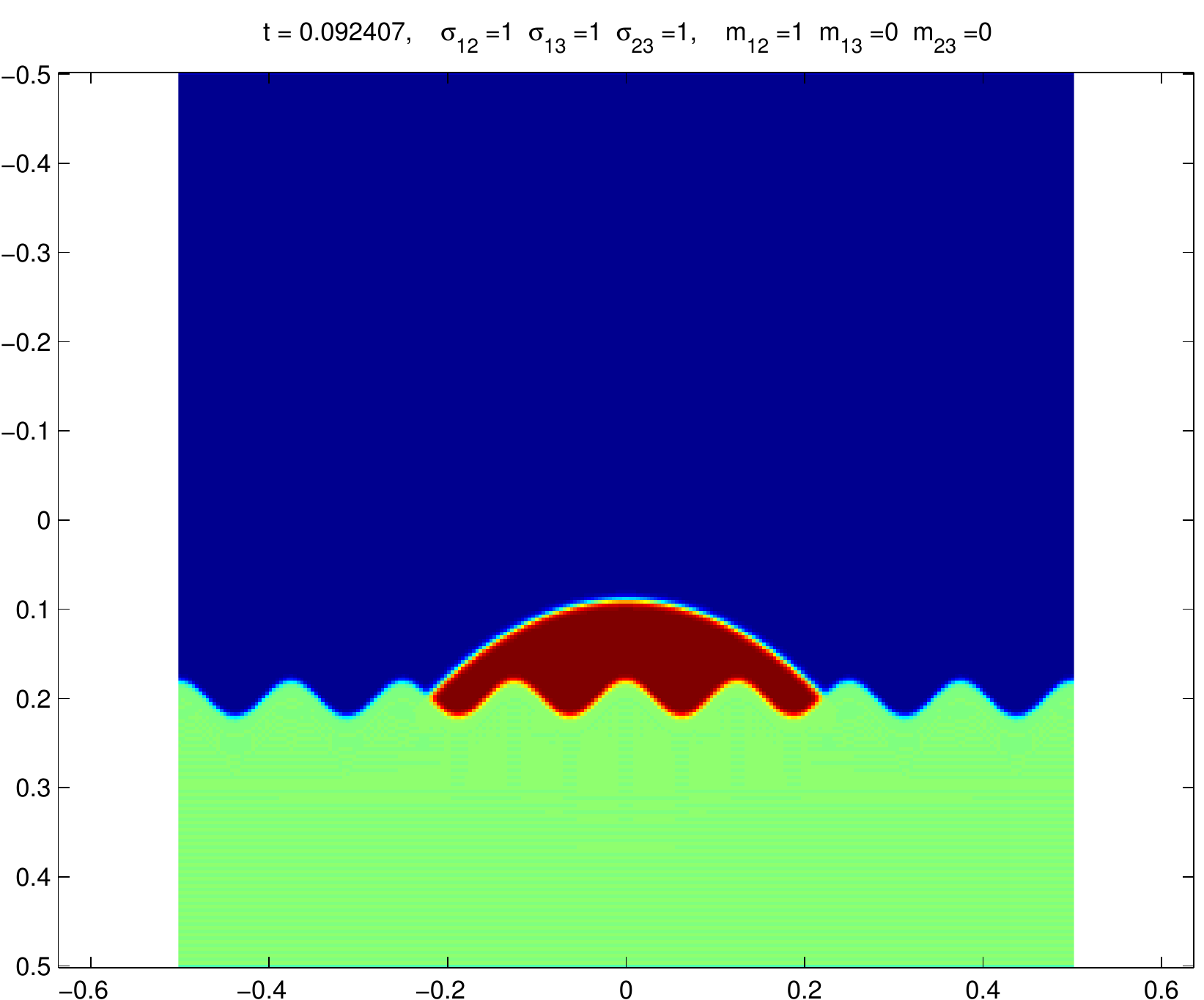}
        \includegraphics[width=3.7cm]{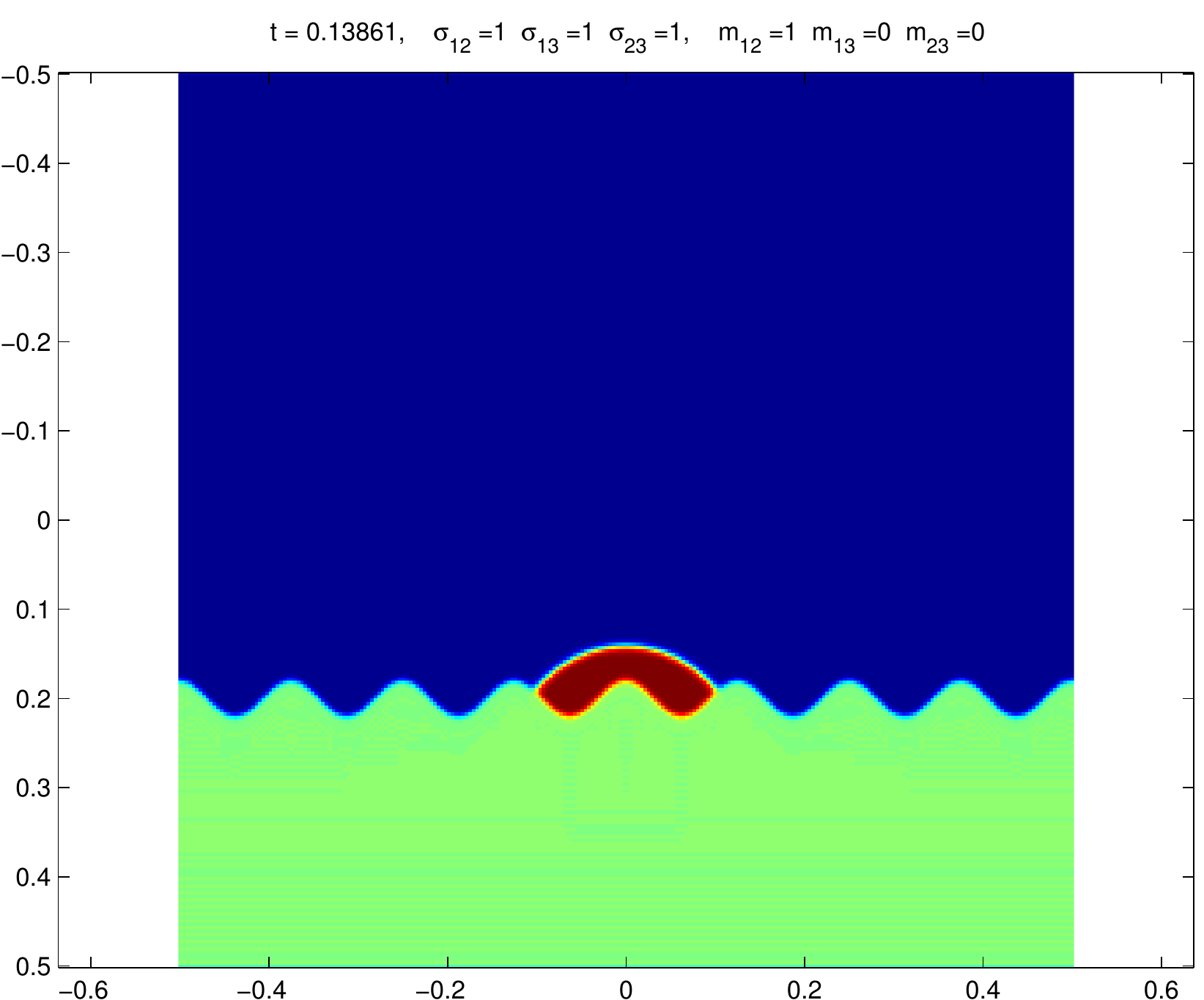} \\
        \includegraphics[width=3.7cm]{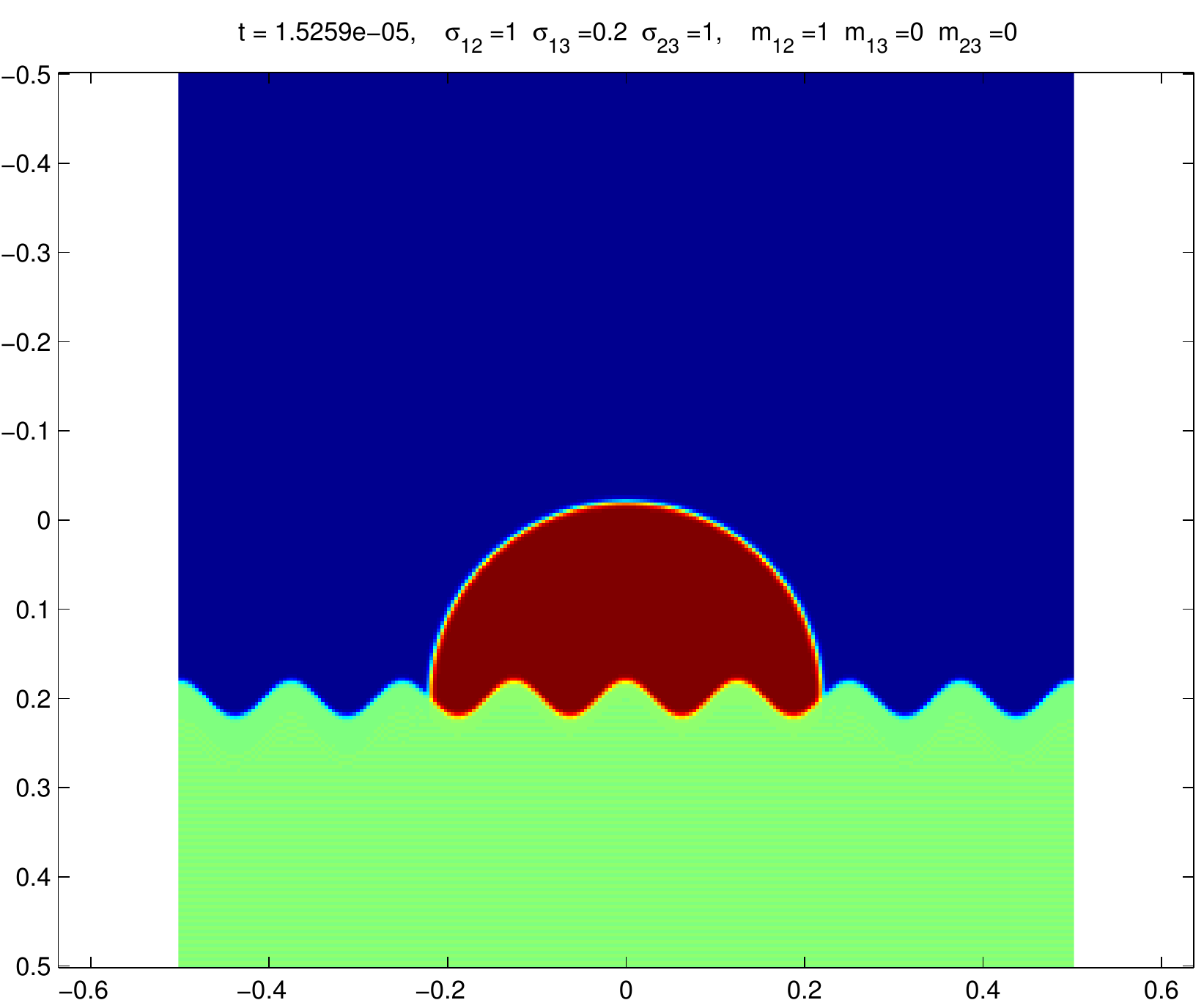}
        \includegraphics[width=3.7cm]{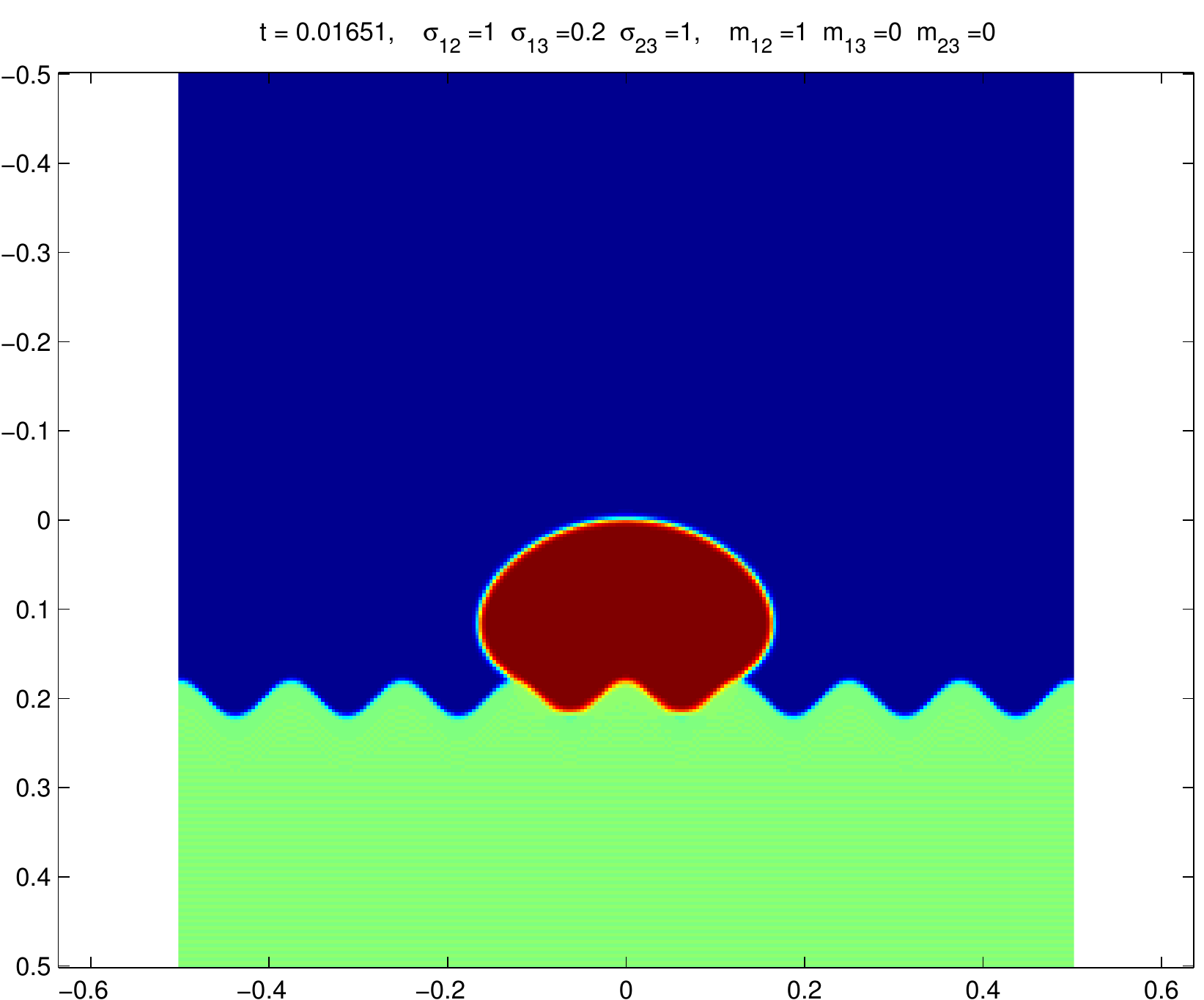}
        \includegraphics[width=3.7cm]{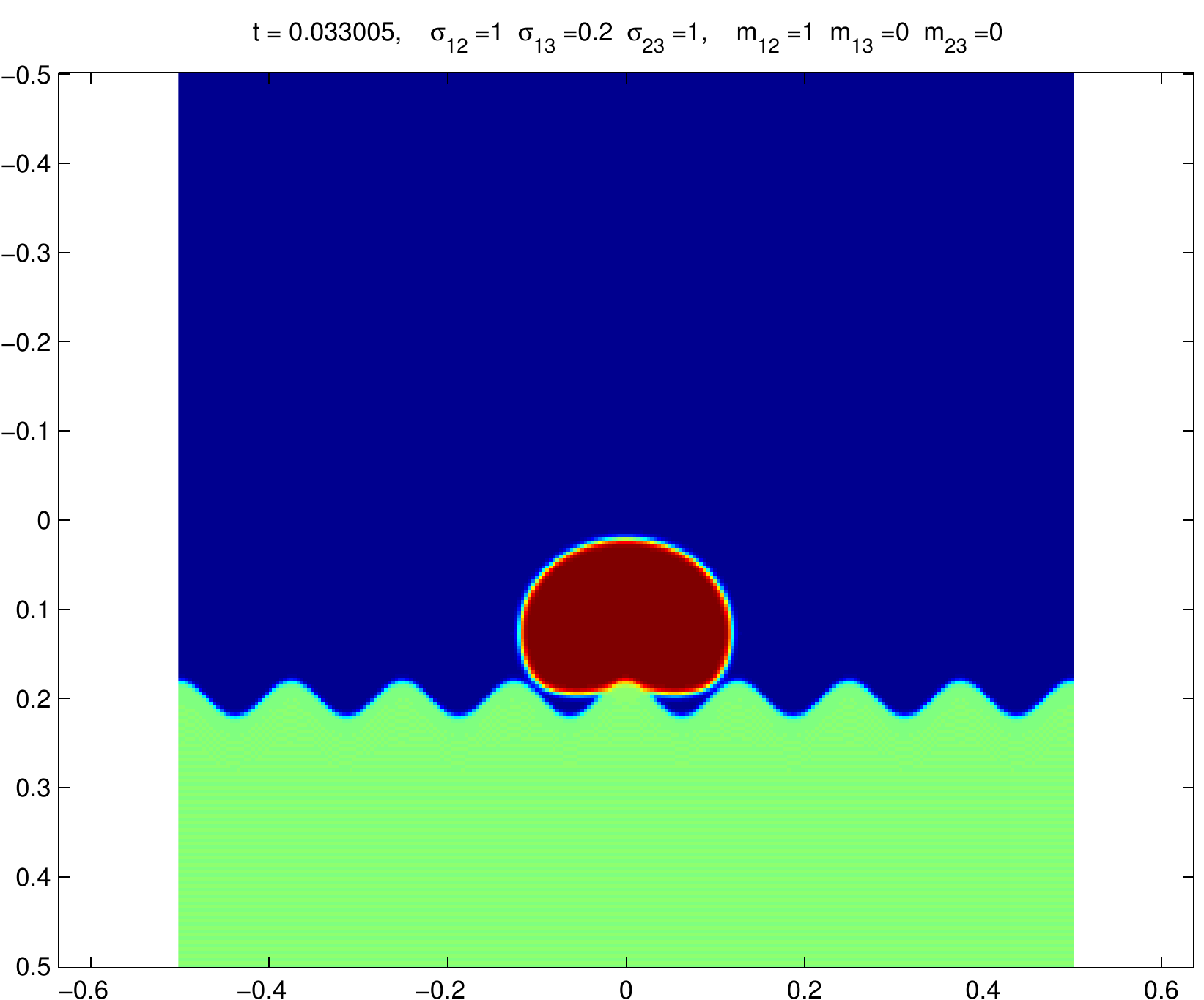}
        \includegraphics[width=3.7cm]{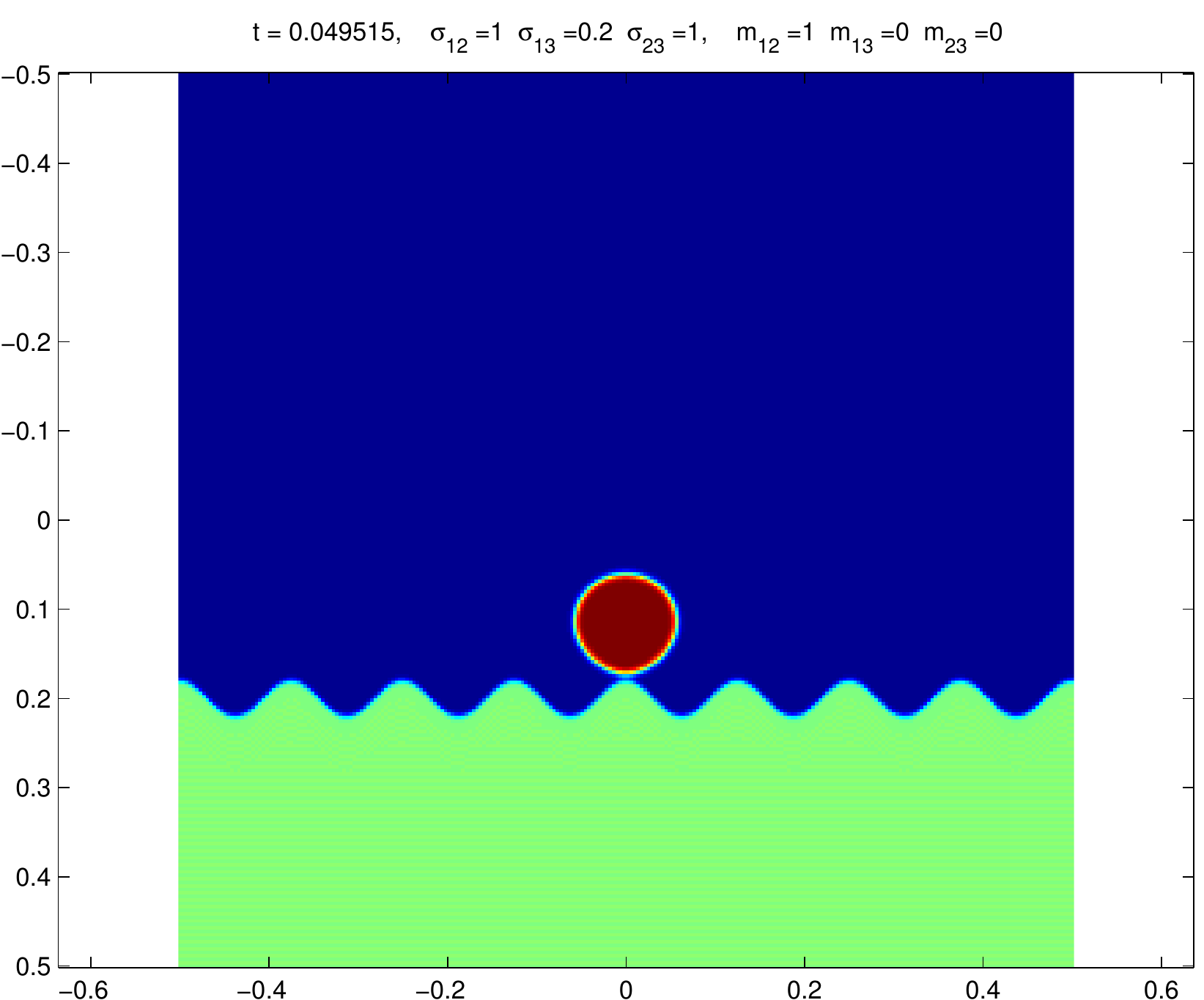} \\
        \includegraphics[width=3.7cm]{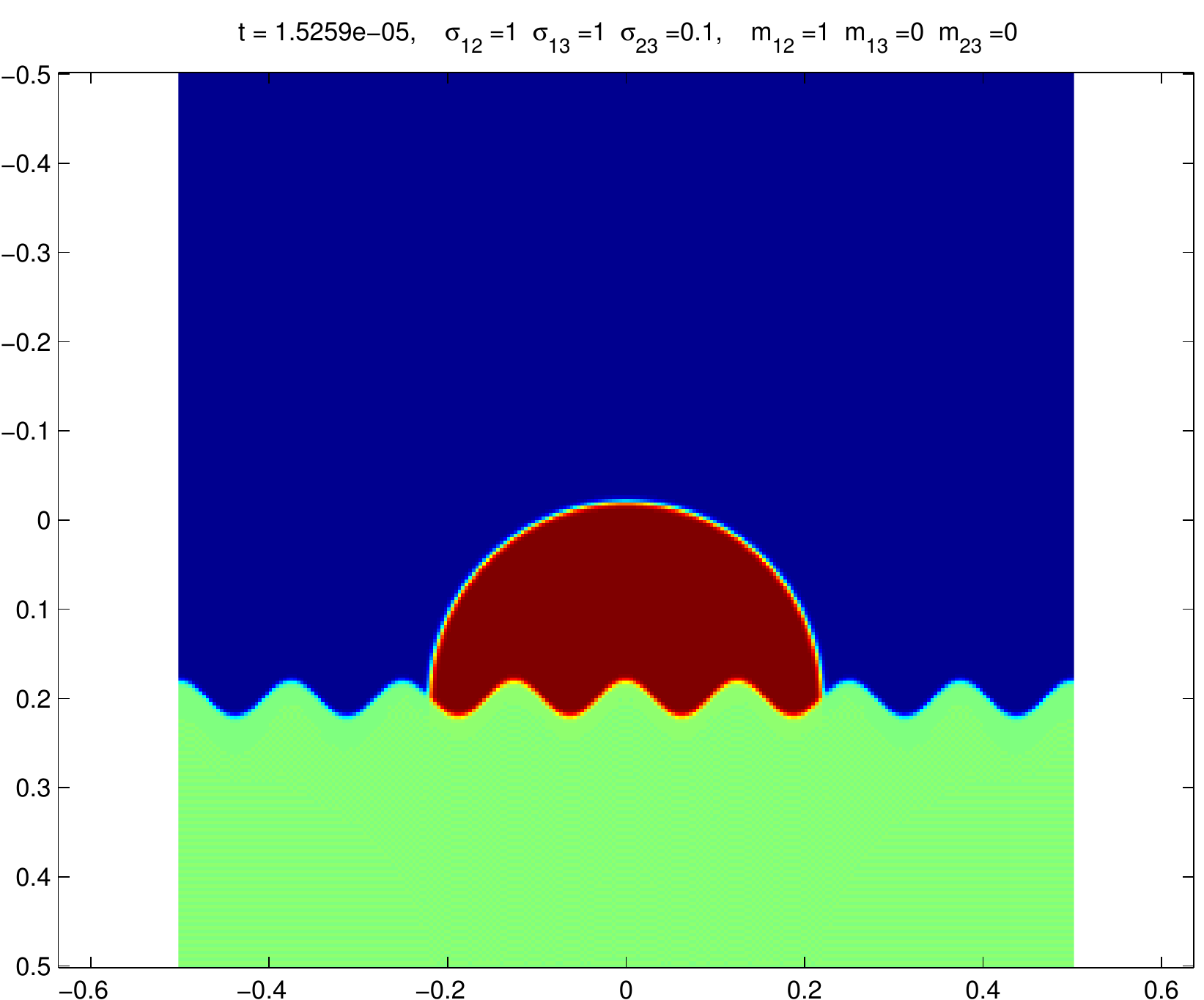}
        \includegraphics[width=3.7cm]{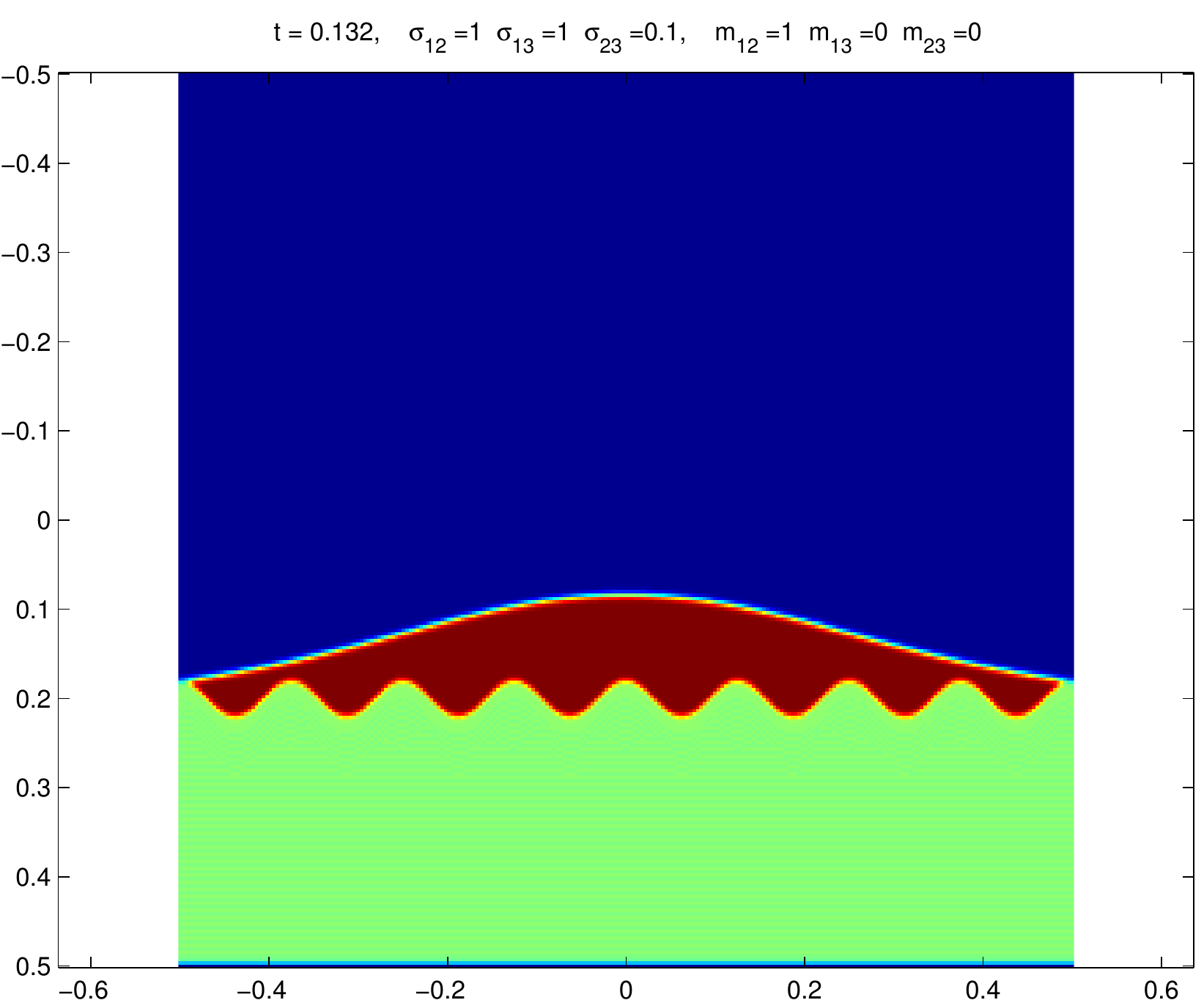}
        \includegraphics[width=3.7cm]{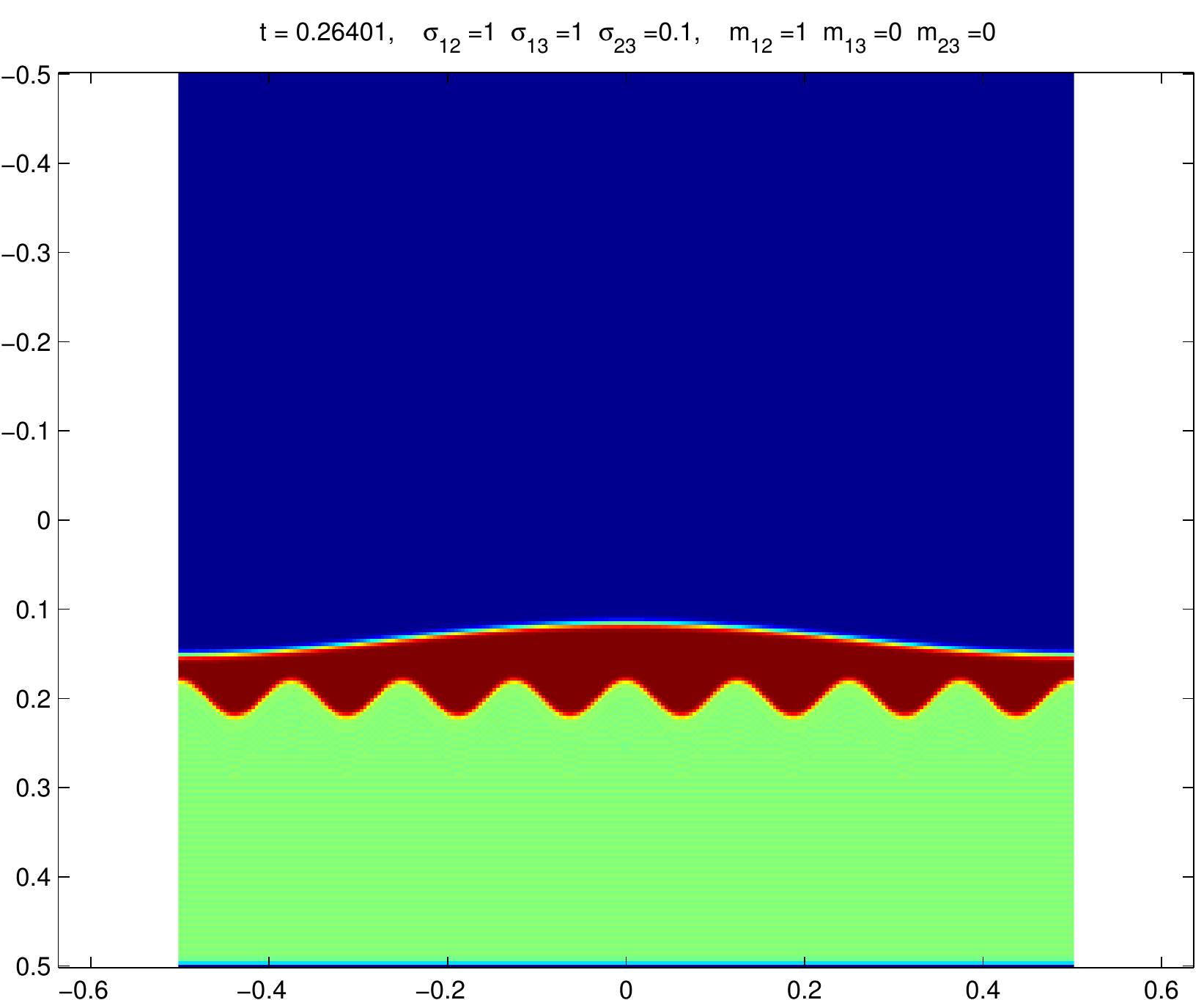}
        \includegraphics[width=3.7cm]{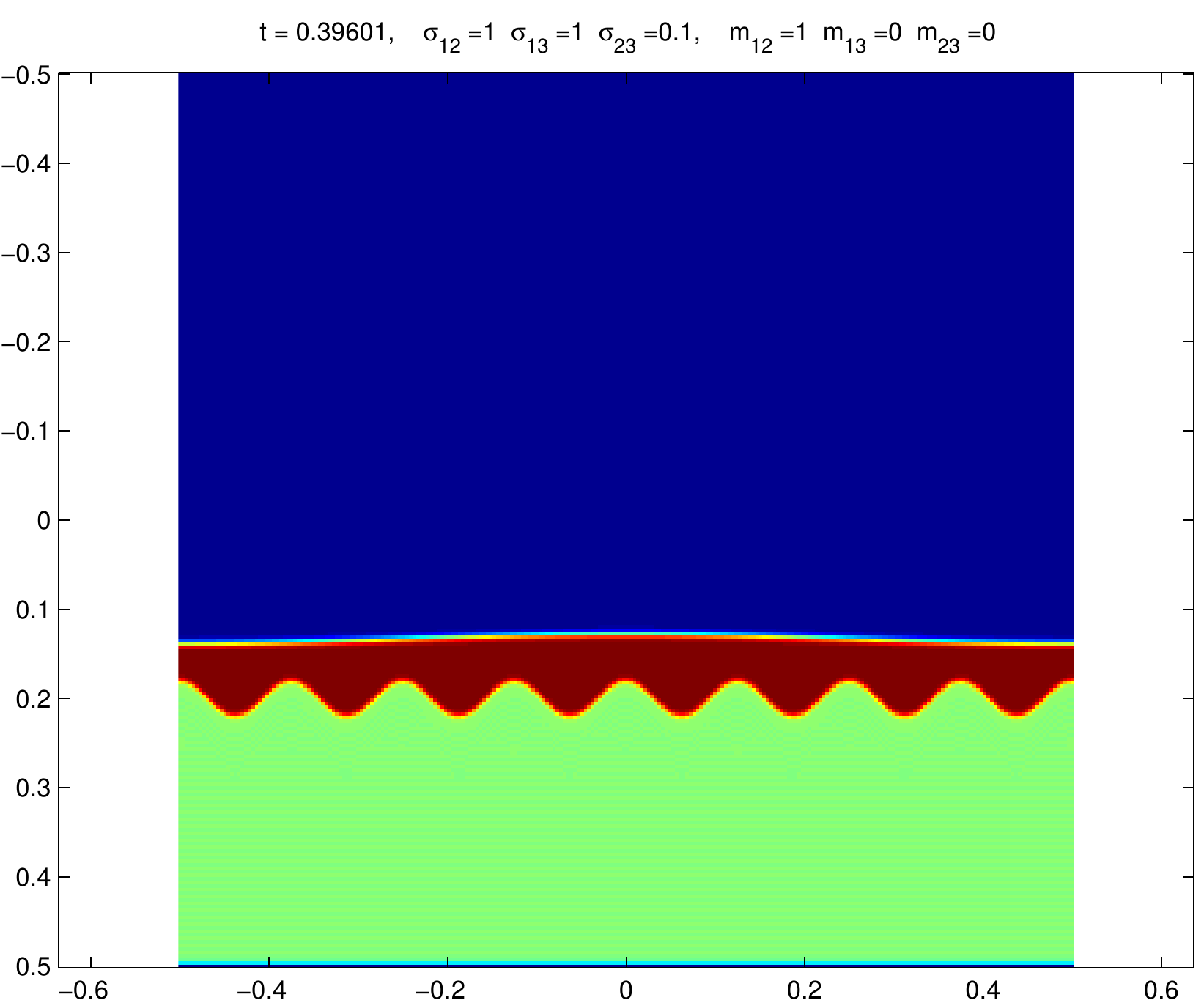} \\
\caption{From left to right, evolution of a liquid phase on a complex solid surface. The mobility coefficients are always defined as ${\boldsymbol m} = (m_{12},m_{13},m_{23}) = (1,0,0)$. The surface tensions are defined from top to bottom as ${\boldsymbol \sigma} = (1,1,1)$, ${\boldsymbol \sigma} = (1,0.2,1)$, and  ${\boldsymbol \sigma} = (1,1,0.2)$, respectively.}
\label{wetting_mean_curvature_flow_2D}
\end{figure}

Figure \eqref{wetting_mean_curvature_flow_2D} illustrates numerical results obtained with mobilities equal to
$${\boldsymbol m} = (m_{12},m_{13},m_{23}) = (1,0,0)$$
and using respectively ${\boldsymbol \sigma} = (1,1,1)$, ${\boldsymbol \sigma} = (1,0.2,1)$ and  ${\boldsymbol \sigma} = (1,1,0.2)$. The liquid, vapor, and solid phases are represented in
red, blue and green colors, respectively.  Numerical computations were performed using $N = 2^8$, $\epsilon = 1/N$, $\delta_t = 1/N^2$, and $L_1, = L_2 = 1$ for each experiment.
In particular, we notice the ability of our model to treat the case of null mobilities and these experiments show the high influence of the contact angle
on the evolution of the liquid phase. We emphasize, again, that our model does not prescribe the contact angle. Its value is rather a straightforward consequence of the multiphase interface energy considered
in each simulation.

 \subsection{Application to nanowire growth}
\paragraph{The explicit shape for nanowire in quasi-static growth} ~\\
The first stage of nanowire growth is much less documented than the {\em stationary} growth regime, when the nanowire length and thickness evolve at relative constant rates. This is because, in the first stages of nanowire growth, both the surface tension at the SV interface and the liquid droplet geometry evolve from the horizontal direction (when the liquid droplet lays between the S and V phases) to the vertical one (when the droplet is pinned at the top of the nanowire). The following derivation (i) assumes the isotropy of all interfaces, and thus represents only a rough approximation of the real physical situation, and (ii) generalizes an early attempt in \cite{gosele} by including the vertical component of the force balance at the triple point.

\vspace*{1cm}
\begin{figure}[h!]
\centering
\begin{tikzpicture}[scale=0.25]
\draw[ultra thin,lightgray](-8,-14) grid (10,9);

\draw[thick,red](4,0) arc(56.3: 123.7: 7.21);
\draw[thin,dashed,red](-4,0) arc(123.7: 416.3: 7.21);
\draw[black,line width = 0.5pt](-8,0)--(-4,0);
\draw[black,line width = 0.5pt](4,0)--(10,0);
\draw[thick,color=red,->,>=latex](4,0)--(1,2);
\draw[thick,color=blue,->,>=latex](4,0)--(2.5,-2);
\draw[thick,blue](-4,0) arc(216.8: 323.2: 5);
\draw[thin,dashed,blue](0,3) circle(5);
\draw[thick,color=black,->,>=latex](4,0)--(8.5,0);
\draw[thick,color=black,<->,>=latex](-4,0)--(0,0);
\draw(-3.5,-0.7) node[right]  {$R_0$};
\draw(-1,-10) node[left] {${\bf S}$};
\draw(2.5, -0.2)  node[left] {${\bf L}$};
\draw(-1, 4)  node[left] {${\bf V}$};
\draw(8.5, -1.8)  node[above] {${{\sigma}_{SV}T_{SV}}$};
\draw(2.3, 1.8)  node[above] {${{\sigma}_{LV}T_{LV}}$};
\draw(2, -4)  node[above] {${{\sigma}_{LS}T_{LS}}$};
\draw[thick,red,->](5,0) arc(0: 146.3: 1);
\draw[thick,blue,<-](5,0) arc(0: -126.9: 1);
\draw(4.5,1.2) node[right]  {$\theta_V$};
\draw(4,-1.7) node[right] {$\theta_S$};

\draw[ultra thin,lightgray](22,-14) grid (42,9);
\draw[thick,red](36,-10) arc(56.3: 123.7: 7.21);
\draw[black,line width = 0.5pt](24,-10)--(28,-10);
\draw[black,line width = 0.5pt](36,-10)--(42,-10);
\draw[thick,color=red,->,>=latex](36,-10)--(33,-8);
\draw[thick,color=blue,->,>=latex](36,-10)--(34.5,-12);
\draw[thick,blue](28,-10) arc(216.8: 323.2: 5);
\draw[thick,color=black,->,>=latex](36,-10)--(40.5,-10);
\draw[thick,blue](34,0) arc(53: 127: 3.33);
\draw[thick,red](34,0) arc(-33.7: 213.7: 2.4);
\draw[thick,green] (36,-10) ..controls +(-2,0) and +(0,-1).. (34,0);
\draw[thick,green] (28,-10) ..controls +(2,0) and +(0,-1).. (30,0);
\draw[thick,dotted,black,->,>=latex](34,0)--(34,-4.5);
\draw[thick,dotted,red,->,>=latex](34,0)--(36,3);
\draw[thick,dotted,blue,->,>=latex](34,0)--(32,1.5);
\draw[thick,black,<->](32,5)--(34,5);
\draw[thick,black,<->](32,-7)--(34.5,-7);
\draw(33,4.9) node[above] {$R_{\frac{\pi}{2}}$};
\draw(33,-7) node[above] {$R_\alpha$};
\end{tikzpicture}
\caption{Left: a typical VLS configuration when the liquid droplet reacts with the substrate (as is the case for Si-Au). In the equilibrium configuration where the balance of interfacial forces holds, the interfaces SL and LV have constant curvatures. The dotted circles illustrate that interfaces LS and LV are contained in curves of constant curvature. Right: During the nanowire growth the three surface tension vectors turn from the initial position where the superficial tension vector $\sigma_{SV}T_{SV}$ is horizontal toward the stationary regime where it is vertical. Relative positions of surface tensions (at equilibirum) remain unchanged.}
\label{fig:VLS_initial_volume}
\end{figure}
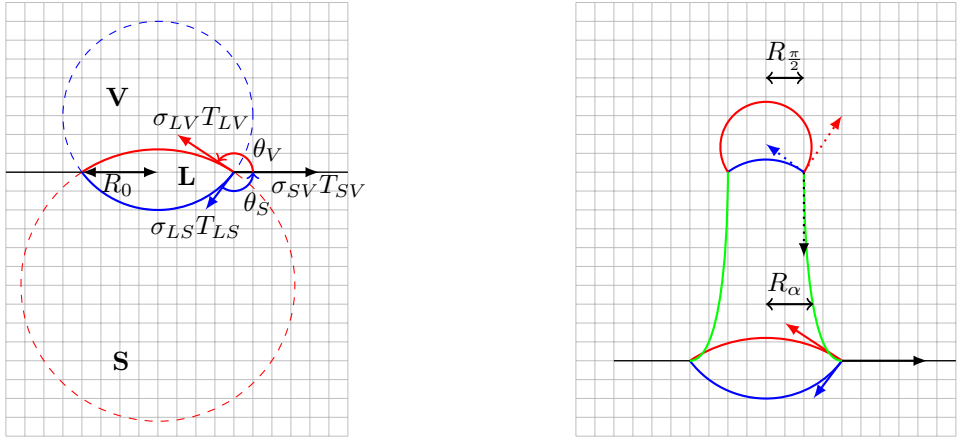

We adopt the quasi-static point of view of sharp interface theory, i.e., we assume that the time increase rate of the nanowire volume is not comparable with the action of geometric energies toward spatial equilibrium. In other words, all interfaces are at space equilibrium at every time. We start from an initial configuration where the triple point is at the equilibrium illustrated in Fig \ref{fig:VLS_initial_volume}, left. In this configuration, both SV and VL interfaces have constant curvature. Let $\theta_V$ denote the angle between the surface tension vectors $\sigma_{SV}T_{SV}$ and $\sigma_{VL}T_{VL}$ (see Fig. \ref{fig:VLS_initial_volume}, left). From the balance of surface tensions at the triple point we obtain
$$
\cos\theta_V = \frac{b^2-a^2-1}{2a},\qquad \cos\theta_S = \frac{a^2-b^2-1}{2b},
$$
where we used $a=\frac{\sigma_{VL}}{\sigma_{SV}}$ and $b=\frac{\sigma_{LS}}{\sigma_{SV}}.$
The volume in the liquid phase can be expressed as the sum of the volumes of two spherical caps and is given by
\begin{equation}
V_L(0) = \frac{\pi R_0^3}{3}\left[f(\theta_V) + f(\theta_S)\right]
\label{eqn:volume_0}
\end{equation}
where $f(\theta)= \frac{(1+\cos\theta)(2-\cos\theta)}{\sin^3\theta}$ and $R_0$ denotes in Fig. \ref{fig:VLS_initial_volume}, left, the half-distance between the triple points in the initial configuration.

In the first stages of nanowire growth, the frame consisting of the three vectors $\sigma_{SV}T_{SV},$ $\sigma_{VL}T_{VL}$, and $\sigma_{LS}T_{LS}$ turns from a position where the surface tension $\sigma_{SV}T_{SV}$ is horizontal (left picture, or equivalently right picture in Fig. \ref{fig:VLS_initial_volume}, lower part) to a position where it is vertical (right picture in Fig. \ref{fig:VLS_initial_volume}, upper part).
If $\alpha$ denotes the angle between the surface tension vector $\sigma_{SV}T_{SV}$ and the horizontal axis, then, during the first stages of nanowire growth, the angles $\theta_S$ and $\theta_V$ become
$$
\theta_S \rightarrow \theta_S + \alpha,\qquad \theta_V \rightarrow \theta_V - \alpha.
$$
Assuming the liquid volume is constant, and using relation (\ref{eqn:volume_0}), one obtains for the radius of the liquid droplet
\begin{equation}
R_\alpha(\alpha) = R_0 \left[\frac{f(\theta_V-\alpha)+f(\theta_S+\alpha)}{f(\theta_V)+f(\theta_S)}\right]^{1/3},
\label{eqn:radius}
\end{equation}
so that the radius $r$ of the droplet in stationary regime is determined completely by the size of the initial section containing the triple points (which can be measured experimentally) and by relative ratios $a$ and $b.$
It follows that the radius of the droplet is in $[R_{\alpha}(\pi/2),R_0]$ and, in the stationary growth regime, it is given by (\ref{eqn:radius}) for $\alpha=\pi/2$.
The inverse function of $R_{\alpha}$ can be also  defined for all $r \in[R_\alpha(\pi/2),R_\alpha(0)]$ by solving the equation
$$ R_\alpha(R^{-1}_{\alpha}(r)) = r,$$
which can be done in practice with a Newton-type algorithm.
Moreover, if we introduce a vertical height $h$ measured from the initial triple point plane, the boundary of the solid phase (in green in Fig. \ref{fig:VLS_initial_volume})
is described by $(r,h(r)).$ Taking the derivative with respect to $r$, we obtain for the tangent unit vector to the SV interface
$$
T_{SV}= (-\cos(\alpha(r)),\sin(\alpha(r))) = \frac{1}{\sqrt{1+(h'(r))^2}}(1,h'(r)),
$$
where $\alpha(r) = R^{-1}_\alpha(r)$. Finally,  this shows that  $h$ satisfies
\begin{equation}
h'(r)  =  - \tan \left(R^{-1}_\alpha(r)\right),  \text{ with }   h(R_0) = 0,
\label{function-h}
\end{equation}
and numerical integration of this differential equation provides an approximation of the nanowire shape $r\mapsto h(r)$ in the first stages of growth.

For both approximations, theoretical and numerical, surface tensions play a key role. In particular, they determine fully the evolution of both the droplet radius and position during the first stages of the nanowire growth. Table~\ref{ref:table} below summarizes some of these values for various catalysts, adatoms, and crystalline orientation~\cite{gosele}. We will use these values in our simulations of nanowire growth to compare the theoretical profile $r\mapsto h(r)$ with the shape provided numerically by our phase field model.
\begin{table}[!h]
\begin{center}
\begin{tabular}{|c|c|c|c|}
\hline
 & $\sigma_{LS} (J/m^2)$  & $\sigma_{LV} (J/m^2)$ & $\sigma_{SV} (J/m^2)$\cr
\hline\hline
Au-Si(111) & 0.62 & 0.85 & 1.24 \cr\hline
Au-Si(100) & 0.62 & 0.85 & 1.36 \cr\hline
Au-Si(311) & 0.62 & 0.85 & 1.38 \cr\hline
Au-Si(110) & 0.62 & 0.85 & 1.43 \cr\hline
Au-Ge(111) & 0.55 & 0.94 & 1.06 \cr\hline
\end{tabular}
\end{center}
\caption{Known surface tensions for various catalysts, adatoms, and crystalline orientations.}\label{ref:table}
\end{table}

\begin{figure}[hbp]
\centering
\includegraphics[scale=0.3]{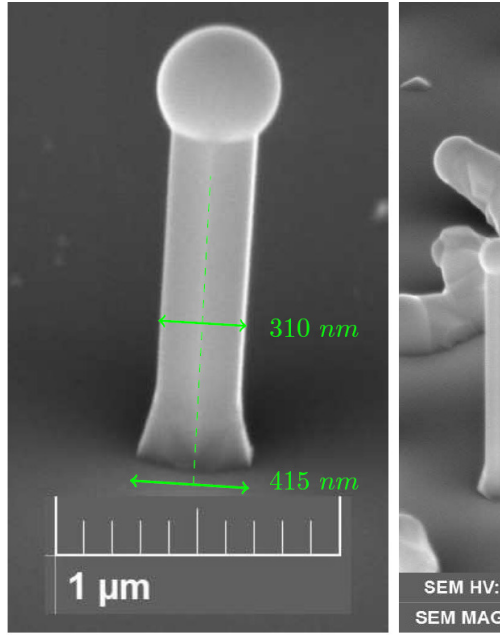}
\caption{Left: estimates of the droplet's and nanowire's diameters in the stationary growth (fixed diameter) regime for the experiment introduced in Figure~\ref{fig:fil_isole}. Right: at a larger scale, the reduction in the nanowire diameter is typical for the first stages of the nanowires growth by VLS~\cite{Benali2017}.}
\label{fig:diameter_evolution}
\end{figure}

\begin{figure}[htbp]
\centering
	 \includegraphics[width=7cm]{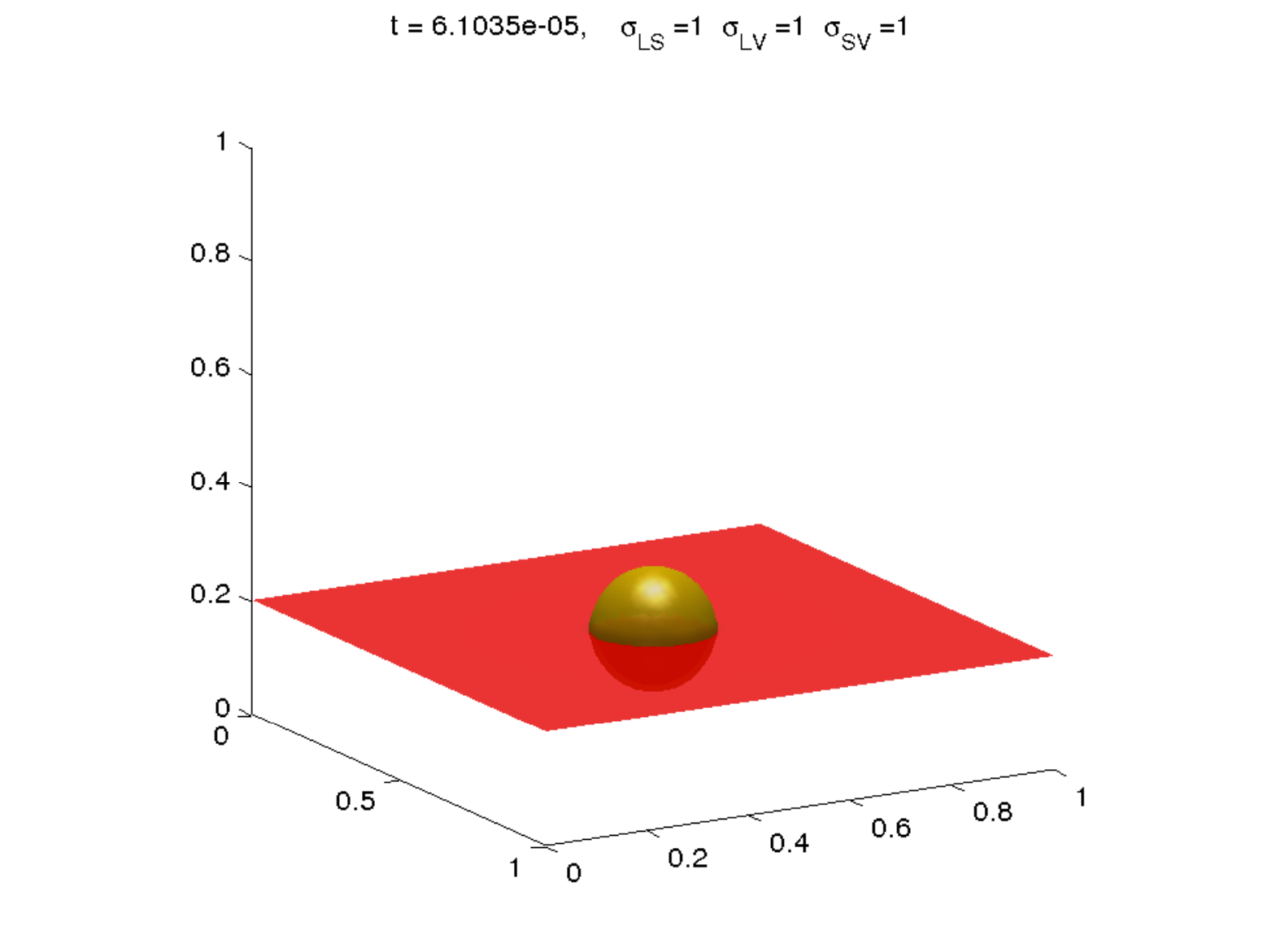}
        \includegraphics[width=7cm]{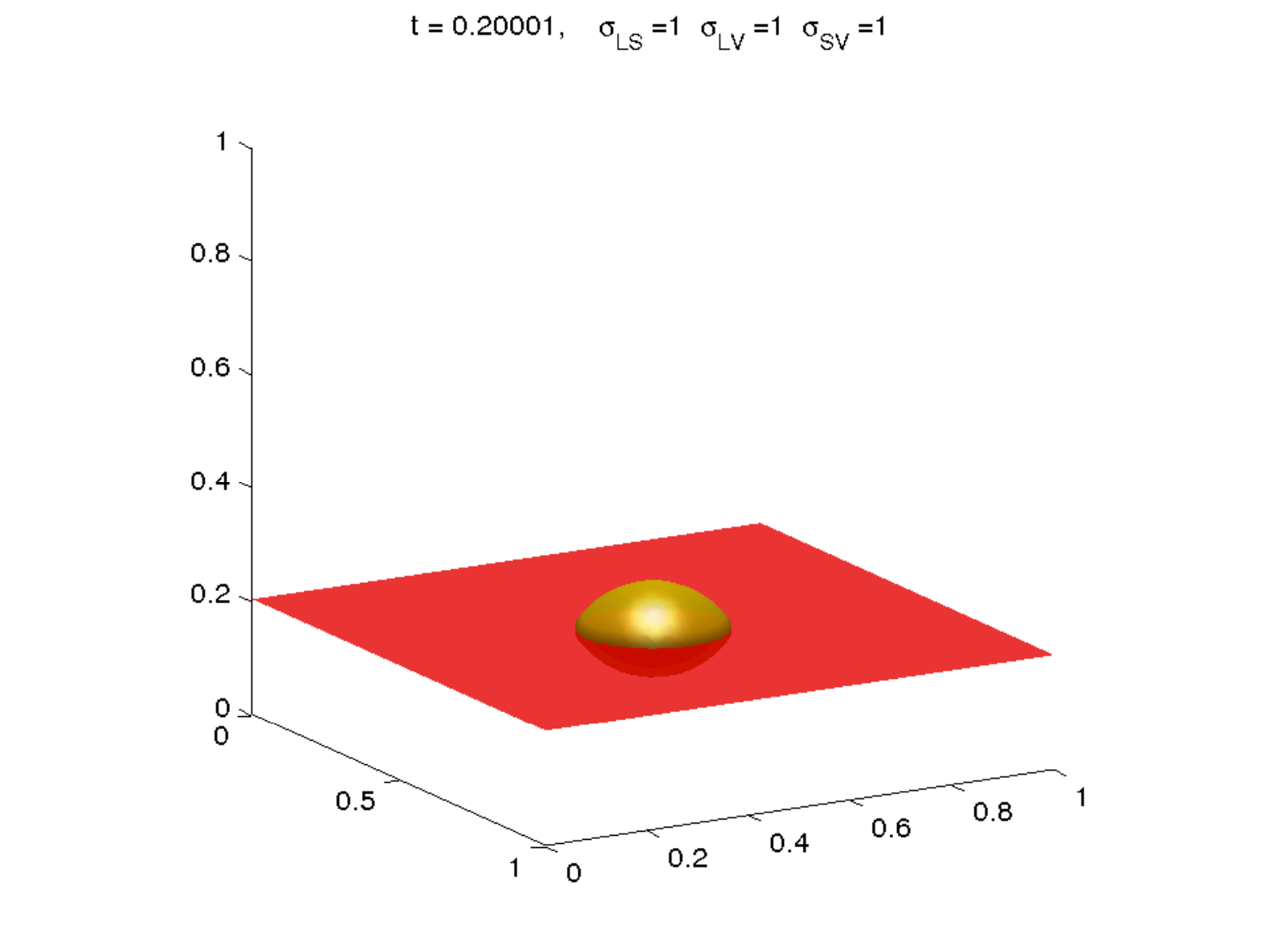}
        \includegraphics[width=7cm]{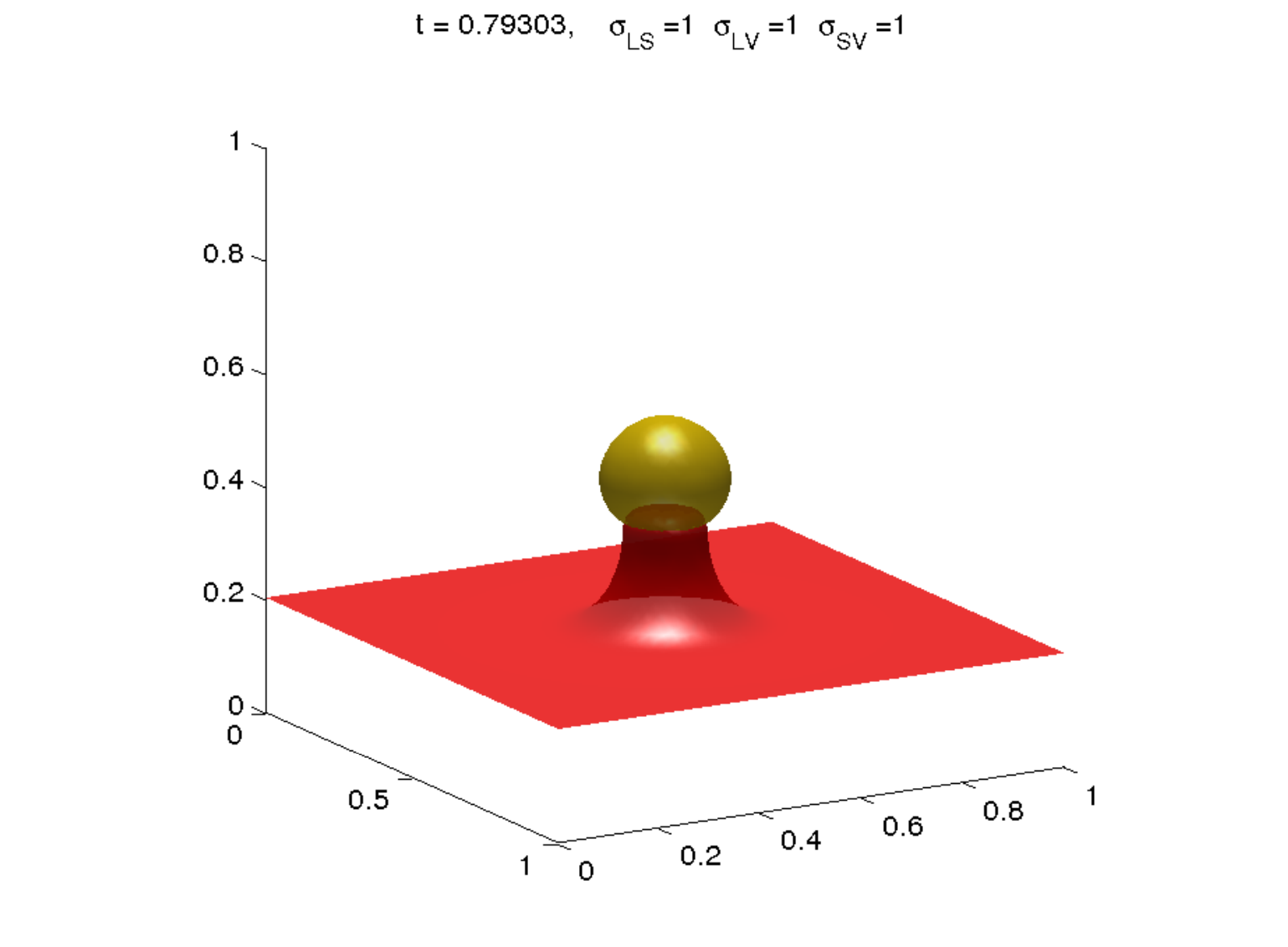}
        \includegraphics[width=7cm]{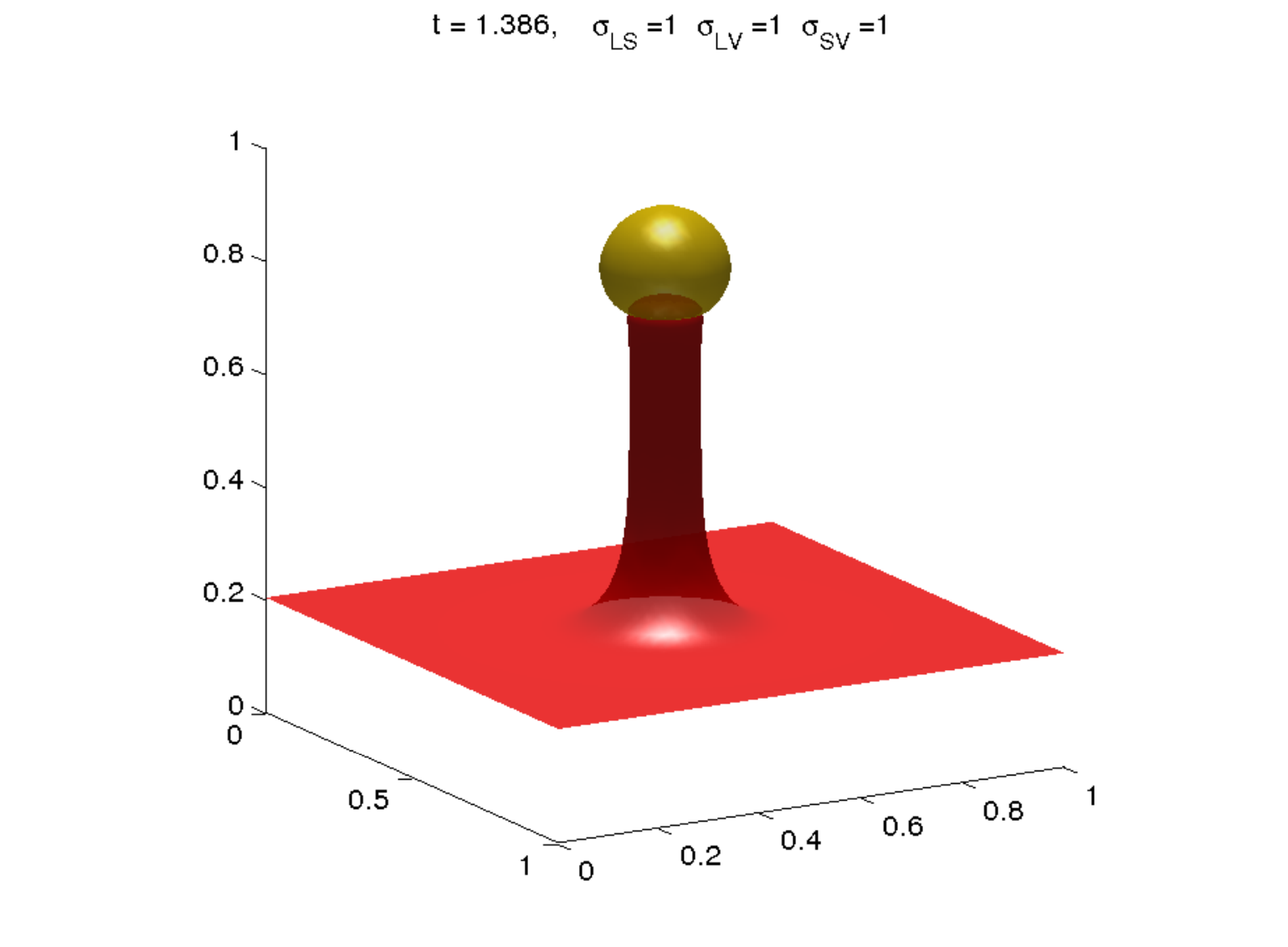}
\caption{Simulation of a $3D$ nanowire evolution using $\sigma_{iso} = (1,1,1)$, captured at different times (from left to right and top to bottom). Starting from the initial configuration, the second image shows the result of the evolution without increase of the solid phase volume and with equal mobilities. Then, the quasi-static model with volume increase and inhomogenous mobilities is used for simulating the two last configurations.}
\label{Nanowire1_isotropic_3D}
\end{figure}

As a quantative illustration of the evolution of nanowire's diameter in the first stages of growth, we give in the left picture of Figure~\ref{fig:diameter_evolution} numerical values obtained in an experimental situation already presented in Figure~\ref{fig:fil_isole}, i.e. using GaAs catalyzed by Ga droplets on a Si(111) substrate covered by its native oxide. Although the precise situation described in Figure \ref{fig:VLS_initial_volume} is actually typical for Au-Si droplets on Si substrates, the discussion above remains valid for predicting the decrease of the nanowire diameter in the first stages of the growth in the experiment shown in Figure~\ref{fig:diameter_evolution}. Using the length scale of the scanning electron microscopy image and neglecting anisotropy, we can estimate in Figure \ref{fig:diameter_evolution} the initial diameter at 415 nm, while the nanowire diameter in the stationnary growth regime is 310 nm. This represents a diameter reduction of $25\%$ of the initial diameter, which has the good order of magnitude when compared to the numerical values predicted by using values in Table \ref{ref:table}. The right picture in Figure \ref{fig:diameter_evolution} shows the same phenomenon on several nanowires at a different picture scale, enforcing the generic character of the diameter reduction in the first stages of nanowire growth~\cite{Boudaa2015, Mavel2017}. It is interesting to compare these experimental results with the numerical simulations of the next section (see, in particular, Figure~\ref{Nanowire1_isotropic_3D}): obviously, the method we propose is able to reproduce with accuracy the physical reality.

\begin{figure}[htbp]
\centering
	 \includegraphics[width=5cm]{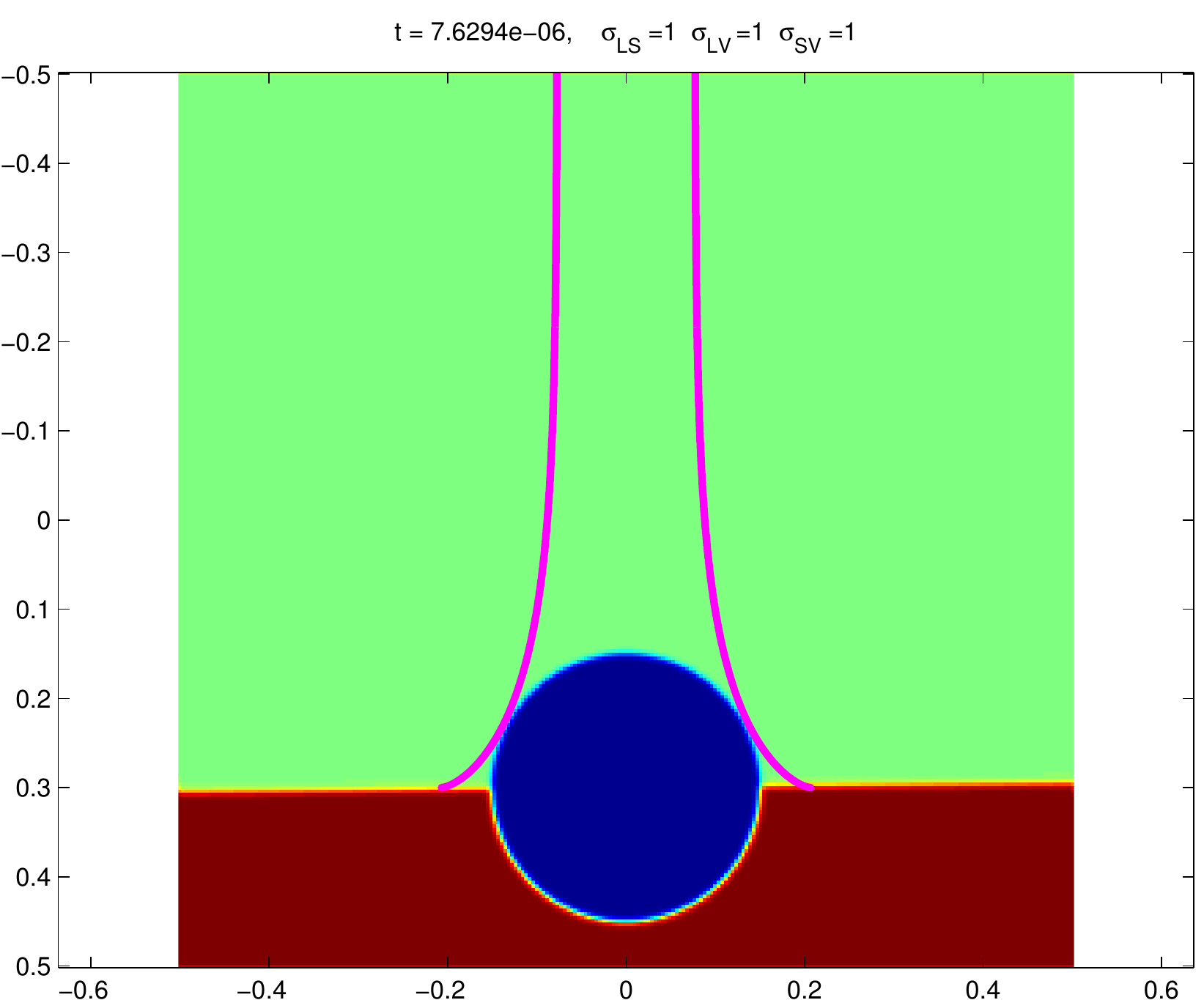}
	 \includegraphics[width=5cm]{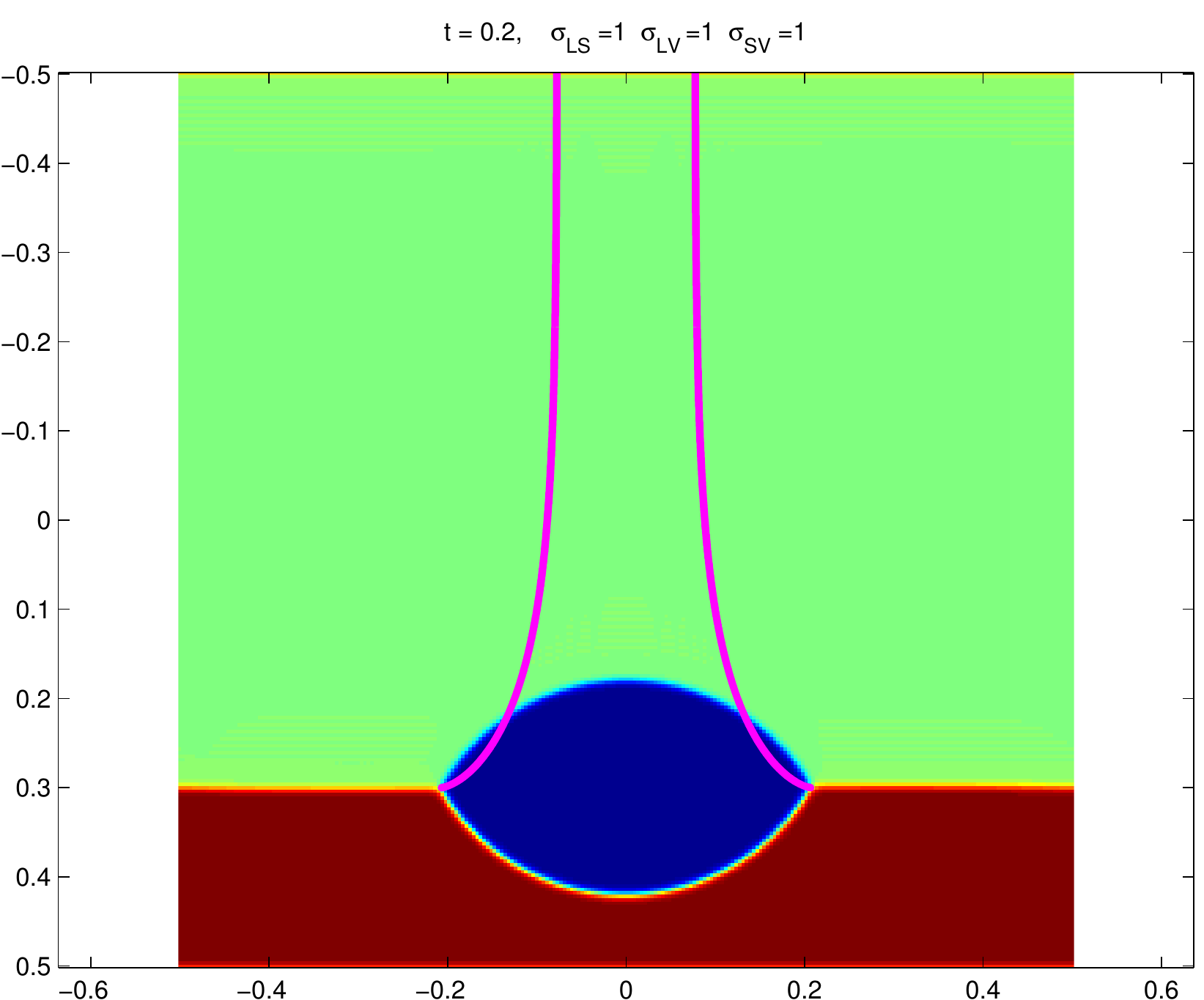}
	 \includegraphics[width=5cm]{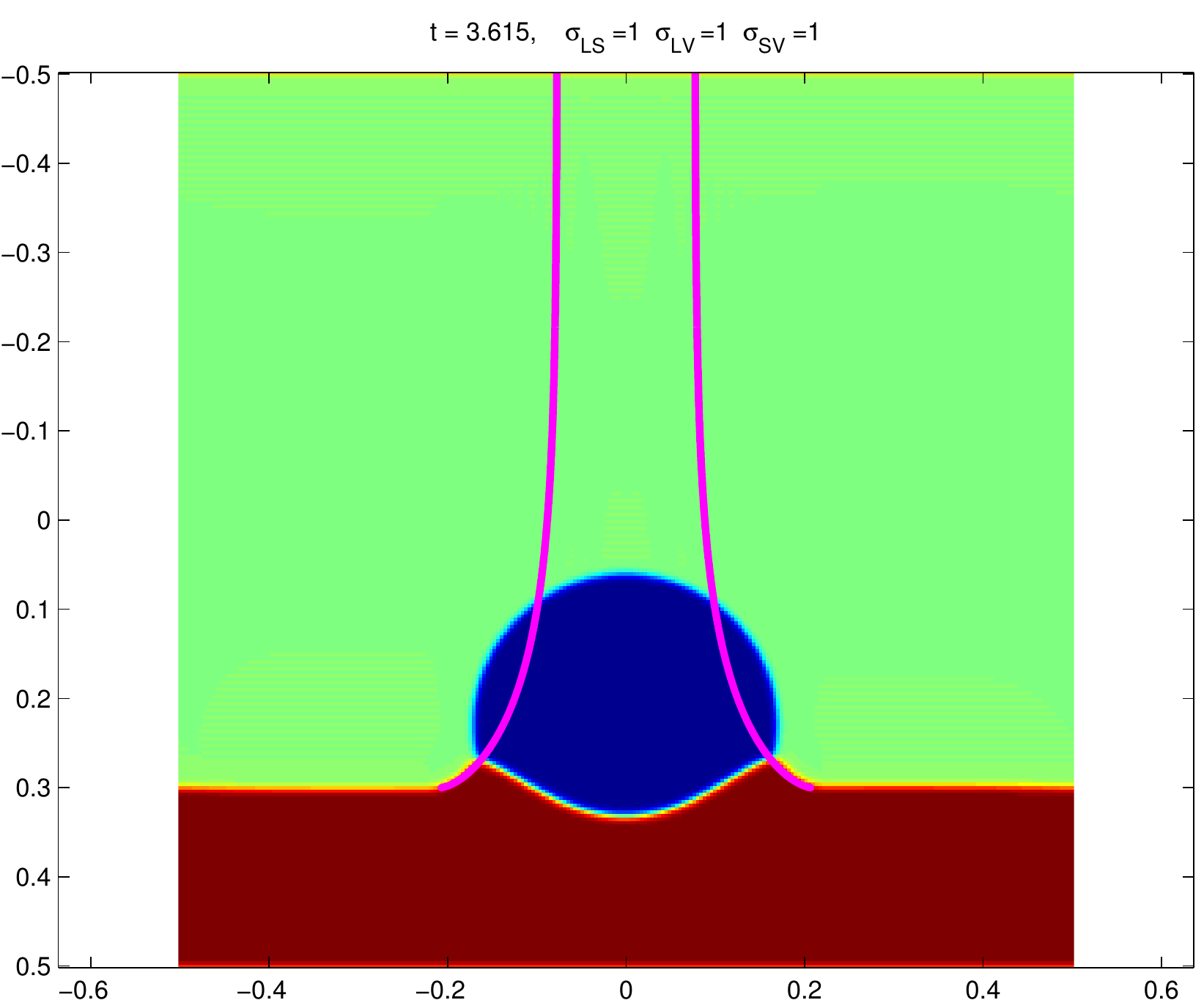} ~\\
	 \includegraphics[width=5cm]{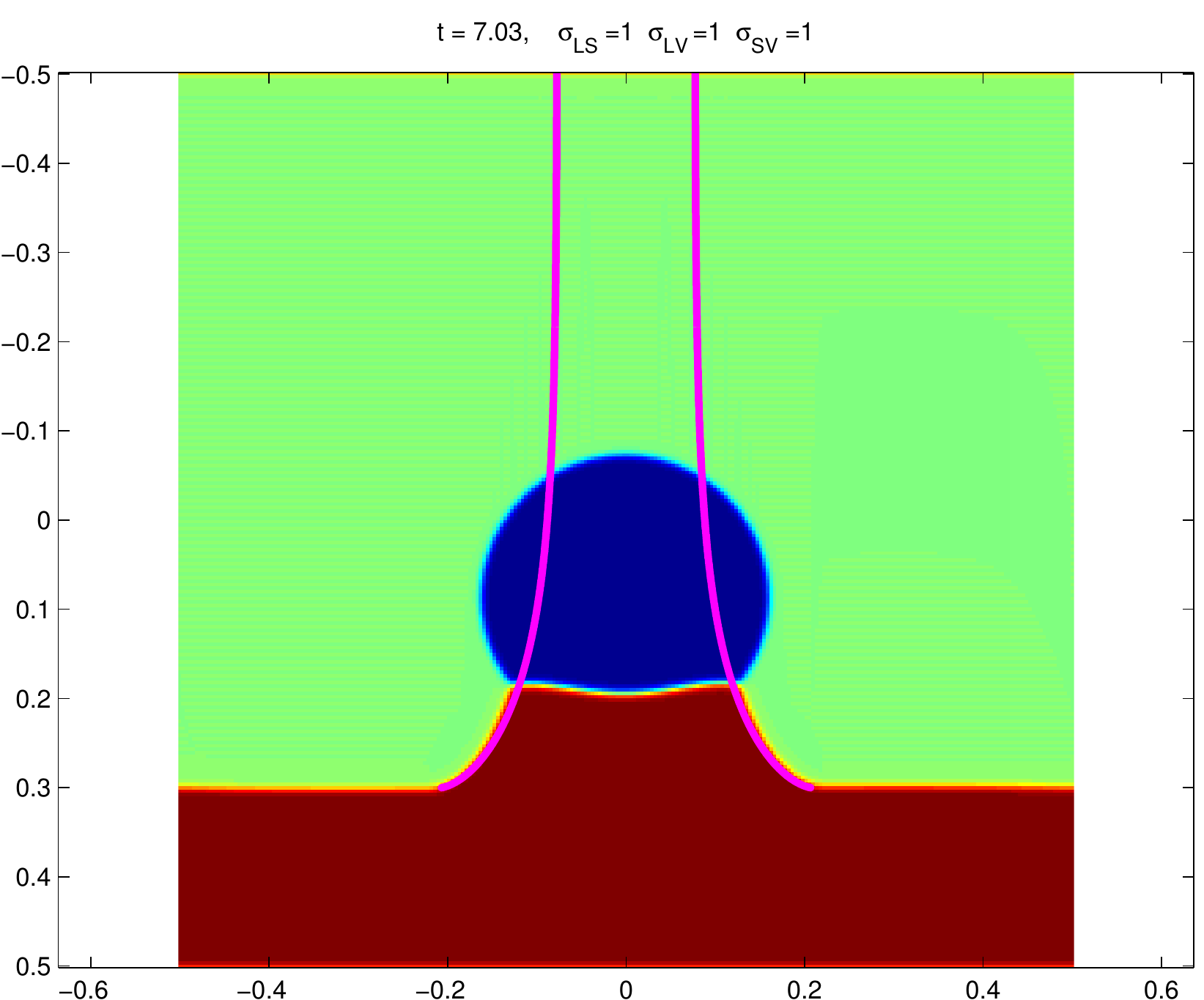}
	 \includegraphics[width=5cm]{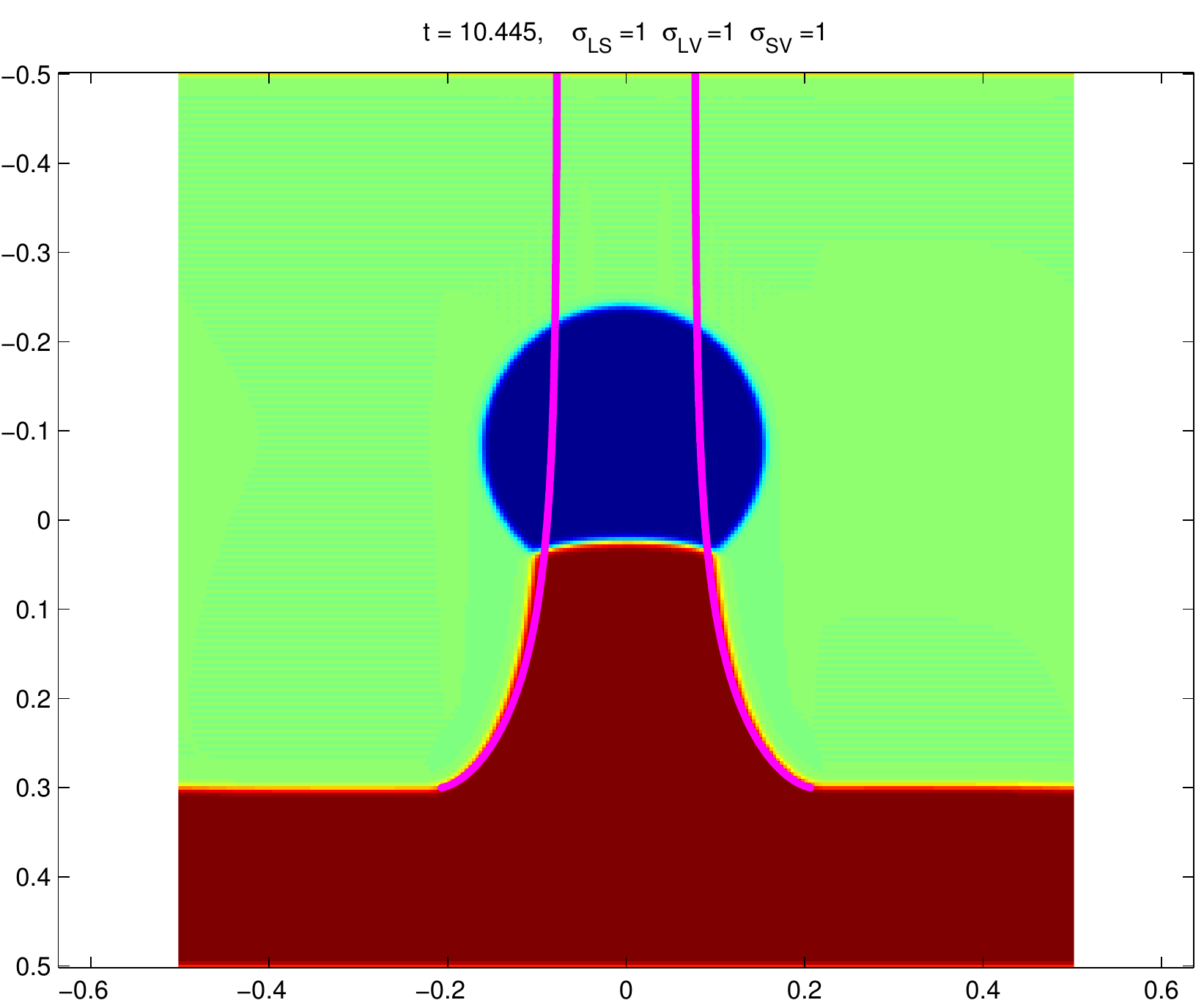}
	 \includegraphics[width=5cm]{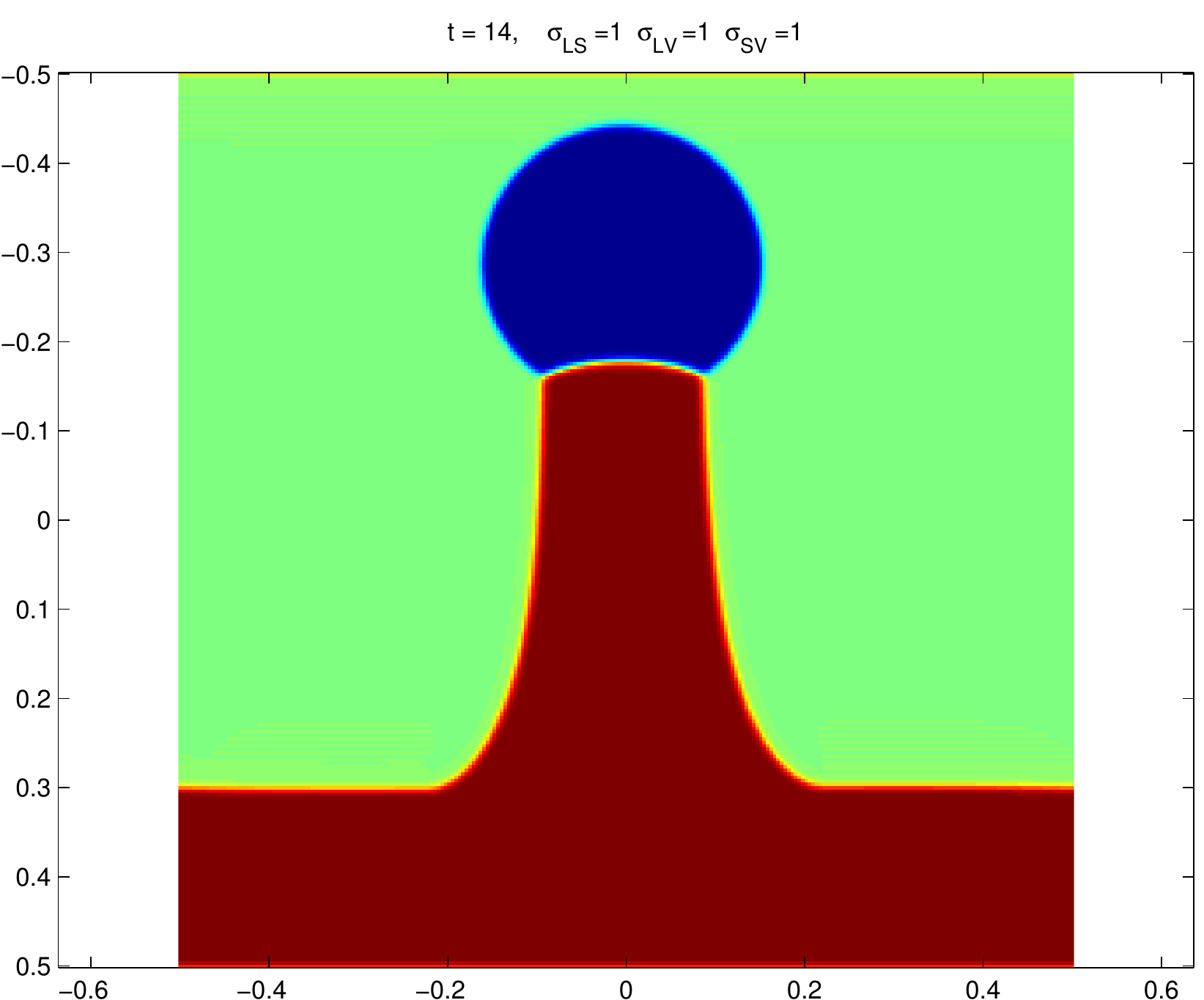}
\caption{Simulation of a $2D$ nanowire evolution using $\sigma_{iso} = (1,1,1)$, captured at different times (from left to right and top to bottom). Starting from the initial configuration, the second image shows the result of the evolution without increase of the solid phase volume and with equal mobilities. Then, the quasi-static model with volume increase and inhomogenous mobilities is used for simulating the other configurations. The purple line represents the theoretically expected profile obtained with the numerical integration of ODE~\eqref{function-h}.}
\label{Nanowire1_isotropic_2D}
\end{figure}

\begin{figure}[htbp]
\centering
	 \includegraphics[width=5cm]{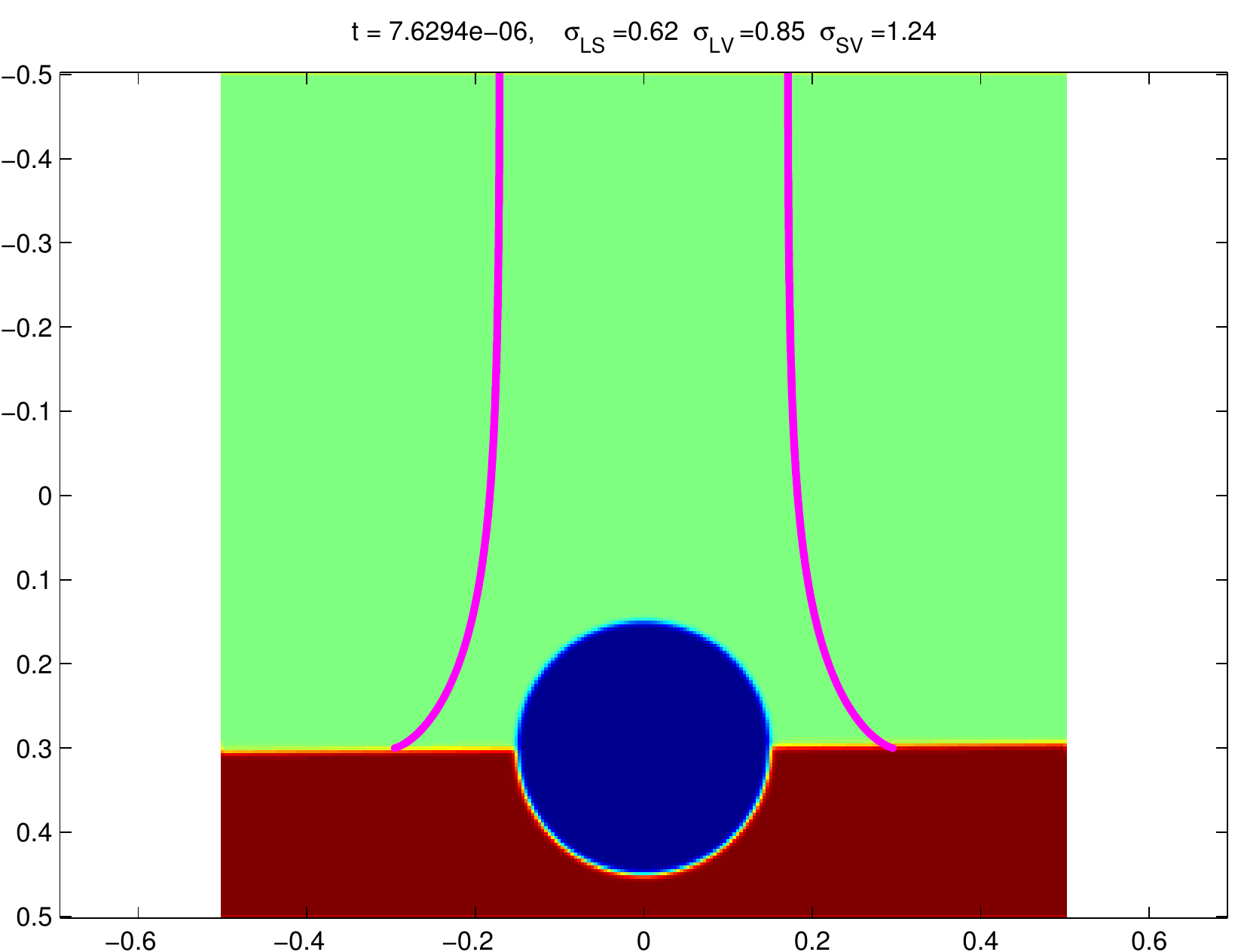}
	 \includegraphics[width=5cm]{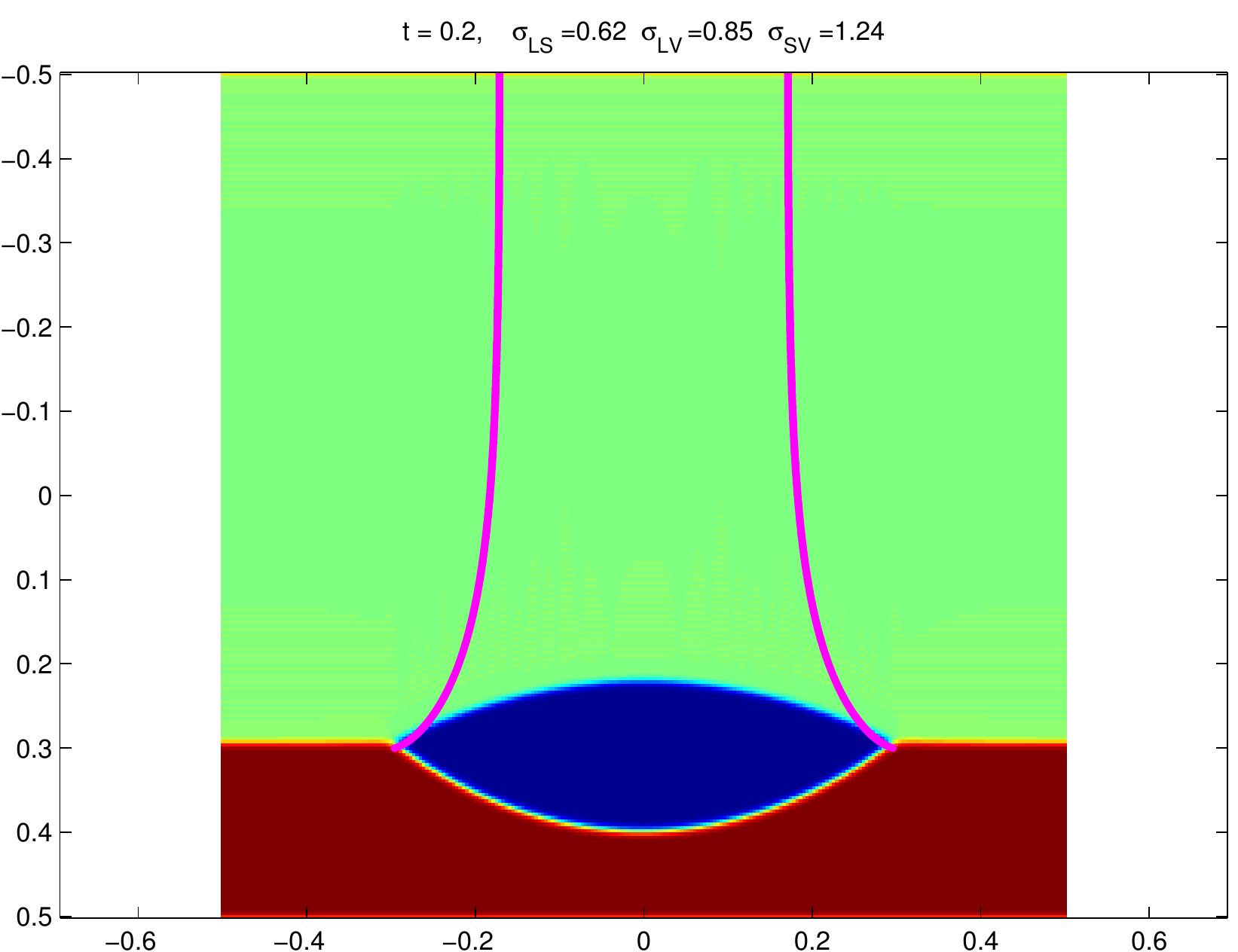}
	 \includegraphics[width=5cm]{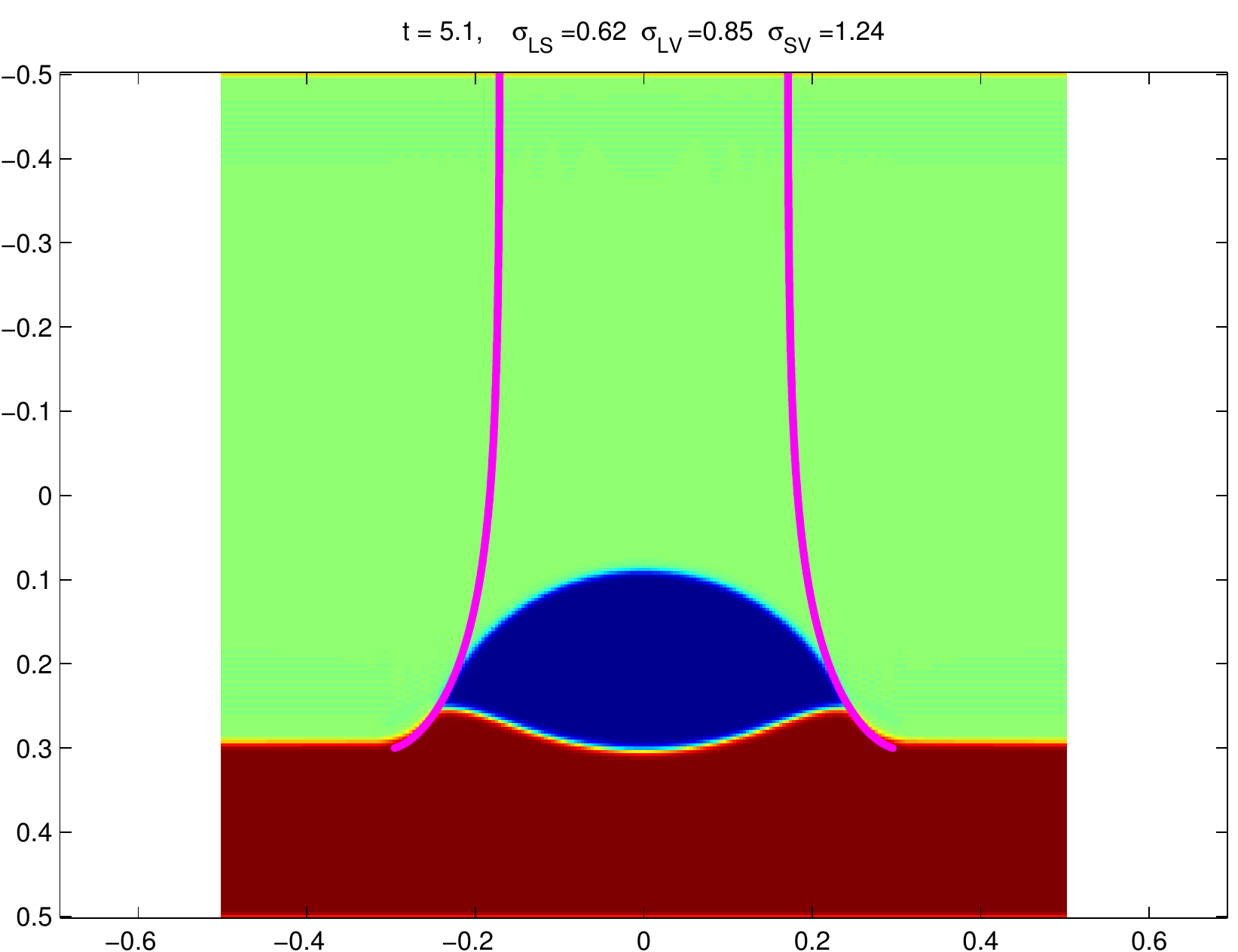} ~\\
	 \includegraphics[width=5cm]{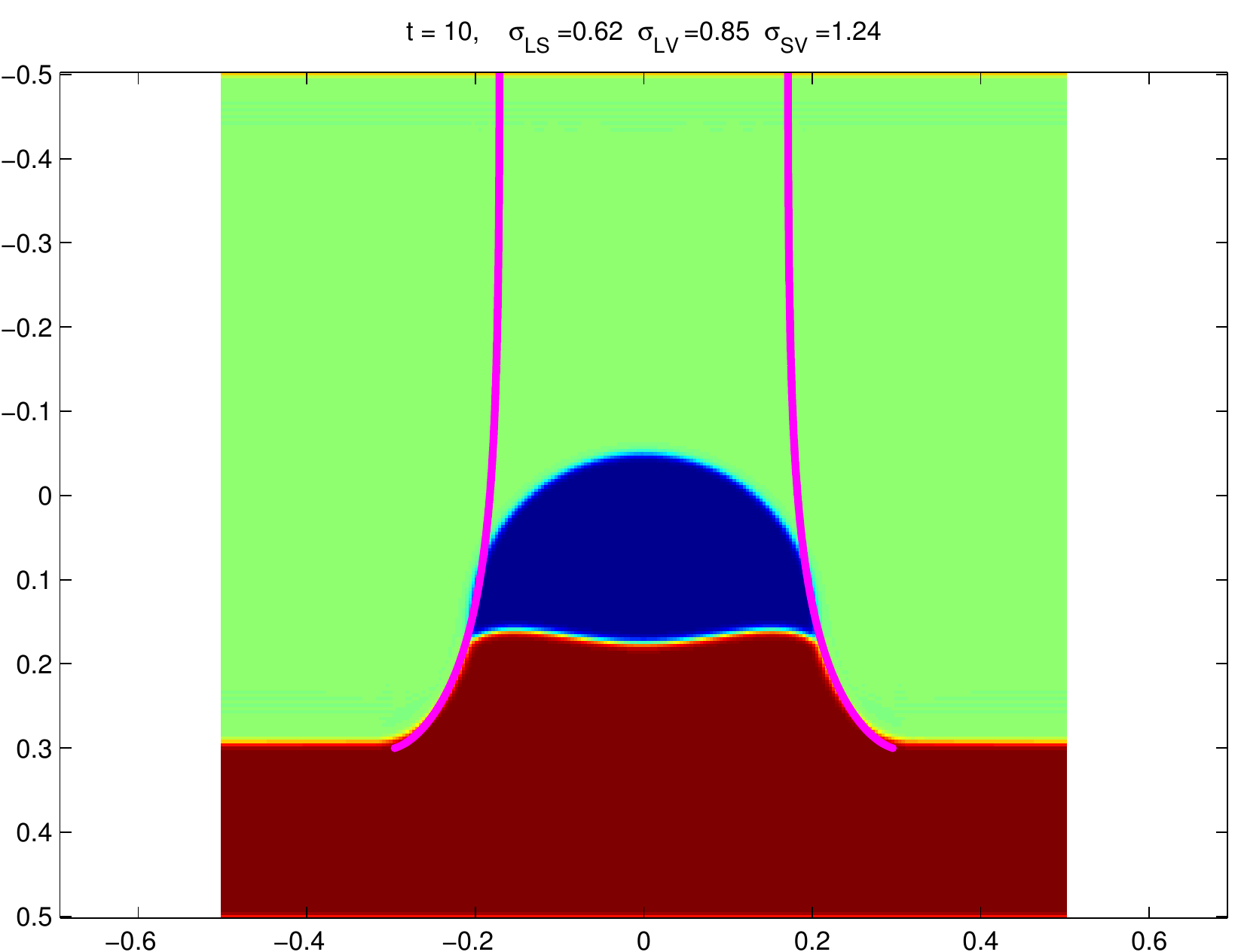}
	 \includegraphics[width=5cm]{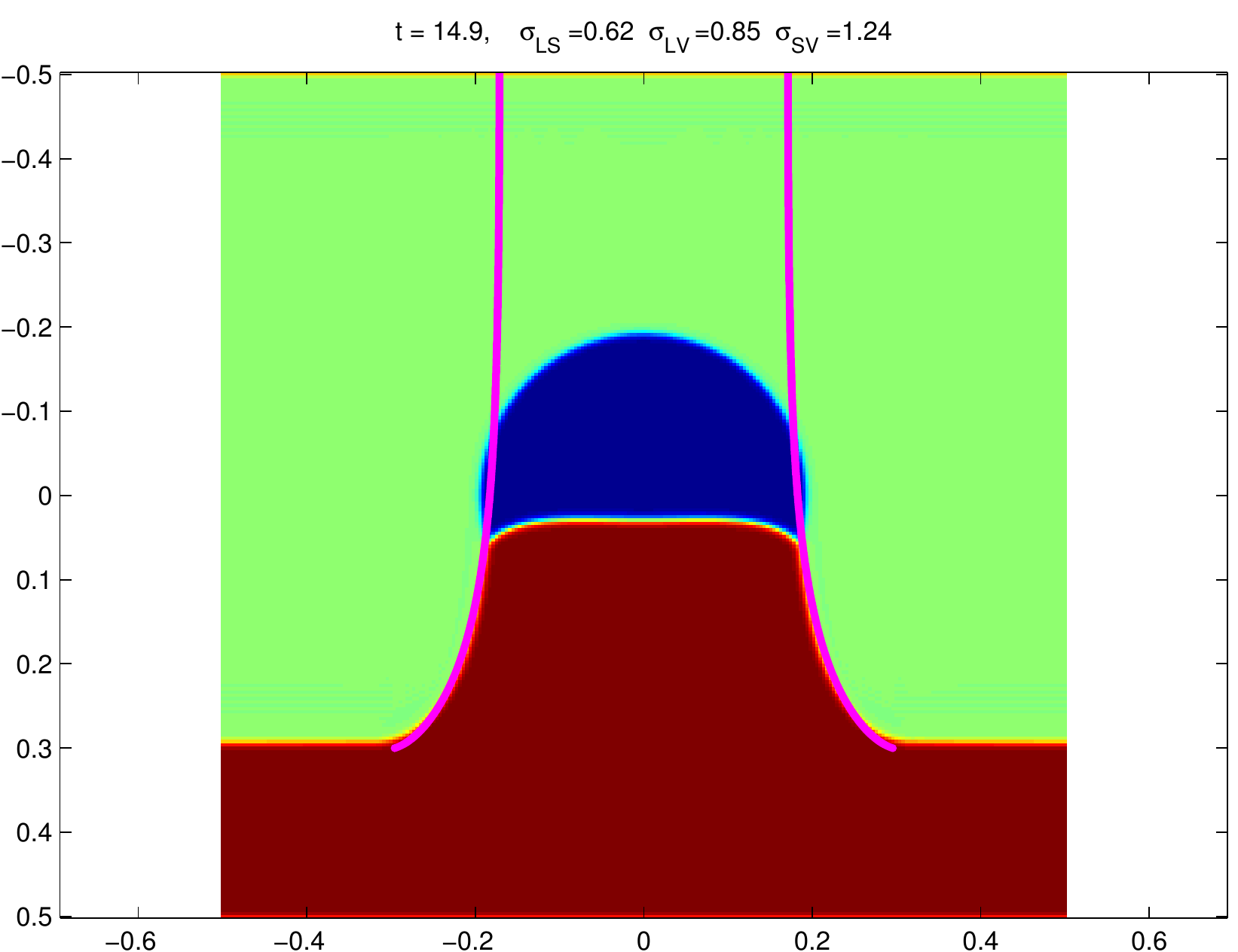}
	 \includegraphics[width=5cm]{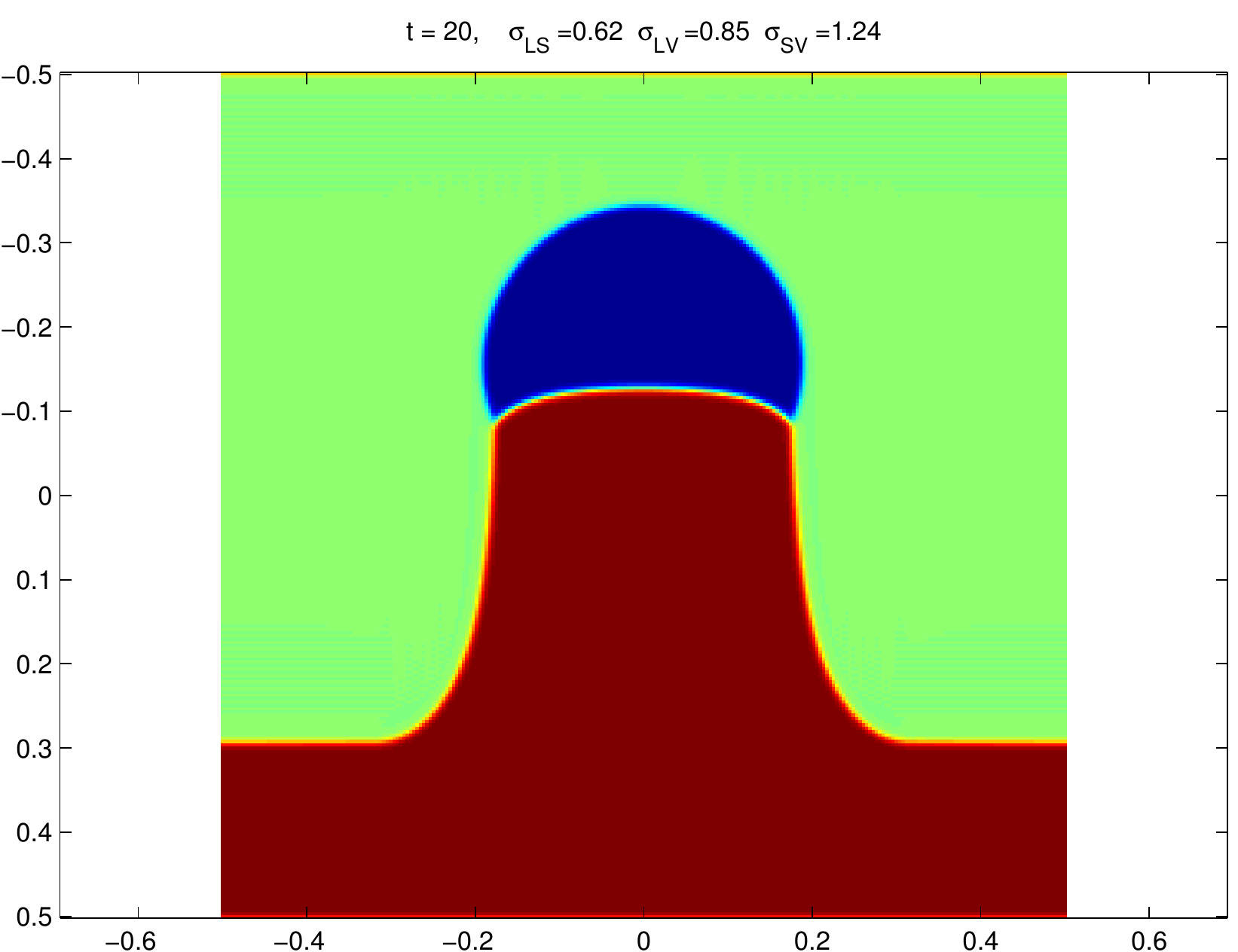}
\caption{Simulation of a $2D$ nanowire evolution using $\sigma_{Au-Si(111)} = (0.62,0.85,1.24)$, captured at different times (from left to right and top to bottom). Starting from the initial configuration, the second image shows the result of the evolution without increase of the solid phase volume and with equal mobilities. Then, the quasi-static model with volume increase and inhomogenous mobilities is used for simulating the other configurations. The purple line represents the theoretically expected profile obtained with the numerical integration of ODE~\eqref{function-h}.}
\label{Nanowire2_isotropic_2D}
\end{figure}

\paragraph{Phase field numerical approximation of nanowire growth} ~\\
Figures \eqref{Nanowire1_isotropic_2D} and \eqref{Nanowire2_isotropic_2D} illustrate some numerical experiments obtained using respectively
$\sigma_{iso} = (1,1,1)$ in the isotropic case and $\sigma_{Au-Si} = (0.62,0.85,1.24)$ for the Au-Si-(111) case. The other parameters are identical with those used previously, {\em i.e.,}
$N = 2^8$, $\epsilon = 1/N$, $\delta_t = 1/N^2$, $L_1, = L_2 = 1$, $c_s = 1/4$. Mobilities are defined as
$${\boldsymbol m} = (m_{LS},m_{VL}, m_{SV}) = \left( \frac{\delta}{1 + \delta},\frac{1}{2}, \frac{\delta}{1 + \delta} \right),$$
with $\delta = 1/(2N)$. As already mentioned, this particular choice is physically sound, and guarantees that mobilities are harmonically additive, see the discussion in Section~\ref{sub:appmodel}. The evolution process is split in two steps :
\begin{itemize}
 \item For $ t \leq T_{growth}=0.2$  the multiphase Cahn-Hilliard energy is minimized with all equal mobility coefficients and without increase of the solid phase, i.e. $c_s = 0$.
 The aim of this first part is to recover the wetting phenomena and to approximate the optimal initial shape of the liquid phase, as in the sharp-theorerical calculations above.
 \item For $t\geq T_{growth}=0.2$ the nanowire begins to grow: the multiphase Cahn-Hilliard energy functional is minimized with inhomogenous mobilities $(m_{LS},m_{VL}, m_{SV}) = \left( \frac{\delta}{1 + \delta},\frac{1}{2}, \frac{\delta}{1 + \delta} \right)$ and with a prescribed increase rate of the solid phase $c_s = 0.25$. In the experimental setting this rate is provided by the external adatom flux, which in turn is fixed by the temperature of the solid source.
\end{itemize}
Figures \ref{Nanowire1_isotropic_2D} and \ref{Nanowire2_isotropic_2D} represent the nanowire shapes at different times, for two different choices of surface tensions, either $\sigma_{iso} = (1,1,1)$ in Figure~\ref{Nanowire1_isotropic_2D} and $\sigma_{Au-Si(111)} = (0.62,0.85,1.24)$ in Figure~\ref{Nanowire2_isotropic_2D}. In each series of images, the first one is the initial shape, and the second image shows the optimal shape at $t = T_{growth}$ (see above). In both experiments, the magenta curve represents the optimal nanowire shape as derived in our sharp-theoretical calculations above. These two numerical simulations illustrate clearly the ability of our numerical approach to approximate in a very realistic way the quasi-static nanowire growth. \\

Lastly, Figure \ref{Nanowire1_isotropic_3D} shows a full 3D numerical simulation obtained using $\sigma_{iso} = (1,1,1)$, $N = 2^8$, $\epsilon = 1/N$, $\delta_t = 1/N^2$, $L_1, = L_2 = 1$, $c_s = 1/4$ and, again,
$${\boldsymbol m} = (m_{LS},m_{VL}, m_{SV}) = \left( \frac{\delta}{1 + \delta},\frac{1}{2}, \frac{\delta}{1 + \delta} \right),$$
with $\delta = 1/(2N)$.

\section{Conclusion}
We showed that multiphase mean curvature flows can be approximated consistently with a phase field method even when highly contrasted, or even degenerate mobilities are involved. The key is to incoporate the mobilities in the metric used for computing the gradient flow of the multiphase perimeter energy. We showed, at least formally, that the diffuse approximation obtained with our model converges to the sharp interface solution, and the convergence has the same order, when mobilities are harmonically additive, as the convergence of the binary phase curvature flow. We also proposed a new quasi-static isotropic approximation of nanowires growth. Numerical simulations of both growing nanowires and droplets wetting on a solid surface confirm the quality of our model. Extension of our approach to the more realistic anisotropic case is ongoing work.

\end{document}